\documentclass[9pt,twoside]{article}
\usepackage[colorlinks]{hyperref}
\usepackage[margin=1in]{geometry}
\usepackage{times,amsmath,stackengine,amsthm,amsfonts,amssymb,mathtools,color,latexsym,amsxtra,layout}
\usepackage{graphicx,subfig}
\usepackage{algorithmicx, algpseudocode}
\usepackage{algorithm}
\usepackage{multirow}

\newcommand{\abs}[1]{\left | #1 \right |}
\newcommand{\para}[1]{\left ( #1 \right )}

\newcommand{\inner}[1]{\left \langle #1  \right \rangle}

\newcommand{\honeinner}[1]{\left \langle #1  \right \rangle_{H^1}}

\newcommand{\ltwoinner}[1]{\left \langle #1  \right \rangle_{L^2}}
\newcommand{\norm}[1]{\left \| #1 \right \|}

\newcommand{\set}[1]{\left \{ #1 \right \}}

\newcommand{\ceil}[1]{\left \lceil #1 \right \rceil}

\newcommand{\eqn}[1]{\begin{equation}\begin{split} #1 \end{split}\end{equation}}
\newcommand{\eqnnolabel}[1]{\begin{equation*}\begin{split} #1 \end{split}\end{equation*}}

\newcommand{\R}{\mathbb{R}}

\newcommand{\p}{\mathcal{P}}

\newcommand{\too}{\longrightarrow}

\theoremstyle{plain}
\newtheorem{theorem}{Theorem}[section]
\newtheorem{lemma}[theorem]{Lemma}
\newtheorem{o-thm}[theorem]{Theorem}

\newtheorem{corollary}[theorem]{Corollary}

\theoremstyle{definition}

\author{Hagop Karakazian \thanks{Hagop Karakazian, Department of Mathematics, American University of Beirut, Lebanon (hk93@aub.edu.lb)} \and Sophie Moufawad \thanks{Sophie Moufawad, Department of Mathematics, American University of Beirut, Lebanon (sm101@aub.edu.lb)}\and Nabil Nassif \thanks{Nabil Nassif, Department of Mathematics, American University of Beirut, Lebanon (nn12@aub.edu.lb)}}

\title{\bf A Finite-Element Model for the Hasegawa-Mima Wave Equation}
\date{\today}

\begin{document}

\maketitle

\begin{abstract}
\noindent In a recent work \cite{kn}, two of the authors have formulated the non-linear space-time Hasegawa-Mima plasma equation as a coupled system of two linear PDEs, a solution of which is a pair $(u,w)$, with $w=(I-\Delta)u$. The first equation is of hyperbolic type and the second of elliptic type. Variational frames for obtaining weak solutions to the initial value Hasegawa-Mima problem with periodic boundary conditions were also derived. Using the Fourier basis in the space variables, existence of solutions were obtained. Implementation of algorithms based on Fourier series leads to systems of dense matrices. \\
In this paper, we use a finite element space-domain approach to semi-discretize the coupled variational Hasegawa-Mima model, obtaining global existence of solutions in $H^2$ on any time interval $[0,T],\,\,\forall T$. \\
In the sequel, full-discretization using an implicit time scheme on the semi-discretized system leads to a nonlinear full space-time discrete system with a nonrestrictive condition on the time step. \\
Tests on a semi-linear version of the implicit nonlinear full-discrete system are conducted for several initial data, assessing the efficiency of our approach.\\\\

\noindent \textbf{Acknowledgements}:  The authors would like to express their gratitude to Prof. Ghassan Antar, AUB Physics Department, for the several discussions held on the physics of the Hasegawa-Mima wave phenomena, providing us suitable testing initial conditions and for his positive comments on the overall numerical results. \\

\noindent \textbf{Keywords}: Hasegawa-Mima; Periodic Sobolev Spaces; Petrov-Galerkin Approximations; Finite-Element Method; Semi-Discrete systems; Implicit Finite-Differences\\

\noindent \textbf{AMS Subject Classification}: 35A01; 35M33; 65H10; 65M60; 78M10

\end{abstract}

\section{Introduction}	
Magnetic plasma confinement is one of the most promising ways in future energy production. The Hasegawa-Mima (HM) model is a simplified two-dimensions turbulent system model which describes the time evolution of drift waves caused during plasma confinement. To understand the phenomena, several mathematical models can be found in literature\cite{iterarticle,hm77,hm78,hw83}, of which the simplest and powerful two dimensions turbulent system model is the HM equation that describes the time evolution of drift waves in a magnetically-confined plasma. It was derived by Akira Hasegawa and Kunioki Mima during late 70s\cite{hm77,hm78}. When normalized, it can\cite{chm89,hariri} be put as the following PDE that is third order in space and first order in time: \begin{equation}\label{hm} 
-\Delta u_t+u_t = \{u,\Delta u\} + \{p,u\}
\end{equation}
\noindent where $\{u,v\}=u_xv_y-u_yv_x$ is the Poisson bracket, $u(x,y,t)$ describes the electrostatic potential, $p= \ln \dfrac{n_0}{\omega_{ci}}$ is a function depending on the background particle density $n_0$ and the ion cyclotron frequency $\omega_{ci}$, which in turn depends on the initial magnetic field. In this context, $p=0$ refers to homogeneous plasma, and $p \neq 0$ refers to non-homogeneous plasma. As a cultural note, equation (\ref{hm}) is also referred as the Charney-Hasegawa-Mima equation in geophysical context that models the time-evolution of Rossby waves in the atmosphere\cite{chm89}.\\
In this paper, we deal with the numerical solution to Hasegawa-Mima equation on a rectangular domain with the solution $u$, satisfying periodic boundary conditions (PBCs). For that purpose, we consider $\Omega = (0,L) \times (0,L) \subset \R^2$ and use the frame of periodic Sobolev spaces which are closed subspace of $H^m(\Omega)$, and therefore itself  a Hilbert space. Specifically:
\begin{equation}
\begin{cases}
H_P^0(\Omega)=L^2(\Omega) \ \mbox{ and } \ H_P^\infty(\Omega):= \cap_{m \geq 1}{\{H_P^m\}},\\
H_P^1(\Omega) =\{u \in H^1(\Omega)\,|\,u(x,0)=u(x,L), x \in (0,L)\mbox{ a.e. },\,u(0,y)=u(L,y),\, y \in (0,L) \mbox{ a.e. }\},\\
H_P^2(\Omega)=\set{u \in H^2(\Omega)  |  u,u_x,u_y  \in H_P^1(\Omega)}\\
H_P^3(\Omega)=\set{u \in H^3(\Omega) | u, u_x, u_y, u_{xx}, u_{yy}, u_{xy}, u_{yx} \in H_P^1(\Omega)}
\end{cases}
\end{equation}
In addition, we use for $p > 2$, the periodic Banach-Sobolev Spaces:\\
$$W_P^{1,p}(\Omega) = \{u \in W^{1,p}(\Omega)\,|\,u(x,0)=u(x,L) \, \mbox{ a.e. } x \in (0,L),\,u(0,y)=u(L,y) \, \mbox{ a.e. }y \in (0,L) \},$$


\noindent Given an initial data $u_0:\overline{\Omega} \to \mathbb{R}$, we seek $u: \overline{\Omega}\times [0,T]  \to \mathbb{R}$ such that:
\begin{equation}\label{hmpbc}
\left\{\begin{array}{lll}
-\Delta u_t+u_t = \{u,\Delta u\} + p_xu_y - p_yu_x      &   \mbox{on }  \Omega\times (0,T] & (1) \\
\mbox{PBCs on } u, u_x \mbox{, and } u_y &  \mbox{on }   \partial \Omega\times(0,T] & (2)\\
     u(x,y,0)=u_0(x,y)  &   \mbox{on } \Omega& (3)
\end{array}\right.
\end{equation}
In \cite{kn}, without loss of generality for the proof of existence, we assume that the background particle density $n_0$ is a function of $x$ only, such that $p_x = \hat{k}$ is a constant and $p_y=0$, i.e. $n_0 = e^{Ax+B}$ for $A,B \in \mathbb{R}$.
When dealing with (\ref{hmpbc}.1), the major difficulty to circumvent, both theoretically and computationally, is the Poisson bracket $\{u,\Delta u\}$. To overcome this issue, we have formulated it in \cite{kn}, as a coupled system of linear hyperbolic-elliptic PDEs that will be naturally amenable to provide a Finite Element scheme for obtaining a numerical approximation/simulation.
For this purpose, a new variable $w=-\Delta u + u$ is introduced, leading to the identity 
$$\{u,\Delta u\} =\{u,u-w\}=\{u,u\}+\{u,-w\}=-\{u,w\}=\{w,u\}=w_xu_y-w_yu_x = { - }\vec{V}(u) \cdot \nabla w $$
where $\vec{V}(u)= -u_y \vec{\textbf{i}} + u_x \vec{\textbf{j}}$ is a {\it divergence-free vector field} ($\mbox{div}(\vec{V}(u))=0$). Then system \eqref{hmpbc} with $p_x = \hat{k}$ and $p_y=0$,  becomes equivalent to the coupled hyperbolic-elliptic PDE system, 
\begin{equation}\label{HMC2}
\left\{\begin{array}{lll}
w_t + \vec{V}(u) \cdot \nabla w = \hat{k}u_y & \mbox{on }  \Omega \times  (0,T] &(1)\\
-\Delta u+u=w &\mbox{on }  \Omega \times  (0,T] &(2)\\
\mbox{PBC's on } u,\, u_x,\, u_y,\,w & \mbox{on }  \partial \Omega\times [0,T] & (3)\\
u(0)=u_0 \mbox{ and }w(0)=w_0 & \mbox{ on } \overline{\Omega}.  & (4)\\
\end{array}\right.
\end{equation} 
Using equation \eqref{skeww} which is obtained through Green's formula and the imposition of periodic boundary conditions,
\begin{equation} \label{skeww} \ltwoinner{\vec{V}(u) \cdot \nabla v,w} = -  \ltwoinner{\vec{V}(u) \cdot \nabla w,v}
\end{equation}
system \eqref{HMC2} can be put in the following strong semi-variational form (on the space variables) whereby one seeks a pair $\{u,w\} \in C([0,T], H^2(\Omega) \cap H_P^1(\Omega)) \times \left [C([0,T], L^2(\Omega)) \cap C^1((0,T), L^2(\Omega))\right ]$ such that


\begin{equation}\label{HMC-V1}
\left\{\begin{array}{ll}
\ltwoinner{w_t,v} = \ltwoinner{\vec{V}(u) \cdot \nabla v,w} + \ltwoinner{\hat{k}u_y,v},  &(1)\\
\honeinner{u,v}=\ltwoinner{w,v}, &(2)\\
 \forall v \in  W^{1,\infty}_P(\Omega)\cap H_P^1(\Omega), \quad \forall t \in (0,T]  & \\
\end{array}\right.
\end{equation} 
with $u(0)=u_0\in H_P^2(\Omega),\,w(0)=w_0=u_0-\Delta u_0\in L^2(\Omega)  $.\vspace{2mm}

\noindent A similar formulation to \eqref{HMC-V1} has been handled in \cite{kn} where using Fourier series, we prove the existence of a solution  $$\{u,w\} \in C([0,T],H_P^1(\Omega)) \cap L^2(0,T; H^2) \times \left [C([0,T], L^2(\Omega)) \cap C^1((0,T), L^2(\Omega))\right ].$$ 
\noindent Note that the restriction $v \in  W^{1,p}_P(\Omega)$ ensures that $\vec{V}(u) \cdot \nabla v \in L^2(\Omega)$, which justifies this formulation. This allows reducing the regularity of the initial conditions $(u_0,w_0)$ from $H^3_P(\Omega) \times H^1_P(\Omega)$ in \cite{kn} to $H^2_P(\Omega) \times L^2(\Omega)$. \\
The formulation used in this paper intends to avoid dealing with the time-derivative $w_t$. For that purpose, we derive  \eqref{HMC-V1} a time integral semi-variational formulation.
\subsection*{Time Integral Formulation}
By integrating (\ref{HMC-V1}.1) over the temporal interval $[t,t+\tau]$, with $0\leq t\leq T-\tau$, one reaches the following L\textsuperscript{2} Integral Formulation:
 \begin{equation}\label{HMC2-tau}
\left\{\begin{array}{ll}
u\in L^2(0,T;H^2(\Omega)\cap H^1_P(\Omega)), \quad w \in L^2(0,T;L^2(\Omega)) & \\
\ltwoinner{w(t+\tau)-w(t),v} =\int_t^{t+\tau} \ltwoinner{\vec{V}(u(s)) \cdot \nabla v,w(s)}+ \ltwoinner{\hat{k}u_y(s),v} ds&(1)\\
\honeinner{u(s),v}=\ltwoinner{w(s),v}, \ &(2)\\
\forall v \in W^{1,\infty}(\Omega) \cap H^1_P(\Omega), \quad \forall t,\,\tau,\, 0\le t<t+\tau\le T, \quad \forall s \in [t,t+\tau]  &\\ 
\end{array}\right.\vspace{-2mm}
\end{equation}
with $u(0)=u_0\in H^2\cap H^1_P,$ and $w(0)=w_0=u_0-\Delta u_0 $.\\
Such formulation is well-suited for semi and full discretization of the original system \eqref{hmpbc} and \eqref{HMC2}. \\\\
For semi-discretization, we start defining the finite element spaces as follows.

\subsection*{Finite-Element Space Semi-Discretization}
\begin{figure}[h]
\centering
\includegraphics[scale=0.15]{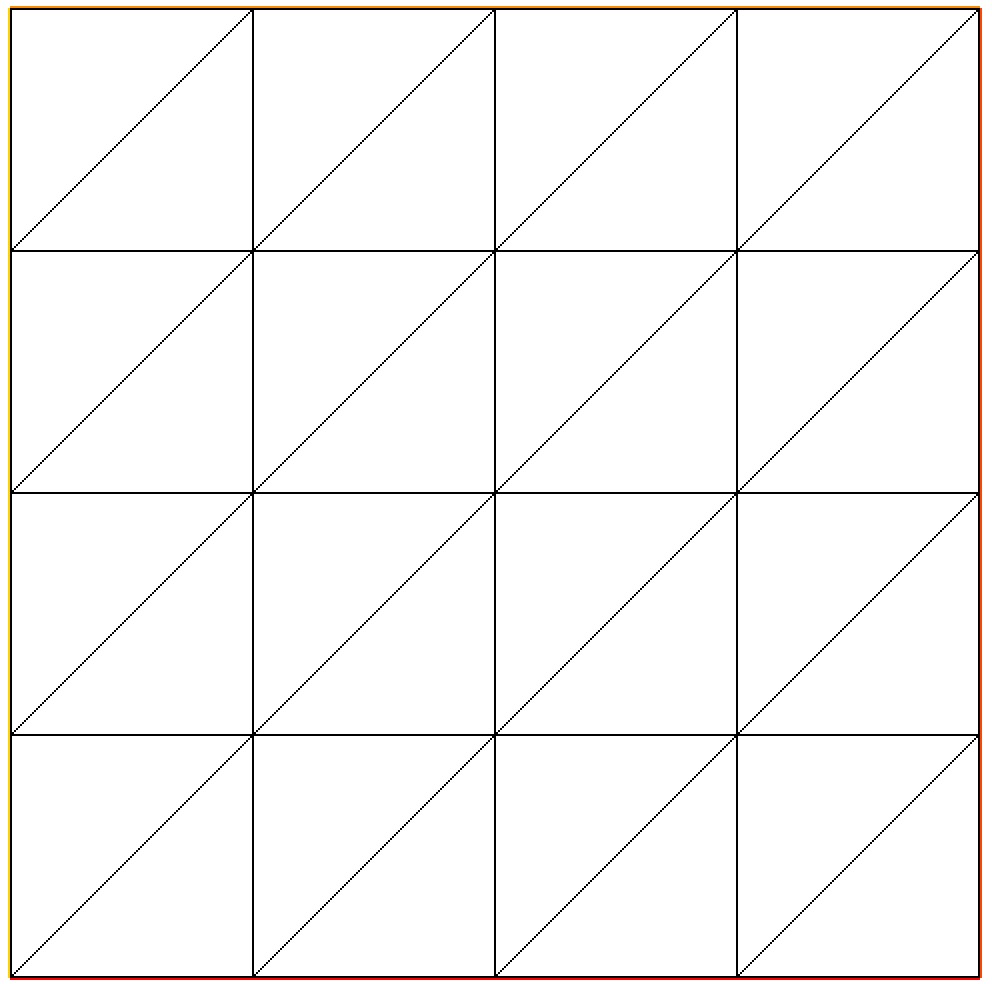}
\caption{\it\small A two-dimension meshing of $\Omega$}\label{fig:mesh1}
\end{figure}

Systems (\ref{HMC-V1}) and \eqref{HMC2-tau} lead to equivalent $\mathbf{\mathbb{P}_1}$ { Finite-element space semi-discretization} constructed as follows:\\
Let 
$\p_x=\{x_i|i=1,...,n\}$ be a partition of $(0,L)$: $0=x_1<x_2<...<x_n=L$ in the $x$ direction and similarly in the $y$ direction, $\p_y=\{y_j|j=1,...,n\}$. Let now:
$$\mathcal{N}=\{P_I(x_i,y_j)|I=1,2,...,N=n^2\}=\p_x\times\p_y,$$
be a structured set of nodes covering $\overline{\Omega}$. Based on $\mathcal{N}$, and as indicated in Figure \ref{fig:mesh1}, 
one obtains a conforming (Delaunay) structured triangulation $\mathcal{T}$ of $\overline{\Omega}$, i.e., $\mathcal{T}=\{E_J|J=1,2,...,M\},\,\,\overline{\Omega}=\cup_{J}{E_J}$. The $\mathbb{P}_1$ finite element subspace $X_N$ of $H^1(\Omega)$ is given by:
 $$X_N=\{v\in C(\overline{\Omega})|v\mbox{ restricted to }E_J\in\mathbb{P}_1,\,J=1,2..,M\} \subset W_P^{1,p}, \quad  1 \leq p \leq \infty$$ 
  with ${\cup_{N\geq 1}\{X_N\}}$  approximating functions in $H^1(\Omega)$.
 For that purpose, we let ${B}_N=\{\varphi_I|\,I=1,2,...N\}$ be a finite element basis of functions with compact support in $\Omega$, i.e.,: 
 $$\forall v_N \in X_N:\,v_N(x,y)=\sum_{I=1}^N{V_{I}\varphi_I(x,y)},\, V_I=v_N(x_I,y_I).$$
We state now useful estimates used in this paper. 
\subsection*{Approximation properties of $X_N$ in $H^1(\Omega)$}
These can be found in section 3.1 of Ciarlet \cite{ciarlet}, specifically:
\eqn{\label{definterpolant} \forall v \in W^{1,p} \mbox{ we define } \pi_N(v):=\sum_{I=1}^N{V_{I}\varphi_I(x,y)}\in X_N \mbox{ to be the interpolant of } v \mbox{ in } X_N}
One has the following estimates:
 \begin{equation}\label{density1}
\forall v \in W^{1,p},\quad  |v-\pi_N(v)|_{0,p}\le  C \; \frac{1}{n}\; |v|_{1,p}, \quad p \in (1,\infty]
 \end{equation}
 and
 \begin{equation}\label{density12}
\forall v \in W^{2,\infty},\quad  |v-\pi_N(v)|_{1,\infty}\le  C \; \frac{1}{n}\; |v|_{2,\infty},
 \end{equation}
 

\noindent Since the Hasegawa-Mima equation is set in $H^1_P(\Omega)$, we let  $$X_{N,P}=X_N\cap H^1_P(\Omega).$$
To discretize \eqref{HMC2-tau},
we start with $u_N(0) = \pi_N(u_0)$, $w_N(0) = \pi_N(w_0)$, and $w_0=u_0-\Delta u_0$, then given \linebreak $(u_N(t),w_N(t)) \in X_{N,P}\times X_{N,P}$ where 
$u_N(t)=\pi_N(u) = \sum_{J=1}^N{U_J(t)\varphi_J(x,y)}$, $U_J(t)=u_N(x_J,y_J,t)$,\linebreak $w_N(t)=\pi_N(w) = \sum_{I=1}^{N}{W_I(t)\varphi_I(x,y)}$, and $W_I(t)=w_N(x_I,y_I,t)$,  one seeks:
\begin{equation}\label{HMC-Comp-Mod}
\left\{\begin{array}{ll}
u_N(t+\tau)\in  X_{N,P}, \quad w_N(t+\tau) \in \times X_{N,P}& \\
\ltwoinner{w_N(t+\tau) - w_N(t),v} = \int_t^{t+\tau}\ltwoinner{\vec{V}(u_N(s)) \cdot \nabla v,w_N(s)} +\ltwoinner{\hat{k}u_{N,y}(s),v}ds&(1)\\
\honeinner{u_N(s),v}=\ltwoinner{w_N(s),v},  &(2)\\
\forall v\in  X_{N,P} \quad \forall t,\,\tau,\, 0\le t<t+\tau\le T, \quad \forall s\in\{t,t+\tau\}
\end{array}\right.
\end{equation} 
We now rewrite equation (\ref{HMC-Comp-Mod}.1) by dividing it by $\tau$ giving
$$\dfrac{1}{\tau}\ltwoinner{w_N(t+\tau) - w_N(t),v} = \dfrac{1}{\tau}\int_t^{t+\tau}\ltwoinner{\vec{V}(u_N(s)) \cdot \nabla v,w_N(s)} +\ltwoinner{\hat{k}u_{N,y}(s),v}ds$$
Letting $\tau$ tend  to $0$, implies that every solution $u_N(t)=\pi_N(u)$, and $w_N(t)=\pi_N(w)$ of \eqref{HMC-Comp-Mod}
 is a solution to the semi-discretization of the $H^1$ formulation \eqref{HMC-V1} given by
 \begin{equation}\label{HMC-Disc-V1} 
\left\{\begin{array}{ll}
\ltwoinner{w_{N,t},v} - \ltwoinner{\vec{V}(u_N) \cdot \nabla v,w_N}= \ltwoinner{\hat{k}u_{N,y},v}, \, \forall v \in X_{N,P}, \ \forall t \in (0,T] &(1)\\
\honeinner{u_N(t),v}=\ltwoinner{w_N(t),v}, &(2)\\
\forall v \in X_{N,P}, \quad \forall t \in [0,T] &\\
\end{array}\right.
\end{equation}
with $u_N(0)=\pi_N (u_0)$, $w_N(0)=\pi_N (w_0)$, $w_0=u_0-\Delta u_0$. 

\noindent Defining the associated vectors  $W(t) = \{W_I(t) \}_I $ and  $U(t) = \{U_J(t) \}_J$,  system \eqref{HMC-Disc-V1} is equivalent in vector form to \begin{equation}\label{HMC-Disc-V1-1}
\left\{\begin{array}{ll}
M\frac{dW}{dt} + S(U)W= RU,\, \qquad  \forall t \in (0,T] &(1)\\
KU(t)=MW(t) \qquad \qquad \forall t \in (0,T] &(2)\\
U(0)=U_0;\,W(0)=W_0 & (3)\\
\end{array}\right.
\end{equation}
with $M$, $K$, $S(U)$ and $R$, $N\times N$ matrices, defined as follows:
\begin{itemize}
\item $M=\{\ltwoinner{\varphi_I,\varphi_J}|\,1\le I,J\le N\}$, \quad $K=\{\honeinner{\varphi_I,\varphi_J}|\,1\le I,J\le N\}$, \quad $R=\left\{\ltwoinner{\hat{k}\varphi_{I,y},\varphi_J}|\,1\le I,J\le N\right\}$\item $S(U)=\left\{ -\ltwoinner{\vec{V}(u_N) \cdot \nabla\varphi_J,\varphi_I}|\,1\le I,J\le N\right\} = \left\{ \ltwoinner{\vec{V}(u_N) \cdot \nabla\varphi_I,\varphi_J}|\,1\le I,J\le N\right\}$,\\  where $u_N(t)=\sum_{K=1}^N{U_K(t)\varphi_K(x,y)}$ and
\begin{equation}\label{Skew-S}
   ( \vec{V}(u_n) \cdot \nabla\varphi_I) \; \varphi_J = - \dfrac{\partial \varphi_I }{\partial x} \varphi_J \sum_{K=1}^N \para{U_K(t)\dfrac{\partial\varphi_K}{\partial y}} + \dfrac{\partial \varphi_I }{\partial y} \varphi_J \sum_{K=1}^N \para{U_K(t)\dfrac{\partial \varphi_K}{\partial x}}
\end{equation}
\end{itemize}
 When implementing system \eqref{HMC-Disc-V1-1} one takes periodicity into account, reducing the degrees of freedom from $N = n^2$ to $N_1 = (n-1)^2$.

\subsection*{Statement of Results}
Using a compactness technique, we prove in Section \ref{Sec-Discretization}, the existence  of a limit point $(u,w)$ to the pair $(u_N,w_N)$ and accordingly, the  existence of a solution to the Hasegawa-Mima coupled system which states as follows.
\begin{theorem}\label{main-theorem} Let $u_0\in H^2(\Omega)\cap H^1_P(\Omega)$ and $w_0=u_0-\Delta u_0\in L^2(\Omega)$. Then for all $T>0$, there exists a unique solution $(u_N(t),w_N(t))$ to (\ref{HMC-Disc-V1}). Furthermore the sequence $\{(u_N(t),w_N(t)),\,N>0\}$ admits a subsequence that converges to a solution pair $(u(t),\,w(t))\in L^2(0,T;H^2(\Omega)\cap H^1_P(\Omega))\times L^2(0,T;L^2(\Omega))$, such that:
\begin{equation}\label{HMC-V2-sub}
\left\{\begin{array}{ll}
\ltwoinner{w(t)-w_0,v} =\int_0^t \ltwoinner{\vec{V}(u(s)) \cdot \nabla v,w(s)}+ \hat{k}\ltwoinner{u_y(s),v} ds& (1)\\
\honeinner{u(t),v}=\ltwoinner{w(t),v}, &(2)\\
\forall v \in H^1_P(\Omega)\cap W^{1,\infty}(\Omega), \quad \forall t \in [0,T]
\end{array}\right.
\end{equation}
\end{theorem}
\noindent In Section \ref{sec:Full-Discrete} we introduce the fully implicit nonlinear discrete scheme \eqref{HMC-Comp}, and prove the existence and uniqueness of its solution under a restriction on the time step.
 In Section \ref{Sec:Algorithm}, we complete the discretization cycle by presenting  the resulting algorithm and its implementation using FreeFem++ software, generating in the sequel the system matrices ${M}, {K}, S(U), R$. The obtained numerical results indicate the robustness of this new software, particularly in handling the complex discretization of the Poisson bracket and circumventing the difficulties encountered in \cite{hariri}.
Concluding remarks are provided in Section \ref{sec:Conclude}.
\section{Proof of Theorem \ref{main-theorem}}\label{Sec-Discretization}
At the core of the proof of this result, are
\begin{enumerate}
\item Existence and uniqueness of solutions to \eqref{HMC-Disc-V1} and equivalently to (\ref{HMC-Disc-V1-1}) proven in Section \ref{Sub-Semi-Discrete} which uses a standard existence theorem for systems of ODE's of the form $\vec Y'(t) = \vec{F}(\vec{Y})$, $\vec{F}$ Lipschitzian.
\item A-priori estimates on these solutions shown in Section \ref{sec:app}.
\item In Section \ref{sec:2.3}, a compactness argument would allow passing to the limit for test functions $v \in W^{2,\infty}(\Omega) \cap H^1_P(\Omega)$ in the formulation \eqref{HMC-V2-sub} of Theorem \ref{main-theorem}.
Finally, a density argument of $W^{2,\infty}(\Omega) \cap H^1_P(\Omega)$ in $H^1_P(\Omega)$ allows us to complete the proof.
 
 \end{enumerate}
This thread of items to be proven requires skew-symmetry results obtained in the following section.


\subsection*{Preliminary Result: Skew-Symmetry on $X_{N,P}$}
For that purpose, we start by obtaining a skew-symmetry result, stated in the following proposition.
\begin{theorem}\label{skewxnP} 
For all $\{v,z,\phi\}\in X_{N,P}\times X_{N,P}\times X_{N,P}$, one has:
\begin{equation}
\ltwoinner{\vec{V}(v) \cdot \nabla z,\phi}= - \ltwoinner{\vec{V}(v) \cdot \nabla \phi,z}.
\end{equation}
\end{theorem}
\noindent To prove Theorem \eqref{skewxnP}, we start by breaking $ \ltwoinner{\vec{V}(v) \cdot \nabla z,\phi}$  into a sum of integrals over each triangle \linebreak $E_J\in\mathcal{T},\,J=1,...,M$. Specifically:
\begin{equation}\label{breakskew}
\ltwoinner{\vec{V}(v) \cdot \nabla z,\phi}=\sum_{E_J\in\mathcal{T}}{\left \langle \vec{V}(v) \cdot \nabla z,\phi \right \rangle_{L^2(E_J)}}
\end{equation}
\begin{lemma}\label{lem-est-0}\noindent Let $E_J\in \mathcal{T}$ with vertices $1_J,\,2_J,\,3_J$. Let also  $\nu$ being the unit outer normal to $\partial E_J$, defined piecewise on each of the sides of the triangle $E_J$ and denoted respectively by  $\nu_1$, $\nu_2$, $\nu_3$ on $[1_J,2_J]$, $[2_J,3_J]$ and $[3_J,1_J]$. Then  one has:
$$\left \langle \vec{V}(v) \cdot \nabla z,\phi \right \rangle_{L^2(E_J)}= - \left \langle \vec{V}(v) \cdot \nabla \phi, z \right \rangle_{L^2(E_J)} + \gamma_{1,J}\int_{1_J}^{2_J}{z\phi}+  \gamma_{2,J}\int_{2_J}^{3_J}{z\phi}+ \gamma_{3,J}\int_{3_J}^{1_J}{z\phi},$$
with:
\begin{equation}\label{gammas}
\left\{\begin{array}{c}
\gamma_{1,J}=[(v(2_J)-v(1_J))(\frac{\nu_{1,x}}{y_{2_J}-y_{2_J}}-\frac{\nu_{1,y}}{x_{2_J}-x_{1_J}})]\\
\\
\gamma_{2,J}=[(v(3_J)-v(2_J))(\frac{\nu_{2,x}}{y_{3_J}-y_{3_J}}-\frac{\nu_{2,y}}{x_{3_J}-x_{1_J}})]\\
\\
\gamma_{3,J}=[(v(1_J)-v(3_J))(\frac{\nu_{3,x}}{y_{3_J}-y_{3_J}}-\frac{\nu_{3,y}}{x_{3_J}-x_{3_J}})]
\end{array}
\right.
\end{equation}
 \end{lemma}
\begin{proof} 
Using Green's formula, where $\nu$ is the outer normal on $\partial E_J$, we obtain:
$$\left \langle \vec{V}(v) \cdot \nabla z,\phi \right \rangle_{L^2(E_J)}= \int_{\partial E_J}{z \nu\cdot (\phi\vec{V}(v))-\int_{E_J}{ z\nabla.(\phi\vec{V}(v)) } }$$ 
Since $\vec{V}(v)$ is divergence free, then:
\begin{eqnarray}\label{sumskew}
\left \langle \vec{V}(v) \cdot \nabla z,\phi \right \rangle_{L^2(E_J)}&=&\int_{\partial E_J}{z \nu\cdot (\phi\vec{V}(v))} -\int_{E_J}{ z\nabla.(\phi\vec{V}(v)) }= \int_{\partial E_J}{z \nu\cdot (\phi\vec{V}(v))} - \left \langle \vec{V}(v) \cdot \nabla\phi, z \right \rangle_{L^2(E_J)}\nonumber \\
\end{eqnarray}
Now the triangle boundary integral can be expressed as follows:
$$\int_{\partial E_J}{z \nu\cdot (\phi\vec{V}(v))} =  \int_{1_J}^{2_J}{z\phi \nu_1\cdot \vec{V}(v)}+\int_{2_J}^{3_J}{ z\phi \nu_2\cdot \vec{V}(v)}+\int_{3_J}^{1_J}{z\phi \nu_3\cdot \vec{V}(v)}$$
Handling as a sample one of these line integrals, one has for (example on $[1_J,2_J]$), $\nu\cdot \vec{V}(v)=\nu_{1,x}.v_{y}-\nu_{1,y}.v_x $. Furthermore as $v\in\mathcal{P}_1$, $ \nu_{1,x}.v_{y}-\nu_{1,y}.v_x $ is a constant on $(\overrightarrow{1_J2_J})$ and given by:
$$ \nu_{1,x}.v_{y}-\nu_{1,y}.v_x =\nu_{1,x}\frac{v(2_J)-v(1_J)}{y_{2_J}-y_{1_J}}-\nu_{1,y}\frac{v(2_J)-v(1_J)}{x_{2_J}-x_{2_J}}=(v(2_J)-v(1_J))(\frac{\nu_{1,x}}{y_{2_J}-y_{2_J}}-\frac{\nu_{1,y}}{x_{2_J}-x_{1_J}}).$$ 
Hence:
$$ \int_{1_J}^{2_J}{ z\phi \nu_1\cdot \vec{V}(v)}=(v(2_J)-v(1_J))(\frac{\nu_{1,x}}{y_{2_J}-y_{2_J}}-\frac{\nu_{1,y}}{x_{2_J}-x_{1_J}}) \int_{1_J}^{2_J}{z\phi}=\gamma_{1,J}\int_{1_J}^{2_J}{z\phi},$$
with similar identities obtained for $ \int_{2_J}^{3_J}{ z\phi \nu_2\cdot \vec{V}(v)}$ and $ \int_{2_J}^{3_J}{ z\phi \nu_3\cdot \vec{V}(v)}$. Replacing these integrals by their expressions in (\ref{sumskew}), one obtains the result of this lemma.
\end{proof}
\noindent The next step is to consider the sum $\sum_{E_J\in\mathcal{T}}{\left \langle \vec{V}(v) \cdot \nabla z,\phi \right \rangle_{L^2(E_J)}} $ and demonstrate the identity.
\begin{lemma}\label{sumlocalintegrals}
$$\ltwoinner{\vec{V}(v) \cdot \nabla z,\phi}= - \ltwoinner{\vec{V}(v) \cdot \nabla \phi,z}+\int_{\partial\Omega}{z\phi\,\nu. \vec{V}(v)}.
 $$
\end{lemma}
\begin{proof} On the basis of (\ref{breakskew}) and of lemma \ref{lem-est-0}, one has:
\begin{equation}\label{sumbreakskew}
\ltwoinner{\vec{V}(v) \cdot \nabla z,\phi}=-\sum_{E_J\in\mathcal{T}}{\left \langle \vec{V}(v) \cdot \nabla \phi,z \right \rangle_{L^2(E_J)}}+\sum_{E_J\in\mathcal{T}}{[\gamma_{1,J}\int_{1_J}^{2_J}{z\phi}+  \gamma_{2,J}\int_{2_J}^{3_J}{z\phi}+ \gamma_{3,J}\int_{3_J}^{1_J}{z\phi}]}.
 \end{equation}
 Given that any internal side $[AB]$ of any triangle $E_J$ is also common to another triangle $E_K$, then the corresponding line integral from $E_J$ is given by $ \gamma_{AB}\int_{A}^{B}{z\phi}$ and that coming from $E_K$ is $ \gamma_{BA}\int_{B}^{A}{z\phi}$, with $\gamma_{AB}=-\gamma_{BA} $, leading to a zero sum for integrals on $[AB]$.\\
 Consequently, one is left on the right hand side of (\ref{sumbreakskew}) with line integrals over  $\partial\Omega$, i.e. the result of the lemma.
\end{proof} 
\noindent Using the results of Lemmas \eqref{lem-est-0} and \eqref{sumlocalintegrals}, the we can complete the proof of Theorem \ref{skewxnP} as follows.
\begin{proof} The periodicity of $v$ and $z$ on $\partial\Omega$ results in  $\int_{\partial\Omega}{z\phi\,\nu.\vec{V}(v)}=0$, because of the periodicity of $z\phi$ on $\partial \Omega$ and the fact that $\nu. \vec{V}(v)$ is given by $ \pm v_x$ on the horizontal sides  and by  $ \pm v_y$ on the vertical sides. \\
 For example, handling the vertical sides gives 
\begin{eqnarray}
\int_{(0,0)}^{(1,0)} [z\phi\,\nu.\vec{V}(v)]dx + \int_{(1,1)}^{(0,1)} [z\phi\,\nu.\vec{V}(v)]dx &=& \int_{(0,0)}^{(1,0)}  -z\phi v_x dx + \int_{(1,1)}^{(0,1)}z\phi v_{x'} dx'\nonumber\\
&=& \int_{(0,0)}^{(1,0)}  -z\phi v_x dx + \int_{(0,1)}^{(1,1)}z\phi v_{-x} d(-x) \nonumber\\
&=& \int_{(0,0)}^{(1,0)}  -z\phi v_x dx + \int_{(0,1)}^{(1,1)}z\phi v_{x} d(x) = 0 \nonumber
\end{eqnarray}
Then using Lemma \eqref{sumlocalintegrals} we complete the proof of Skew-symmetry.
\end{proof}
\noindent This obviously leads to the following corollary. 
\begin{corollary}$\forall \{v,z\}\in X_{N,P}\times X_{N,P}$, we have that
$$\ltwoinner{\vec{V}(v) \cdot \nabla z,z}= 0$$
\end{corollary}

\subsection{Existence and Uniqueness of a solution to the Semi-Discrete System \eqref{HMC-Disc-V1}(\& \eqref{HMC-Disc-V1-1})}\label{Sub-Semi-Discrete}

\noindent In the rest of the paper, we will use the following norms in $\R^N$ for $V \in \R^N$:
\begin{eqnarray}
||V||^2_M:= V^T\,M\,V\\
||V||^2_K:= V^T\,K\,V
\end{eqnarray}

\noindent By associating the function $v_N(x,y)=\sum_{I=1}^N{V_I\varphi_I(x,y)}\in X_{N,P}$ to $V \in \R^N$, then we have the isometries:
\begin{eqnarray}
||V||^2_M&=&||v_N||^2 \label{eqnormM}\\
||V||^2_K&=&||v_N||_1^2 \label{eqnormK}
\end{eqnarray}

\noindent We begin by obtaining relations between solutions to (\ref{HMC-Disc-V1}) and (\ref{HMC-Disc-V1-1}). For $W(t),\; U(t) \in\mathbb{R}^N$ with 
$$w_N(x,y,t)=\sum_{I=1}^N{W_I(t)\varphi_I(x,y)}\in X_{N,P}\mbox{ and } u_N(x,y)=\sum_{I=1}^N{U_I(t)\varphi_I(x,y)}\in X_{N,P}$$
one has the following: 
\begin{lemma}\label{lemma-MS} 
Any pair $(U,W)$ that solves (\ref{HMC-Disc-V1-1}.2) satisfies the following:
$$||U(t)||_M\le ||U(t)||_K \leq ||W(t)||_M$$
\end{lemma}
\begin{proof} From the equation $KU(t)=MW(t)\Longleftrightarrow <u_N,v>_1=<w_N,v>,\,\forall v\in X_{N,P}$. Letting $v=u_N(t)$, yields $||u_N(t)||_1^2=<w_N,u_N>\le ||w_N(t)||.||u_N(t)||$, leading to:
$$||u_N(t)||^2\le||u_N(t)||_1^2 =<w_N,u_N>\le ||w_N(t)||.||u_N(t)|| \Longrightarrow ||u_N(t)||\le ||w_N(t)||,$$
which translates to the result via the isometries \eqref{eqnormM} and \eqref{eqnormK} .
\end{proof}

\noindent The existence of a unique solution to the semi-discrete system can be obtained by reducing (\ref{HMC-Disc-V1}) (or  (\ref{HMC-Disc-V1-1})) to a system of non-linear ordinary differential equations in $W(t)$. Specifically, in (\ref{HMC-Disc-V1-1}), we eliminate the variable $U(t)$, using the matrix $A=K^{-1}M$ and obtain  the system: 
\begin{equation}\label{HMC-Disc-V1-2}
\left\{\begin{array}{ll}
M\frac{dW}{dt} + S(AW)W= RAW,\, \forall t \in (0,T] &(1)\\
W(0)=W_0 & (2)\\
\end{array}\right.
\end{equation} 
which is equivalent to: 
\begin{equation}\label{HM-ODE}
\left\{\begin{array}{ll}
\frac{dW}{dt}= F(W(t)),\, \forall t \in (0,T] &(1)\\
W(0)=W_0 & (2)
\end{array}\right.
\end{equation}
where  
$$F(W)=M^{-1}[RAW-S(AW)W].$$

%
%

\noindent Now we show that $F$ is locally Lipschitz on the spaces $$\mathcal{X}_N=\{V \in \R^N \,|\,||V||_{M,C([0,T];\mathbb{R}^N)} \le C_T\}$$
\noindent where $C_T:= e^{||\hat{k}||_{\infty}T}||w_0||$ is determined in Lemma \ref{lem-est-1}.

\begin{lemma} \label{lemlipschitz} For $W_1,\,W_2\in \mathcal{X}_N$, there exists a positive constant $L_T$ independent from $h$ such that 
$$|| F(W_1)-F(W_2)||_M\le \frac{L_T}{h\sqrt{h}}||W_1-W_2||_M $$
\end{lemma}
\begin{proof}
Let $Z_k=F(W_k)$ with  $U_k=AW_k$ for $k=1,2$. \\
Let $z_{k,N}(x,y) =\sum_{I=1}^N{Z_{k,I}(t)\varphi_I(x,y)}$, $u_{k,N}(x,y) =\sum_{I=1}^N{U_{k,I}(t)\varphi_I(x,y)}$, $w_{k,N}(x,y) =\sum_{I=1}^N{W_{k,I}(t)\varphi_I(x,y)}$. 
Then, 
$$M(Z_1-Z_2)=R(U_1-U_2)-\left(S(U_1)W_1-S(U_2)W_2\right) $$ is equivalent to 
$$<z_{1,N}-z_{2,N},\phi> = {-}<\vec{V}(u_{1,N}).\nabla w_{1,N} - \vec{V}(u_{2,N}).\nabla w_{2,N}, \phi> + <\hat{k} (u_{1,N} - u_{2,N})_y, \phi>$$ 
for all $\phi \in X_{N,P}$.

\noindent Let $\phi = z_{1,N} - z_{2,N}$ then by Cauchy-Schwartz we get that
\begin{eqnarray}
||z_{1,N} - z_{2,N}||^2 &=& {-}<\vec{V}(u_{1,N}).\nabla w_{1,N} - \vec{V}(u_{2,N}).\nabla w_{2,N}, z_{1,N} - z_{2,N}> +  <\hat{k} (u_{1,N} - u_{2,N})_y, z_{1,N} - z_{2,N}> \nonumber \\
&\leq & ||\vec{V}(u_{1,N}).\nabla w_{1,N} - \vec{V}(u_{2,N}).\nabla w_{2,N}||\;||z_{1,N} - z_{2,N}||+ ||\hat{k} (u_{1,N} - u_{2,N})_y ||\;||z_{1,N} - z_{2,N}|| \nonumber 
\end{eqnarray}
 Simplifying by $||z_{1,N} - z_{2,N}||$ we get
 \begin{eqnarray}
||z_{1,N} - z_{2,N}|| &\leq & ||\vec{V}(u_{1,N}).\nabla w_{1,N} - \vec{V}(u_{2,N}).\nabla w_{2,N}||+ ||\hat{k} (u_{1,N} - u_{2,N})_y || \label{zn}
\end{eqnarray}
 Note that \begin{equation}||\hat{k} (u_{1,N} - u_{2,N})_y || \leq  ||\hat{k}||_{\infty} \;||(u_{1,N} - u_{2,N}) ||_{H^1} \leq  ||\hat{k}||_{\infty} \;||(w_{1,N} - w_{2,N}) || \label{un}\end{equation}

\noindent where we have used the fact that $K(U_1-U_1) = M(W_1-W_2) \implies ||u_{1,N}-u_{2,N}||_{H^1} \leq ||w_{1,N}-w_{2,N}||$.
 
\noindent On the other hand, using the triangle inequality
\begin{equation}
 ||\vec{V}(u_{1,N}).\nabla w_{1,N} - \vec{V}(u_{2,N}).\nabla w_{2,N}|| \leq ||\vec{V}(u_{1,N}- u_{2,N}).\nabla w_{1,N}|| + ||\vec{V}(u_{2,N}).\nabla (w_{2,N}-w_{1,N})|| \label{un1}\end{equation}
\noindent Note that for any $(\alpha_N, \beta_N) \in X_{N,p}\times X_{N,p}$, one has
\begin{equation}
 ||\vec{V}(\alpha_N).\nabla \beta_N|| \leq  \max\limits_{\phi \in X_{N,P}, ||\phi|| = 1} |<\vec{V}(\alpha_N).\nabla \beta_N, \phi>| \label{alphan}
 \end{equation}
 
\noindent Using skew-symmetry, $<\vec{V}(\alpha_N).\nabla \beta_N, \phi> = - <\vec{V}(\alpha_N).\nabla \phi,\beta_N>$, 
  then \eqref{alphan} becomes 
\begin{eqnarray}
||\vec{V}(\alpha_N).\nabla \beta_N|| &\leq&   \max\limits_{\phi \in X_{N,P}, ||\phi|| = 1} |<\vec{V}(\alpha_N).\nabla\phi, \beta_N>|
\end{eqnarray}
\end{proof}

 \noindent We need now the following Lemma.
 
 \begin{lemma}\label{lemma-skew}
 For some constant $C$ independent from $h$ we have:
 $$||\vec{V}(\alpha_N).\nabla \beta_N|| \leq  \frac{C}{h\sqrt{h}} ||\alpha_N||_1 ||\beta_N ||$$
 \end{lemma}
 \begin{proof}
 $$|<\vec{V}(\alpha_N).\nabla\phi, \beta_N>| \leq ||\vec{V}(\alpha_N).\nabla\phi|| \; ||\beta_N || \leq |\phi|_{1,\infty} ||\alpha_N||_1 ||\beta_N || \leq \frac{C}{h\sqrt{h}} ||\alpha_N||_1 ||\beta_N ||,$$
 as using a result in Ciarlet (\cite{ciarlet}, Theorem 3.2.6), one has 
 $$|\phi|_{1,\infty}\leq Ch^{-3/2}||\phi||.$$

 \end{proof}
\noindent We complete the proof, by applying this lemma twice to the right-hand side of \eqref{un1} thus obtaining
\begin{eqnarray}
 ||\vec{V}(u_{1,N}).\nabla w_{1,N} - \vec{V}(u_{2,N}).\nabla w_{2,N}|| &\leq& ||\vec{V}(u_{1,N}- u_{2,N}).\nabla w_{1,N}|| + ||\vec{V}(u_{2,N}).\nabla (w_{2,N}-w_{1,N})|| \nonumber\\
 &\leq &\frac{C}{h\sqrt{h}} ||u_{1,N}- u_{2,N}||_1 ||w_{1,N}|| + \frac{C}{h\sqrt{h}} ||u_{2,N}||_1 ||w_{2,N}-w_{1,N}||\
 \end{eqnarray}
Plugging this inequality and \eqref{un} in \eqref{zn} leads to:
 \begin{eqnarray}
||z_{1,N} - z_{2,N}|| &\leq & ||\vec{V}(u_{1,N}).\nabla w_{1,N} - \vec{V}(u_{2,N}).\nabla w_{2,N}||+ ||\hat{k}||_{\infty} \;||(w_{1,N} - w_{2,N}) ||\nonumber\\
&\leq & \frac{C}{h\sqrt{h}} ||u_{1,N}- u_{2,N}||_1 ||w_{1,N}|| + \frac{C}{h\sqrt{h}} ||u_{2,N}||_1 ||w_{2,N}-w_{1,N}||+||\hat{k}||_{\infty} \;||(w_{1,N} - w_{2,N}) || \nonumber
\end{eqnarray}
Using the assumption that $W_1,\,W_2\in \mathcal{X}_N$, the last inequality leads to 
\begin{eqnarray}
||z_{1,N} - z_{2,N}|| &\leq & \frac{CC_T}{h\sqrt{h}} ||w_{1,N}- w_{2,N}|| + \frac{CC_T}{h\sqrt{h}} ||w_{2,N}-w_{1,N}||+||\hat{k}||_{\infty} \;||w_{1,N} - w_{2,N} ||\nonumber \\
&\leq & (2\frac{CC_T}{h\sqrt{h}} + ||\hat{k}||_{\infty}  )||w_{1,N} - w_{2,N} || \leq \frac{L_T}{h\sqrt{h}} ||w_{1,N} - w_{2,N} ||
\end{eqnarray}
Now squaring and multiplying both sides by $M^T$ from the left, and using the equivalence of the $L_2$ norm on $X_{N,p}$ and the $M$ norm on $\mathbb{R}^N$ (\ref{eqnormM}) completes the proof of Lemma \ref{lemlipschitz}.

\noindent Thus, the semi-discrete system (\ref{HM-ODE}) has a unique solution $W(t)$ or a solution pair $\{U(t),W(t)\}$, for which we derive some a priori estimates.

\subsection{A Priori Estimates for Solutions to  \eqref{HMC-Disc-V1}} \label{sec:app}
\noindent We may now state some a priori estimates. 

%
%
%



\begin{lemma}\label{lem-est-1} Every unique solution $\{u_N,w_N\}$ to (\ref{HMC-Disc-V1}), satisfies the following estimates:
\begin{equation}\label{estimate-cinf}
\left\{\begin{array}{ll}
||w_N||_{C([0,T];L^2(\Omega))}=\max_{t\in[0,T]}{||w_N||(t)}\le C_T:=e^{||\hat{k}||_{\infty}T}||w_0|| & (1)\\
 ||u_N||_{C([0,T];H_P^1(\Omega))}=\max_{t\in[0,T]}{||u_N||_1(t)}\le C_T:=e^{||\hat{k}||_{\infty}T}||w_0||&(2)
\end{array}\right.
\end{equation}
\end{lemma}
\begin{proof} In (\ref{HMC-Disc-V1}.1), let $v=w_N$. Then one has:
\begin{eqnarray}
\ltwoinner{w_{N,t},w_N} + \ltwoinner{\vec{V}(u_N) \cdot \nabla w_N,w_N}= \ltwoinner{\hat{k}u_{N,y},w_N}\,\forall t \in (0,T], \label{eq1}
\end{eqnarray}
and letting $v=u_N$ in (\ref{HMC-Disc-V1}.2), one obtains also:
\begin{eqnarray}
\honeinner{u_N(t),u_N(t)}=\ltwoinner{w_N(t),u_N(t)},\, \forall t \in (0,T]. \label{eq2} \end{eqnarray}

\noindent Using Lemma \ref{lem-est-0} and the fact that $\ltwoinner{w_{N,t},w_N} = \frac{1}{2}\frac{d}{dt}||w_N||^2(t)$, then equations \eqref{eq1} and \eqref{eq2} lead to:
$$
\left\{\begin{array}{lll}
\frac{d}{dt}||w_N||^2(t) &\le { 2}||\hat{k}||_{\infty}||w_N || (t)||u_N ||_1(t),\, \forall t \in (0,T] &\\
||u_N||_1^2(t)&{\le ||w_N||(t)\,||u_N||(t),}\\
&\le ||w_N||(t)\,||u_N||_1(t),\,\forall t\in [0,T]\mbox{ given the initial choice of }u_N(0) & \\
\end{array}\right.
$$
Hence from these two inequalities, one gets $\forall t\in [0,T] $:
\begin{eqnarray}
 ||u_N||_1(t)&\le & ||w_N||(t),\, \label{ineq1}\\
\frac{d}{dt}||w_N||^2(t)& \le & ||\hat{k}||_{\infty}||w_N ||^2 (t). \label{ineq11} 
\end{eqnarray}
\noindent Integration of the differential inequality \eqref{ineq11} gives:
\begin{eqnarray}
||w_N||^2(t)&\le e^{{}||\hat{k}||_{\infty}t}||w_N||^2(0),\,\forall t\in [0,T] \nonumber\\
&\le e^{{ 2}||\hat{k}||_{\infty}T}||w_N||^2(0),\,\forall t\in [0,T] \nonumber \\
\therefore ||w_N||(t)& \le e^{||\hat{k}||_{\infty}T}||w_N||(0),\,\forall t\in [0,T] \label{ineq2}
\end{eqnarray}
\noindent Thus inequalities \eqref{ineq1} and \eqref{ineq2} give the results of the lemma, where $||w_N||(0) = ||\pi_N(w_0)|| \leq ||w_0||$.
\end{proof}

\subsection{Passing to the limit} \label{sec:2.3}
\noindent At this point we introduce the sequence $\{z_N\}$, defined by:
\begin{equation}\label{eq:reg}
z_N(t)\in H^1_P(\Omega)\cap H^2(\Omega) :\,-\Delta z_N(t)+z_N(t)=w_N(t)\mbox{ with }||z_N||_2(t)\le C||w_N||(t)
\end{equation}
 Note that in this case, the finite element approximation to $z_n$ in  $X_{N,P}$ is $u_N(t)$. Using the well-known Cea's  estimate for elliptic problems (\cite{ciarlet}, Theorem 3.2.2) in addition to \eqref{eq:reg}, one has  $\forall t\in[0,T]$:
 \begin{equation}\label{z-estimates}
||z_N(t)-u_N(t)||_1\le ||z_N(t)-\pi_N(u_N)(t)||_1\le \frac{C}{n}|z_N(t)|_{2,2} \le \frac{C}{n}||w_N(t)||
\end{equation}
\noindent with $C$ a generic constant independent of $N$. Thus, instead of studying the convergence of the pair $(u_N,w_N)$, we study that of $(z_N,w_N)$. \\
Following (\ref{estimate-cinf}.1) and (\ref{z-estimates}.1), one concludes that:
$$||z_N(t)||_{C([0,T];H^2)}\le C e^{||\hat{k}||_{\infty}T/2}||\pi_N(w_0)||.$$
\begin{lemma}\label{strong-conv-z} There exists an element $u\in L^2(0,T;H^1_P)$ and a subsequence $\{z_{N_i}\}\subset \{z_N\}$, such that:
$$\lim_{N_i\to\infty}{||u-z_{N_i}||_{L^2(0,T;H^1)}}=0. $$ 
\end{lemma}
\begin{proof}
This result follows from the Rellich-Kondrachov theorem that stipulates the compact injection of $H^2(\Omega)$ in $H^1(\Omega)$ (\cite{brezis}, p.285). 
.\end{proof}
\noindent Let us now denote $(z_{N_i},w_{N_i})$ by $(z_{N},w_{N})$ and seek first a limit point to the sequence  $\{w_N\}$. Specifically, we have the following result.

\begin{lemma}\label{weak-conv-w-1} There exists $w\in L^2(0,T;L^2(\Omega))$ and a subsequence
 $\{w_{N_j}\} $ of $\{w_N\}$ such that:
\begin{equation}
w_{N_j}(t){\rightharpoonup} w(t) \mbox{ in }L^2(0,T;L^2(\Omega))
\end{equation}
\begin{equation}\label{limit-points}
w_{N_j}(t){\rightharpoonup} w(t) \mbox{ in }L^2(\Omega) \,\, \mbox{ for all } t \in [0,T]
\end{equation} 
Furthermore $||w||_{L^2(0,T;L^2(\Omega))}\le e^{||\hat{k}||_{\infty}T}||w_0|| $.
\end{lemma}
\begin{proof} 

Observe via (\ref{estimate-cinf}.1) that the sequence $\{ w_N\}$ is uniformly bounded in the reflexive space $L^2(0,T;L^2(\Omega))$, and so it has a subsequence $\{ w_{N_j}\}$ which converges weakly, say to $w \in L^2(0,T;L^2(\Omega))$. i.e.
$$\int_0^T \inner{w_N(t),v} dt \too \int_0^T \inner{w(t),v} dt \mbox{ for all } v \in L^2(0,T;L^2(\Omega))$$
Now fix $v \in L^2(\Omega)$ and consider the sequence $F_j(t):=\int_0^t \inner{w_{N_j}(s),v} ds$ of functions on $[0,T]$. Observe that:
\begin{enumerate}
\item $F_j(t) \in C^1(0,T)$ and $F_j'(t)=\inner{w_{N_j}(t),v}$ for all $j$ by the Fundamental Theorem of Calculus.
\item $F_j(t) \too F(t):=\int_0^t \inner{w(s),v} ds$ pointwise on $[0,T]$.
\item $\norm{F_j(t)}$'s are uniformly bounded on $[0,T]$ by (\ref{estimate-cinf}.1).
\item $\{F_j(t)\}$ are uniformly equicontinuous on $[0,T]$ as
\eqnnolabel{\abs{F_j(t) - F_j(s)} &= \abs{\int_s^t \inner{w_{N_j}(\tau), v} d\tau} \\
&\leq \int_s^t \abs{\inner{w_{N_j}(\tau), v} } d\tau \\
&\leq \int_s^t \norm{w_{N_j}(\tau)} \norm{v} d\tau \leq C_T\norm{v}\abs{t-s}}
\end{enumerate}
\noindent so that by Arzel\`a-Ascoli theorem, $F_j(t)$ has a subsequence $F_{j_k}(t)$ that converges to $F(t)$ uniformly on $[0,T]$, and so $\inner{w_{N_{j_k}}(t),v}=F'_{j_k}(t) \too F'(t)=\inner{w(t),v}$ for every $t \in [0,T]$, which gives \eqref{limit-points} after relabelling.\\
\noindent Finally, weakly lower-semicontinuity of norms implies that $$||w||_{L^2(0,T;L^2(\Omega))}\le \liminf_{N\to\infty} ||w_N||_{L^2(0,T;L^2(\Omega))}\le e^{||\hat{k}||_{\infty}T}||w_0||$$
\end{proof}
\vskip 0.1cm
\noindent To complete the proof of Theorem \ref{main-theorem}, we denote the pair $(u_{N_j},w_{N_j})$ by $(u_{N},w_{N})$ and aim at proving  that the limit pair $(u,w)$ satisfies 
\begin{equation}\label{sys:toprove}
\left\{\begin{array}{ll}
\ltwoinner{w(t)-w_0,v} =\int_0^t \ltwoinner{\vec{V}(u(s)) \cdot \nabla v,w(s)}+ \hat{k}\ltwoinner{u_y(s),v} ds, &(1)\\
\honeinner{u,v}=\ltwoinner{w,v},\forall v \in H^1_P(\Omega),  &(2)
\end{array}\right.
\end{equation}
$\forall v \in H^1_P(\Omega)\cap W^{1,\infty}(\Omega), \ \forall t \in (0,T]$ with $w(0)=\pi_N(w_0)$.\\\\
\noindent For that purpose, we replace in \eqref{HMC-Comp-Mod}  $t$ by  $0$ and $t+\tau$ by $t$ and simultaneously  use the skew-symmetry property in Theorem \ref{skewxnP},  getting consequently:
\begin{equation}\label{HMC-Disc-int-w}
\left\{\begin{array}{ll}
\ltwoinner{w_{N}(t)-w_{N}(0),v_N} -\int_0^t{ \ltwoinner{\vec{V}(u_N(s)) \cdot \nabla v_N, w_N(s)}\,ds}= \int_0^t{\ltwoinner{\hat{k}u_{N,y}(s),v_N}ds}, &(1)\\
\honeinner{u_N(t),v_N}=\ltwoinner{w_N(t),v_N},&(2)\\
 \, \forall v_N \in X_{N,P}, \ \forall t \in (0,T]&\\
\end{array}\right.
\end{equation}
with $w_N(0)=\pi_N (u_0-\Delta u_0),$ which is the direct semi-discretization of (\ref{HMC-V2-sub}).
\\\\
To consider limit points when $N\to\infty$ of each of the terms in (\ref{HMC-Disc-int-w}.1) and (\ref{HMC-Disc-int-w}.2), the following sequence of lemmas is needed in which 
$C(T)$ is a generic constant of the form $ (aT+b)e^{d\,||\hat{k}||_\infty T}$ independent from $n$ where $a,b,d \in \mathbb{N}$.
\begin{lemma}\label{lim-w} For all $v\in H^1_P(\Omega)\cap W^{1,\infty}(\Omega)$, and for all $ t\in [0,T]$, one has: \vspace{-3mm}
\begin{equation}\label{lim-w-eq}
\ltwoinner{w_N(t),v}=\ltwoinner{w_N(t),\pi_N(v)} +\epsilon_{N,1}(t)\quad \mbox{ with } \quad |\epsilon_{N,1}(t)|\le C(T)\frac{1}{n} \;||w_0||\;|v|_{1}
\end{equation}
where $\epsilon_{N,1}(t)=\ltwoinner{w_N(t),\pi_N(v)-v}$
\end{lemma} 
\begin{proof} Given the identity:
$$ \ltwoinner{w_N(t),\pi_N(v)}=\ltwoinner{w_N(t),v} +\ltwoinner{w_N(t),\pi_N(v)-v},$$
and letting: $\epsilon_{N,1}(t)=\ltwoinner{w_N(t),v_N-v}$, one has using (\ref{density1}) and (\ref{estimate-cinf}): 
\begin{eqnarray}
|\epsilon_{N,1}(t)| &\le & ||w_N(t)||.||v-\pi_N(v)||, \\
&\le & e^{||\hat{k}||_\infty T} ||w_0|| C \; \dfrac{1}{n} |v|_{1} \quad \leq \quad C(T)\;\dfrac{1}{n}\;||w_0||\;|v|_{1}
\end{eqnarray}
where $C(T) = C \; e^{||\hat{k}||_\infty T}$.
\end{proof}
\noindent Similarly, one has:
\begin{lemma}\label{lim-u} For all $v\in H^1_P(\Omega)\cap W^{1,\infty}(\Omega)$, one has: \vspace{-3mm}
\end{lemma} 
\begin{equation}\label{lim-u-eq}
 \int_0^t{\ltwoinner{\hat{k}u_{N,y}(s),v}ds}=\int_0^t{\ltwoinner{\hat{k}u_{N,y}(s),\pi_N(v)}ds}+\epsilon_{N,2}(t)\quad \mbox{ with }\quad |\epsilon_{N,2}(t)|\le  C(T) \dfrac{1}{n} \,||w_0|| \,|v|_{1} 
\end{equation}
\begin{proof} Given the identity:
$$\int_0^t{\ltwoinner{\hat{k}u_{N,y}(s),v}\,ds}= \int_0^t{\ltwoinner{\hat{k}u_{N,y}(s),\pi_N(v)}\,ds}+\int_0^t{\ltwoinner{\hat{k}u_{N,y}(s),v-\pi_N(v)}},$$
and letting: $\epsilon_{N,2}(t)=\int_0^t{\ltwoinner{\hat{k}u_{N,y}(s),v-\pi_N(v)}}$, one has using (\ref{density1}) and (\ref{estimate-cinf}),: 
\begin{eqnarray}
|\epsilon_{N,2}(t)| &\le & ||v-\pi_N(v)||.\int_0^t{||u_N(s)||_1\,ds}, \\
& \leq & C\, \dfrac{1}{n} |v|_{1} \, T \, e^{||\hat{k}||_\infty T} ||w_0|| \quad \le \quad  C(T) \dfrac{1}{n} \,||w_0|| \,|v|_{1} 
\end{eqnarray}
where $C(T) = T \,C\, e^{||\hat{k}||_\infty T}$
\end{proof}
\noindent We turn now to the key term in (\ref{HMC-Disc-int-w}) and prove the following.
\begin{lemma}\label{key} For all $v\in H^1_P(\Omega)\cap W^{2,\infty}(\Omega)$, one has:
\begin{equation}\label{key-eq} 
\int_0^t{ \ltwoinner{\vec{V}(z_N(s)) \cdot \nabla v, w_N(s)}\,ds}=\int_0^t{ \ltwoinner{\vec{V}(u_N(s)) \cdot \nabla \pi_N(v), w_N(s)}\,ds}+\epsilon_{N,3}(t)
\end{equation}
where \vspace{-2mm}$$\epsilon_{N,3}(t)=\int_0^t{\ltwoinner{\vec{V}(z_N(s)-u_N(s)) \cdot \nabla v, w_N(s)}\,ds}+\int_0^t{\ltwoinner{\vec{V}(u_N(s)) \cdot \nabla (v-\pi_N(v)), w_N(s)}\,ds},$$
{with } \vspace{-5mm}$$|\epsilon_{N,3}(t)|\le   C(T)\,\dfrac{1}{n} ||w_0||^2\,\max \{|v|_{1,\infty},|v|_{2,\infty}\}.$$
\end{lemma}
\begin{proof} Adding up the following two identities and using the definition of $\epsilon_{N,3}(t)$:
$$ \int_0^t{ \ltwoinner{\vec{V}(z_N(s)) \cdot \nabla v, w_N(s)}\,ds}=\int_0^t{ \ltwoinner{\vec{V}(z_N(s)-u_N(s)) \cdot \nabla v, w_N(s)}\,ds}+\int_0^t{ \ltwoinner{\vec{V}(u_N(s)) \cdot \nabla v, w_N(s)}\,ds}$$
$$\int_0^t{ \ltwoinner{\vec{V}(u_N(s)) \cdot \nabla v, w_N(s)}\,ds}=\int_0^t{ \ltwoinner{\vec{V}(u_N(s)) \cdot \nabla \pi_N(v), w_N(s)}\,ds}+\int_0^t{ \ltwoinner{\vec{V}(u_N(s)) \cdot \nabla (v-\pi_N(v)), w_N(s)}\,ds}.$$
 yields:
$$ \int_0^t{ \ltwoinner{\vec{V}(z_N(s)) \cdot \nabla v, w_N(s)}\,ds}=\int_0^t{ \ltwoinner{\vec{V}(u_N(s)) \cdot \nabla \pi_N(v), w_N(s)}\,ds}+ \epsilon_{N,3}(t).$$
Thus, using (\ref{density12}), (\ref{estimate-cinf}) and (\ref{z-estimates}) one concludes the estimate:
\begin{eqnarray}
 |\epsilon_{N,3}(t)| &\le & | v-\pi_N(v)|_{1,\infty}.\int_0^t{ |u_N|_{1}(s).||w_N(s)||\,ds}+| v|_{1,\infty}.\int_0^t{ |u_N-z_N|_{1}(s).||w_N(s)||\,ds} \nonumber \\
 &\le &\dfrac{C}{n}|v|_{2,\infty}\int_0^t{ |u_N|_{1}(s).||w_N(s)||\,ds}+| v|_{1,\infty}.\int_0^t{ |u_N-z_N|_{1}(s).||w_N(s)||\,ds} \nonumber \\
 &\le &|v|_{2,\infty} \dfrac{C}{n} T{ e^{2||\hat{k}||_\infty T}||w_0||^2} + | v|_{1,\infty}\dfrac{C}{n}T{ e^{2||\hat{k}||_\infty T}||w_0||^2} \quad \le \quad   C(T)\,\dfrac{1}{n} ||w_0||^2\, \max \{|v|_{1,\infty},|v|_{2,\infty}\}\nonumber 
 \end{eqnarray}
 where $C(T) = 2\,T\, C\,  e^{2||\hat{k}||_\infty T}$.
\end{proof}
\noindent In a similar way to Lemma \ref{lim-w}, one can prove the following lemma:
\begin{lemma}\label{lim-u-1} For all $v\in H^1_P(\Omega)\cap W^{1,\infty}(\Omega)$, and for all $ t\in [0,T]$ one has: \vspace{-3mm}
\begin{equation}\label{lim-u-1-eq}
\honeinner{u_N(t),v}=\honeinner{u_N(t),\pi_N(v)}+\epsilon_{N,4}(t)\mbox{ with }|\epsilon_{N,4}(t)|\le C(T)\;\frac{1}{n} \;||w_0||\;|v|_{1}
\end{equation}
where $\epsilon_{N,4}(t) = \ltwoinner{u_N(t),\pi_N(v)-v}$ .
\end{lemma}
\subsection*{Synthesis: Completion of Proof of Existence (Theorem \ref{main-theorem})}
\noindent Using (\ref{lim-w-eq}), (\ref{lim-u-eq}), and (\ref{key-eq}) one gets: 
\begin{eqnarray}
\ltwoinner{w_{N}(t)-w_{N}(0),v} &=& \ltwoinner{w_{N}(t)-w_{N}(0),\pi_N(v)} +\epsilon_{N,1}(t)-\epsilon_{N,1}(0) \label{eq:1}\\
-\int_0^t{\ltwoinner{\hat{k}u_{N,y}(s),v}ds} &=& -\int_0^t{\ltwoinner{\hat{k}u_{N,y}(s),\pi_N(v)}ds}-\epsilon_{N,2}(t) \label{eq:3}\\
-\int_0^t{ \ltwoinner{\vec{V}(z_N(s)) \cdot \nabla v,w_N(s)}\,ds}&=& -\int_0^t{ \ltwoinner{\vec{V}(u_N(s)) \cdot \nabla \pi_N(v),w_N(s)}\,ds} - \epsilon_{N,3}(t)\label{eq:2}
\end{eqnarray}
Then, by summing up \eqref{eq:1}, \eqref{eq:2}, \eqref{eq:3},  and using (\ref{HMC-Disc-int-w}.1) , we get $\forall t\in (0,T)$:
\begin{equation}\label{before-final}
\ltwoinner{w_{N}(t)-w_{N}(0),v}-\int_0^t{ \ltwoinner{\vec{V}(z_N(s)) \cdot \nabla v,w_N(s)} + \ltwoinner{\hat{k}u_{N,y}(s),v}ds}=  \epsilon_{N,1}(t)-\epsilon_{N,1}(0)-\sum_{k=2}^3\epsilon_{N,k}(t)
\end{equation}
We may now let $N\to\infty$.
 Using the previous lemmas in this section,  we have subsequently:
\begin{enumerate}
\item For the right hand side of (\ref{before-final}), using Lemmas \ref{lim-w}, \ref{lim-u} and \ref{key}:
$$\lim_{N\to\infty}\epsilon_{N,1}(t)-\epsilon_{N,1}(0)-\sum_{k=2}^3\epsilon_{N,k}(t) =0, \quad \forall v \in W^{2,\infty}(\Omega)\cap H^{1}_P(\Omega)$$
\item For the left-hand side of (\ref{before-final}):
\begin{itemize}
\item $\lim_{N\to\infty}{\ltwoinner{w_{N}(t)-w_{N}(0),v}}=\ltwoinner{w(t)-w(0),v}$, using Lemma \ref{lim-w}.
\item $ \lim_{N\to\infty}\int_0^t{\ltwoinner{\hat{k}u_{N,y}(s),v}ds}=\int_0^t{\ltwoinner{\hat{k}u_{y}(s),v}ds}$,  using Lemma \ref{lim-u}.
\item For the term $\int_0^t{ \ltwoinner{\vec{V}(z_N(s)) \cdot \nabla v,w_N(s)}\,ds}$, note that for all $s \in [0,t]$ and for all $v \in W^{1,\infty}(\Omega)$, $\vec{V}(z_N(s)).\nabla v$ converges to $\vec{V}(u(s)).\nabla v$ strongly in $L^2(0,T;L^2(\Omega))$, as
$$\int_0^t \inner{\vec{V}(z_N(s)-u(s)) \nabla v, \varphi} ds \leq \int_0^t \norm{z_N(s)-u(s)}_1 \norm{v}_{W^{1,\infty}}  \norm{\varphi} ds \to 0 \mbox{ as } N \to 0$$

Combining this with the weak convergence of $w_N$ to $w$ in $L^2(0,T;L^2(\Omega))$, we obtain
$$\lim_{N\to\infty}{\int_0^t{ \ltwoinner{\vec{V}(z_N(s)) \cdot \nabla v,w_N(s)}\,ds}}=\int_0^t{ \ltwoinner{\vec{V}(u(s)) \cdot \nabla v,w(s)}\,ds} $$ 
\end{itemize}
\item By applying similar techniques with respect to $$\honeinner{u_N(t),v_N}=\ltwoinner{w_N(t),v_N},\forall v_N \in X_{N,P}, \, \forall t \in [0,T] $$ in (\ref{HMC-Disc-int-w}.2), using Lemmas \ref{lim-w} and \ref{lim-u-1} one obtains as $N\to\infty$:
$$\honeinner{u(t),v}=\ltwoinner{w(t),v},\forall v \in H^1_P(\Omega), \, \forall t \in [0,T]. $$
\end{enumerate}
Thus these last 3 consecutive points prove formulation \eqref{HMC-V2-sub} for test functions $v \in W^{2,\infty}(\Omega) \cap H^1_P(\Omega)$. \\\\
\noindent Finally, the density of $W^{2,\infty}(\Omega) \cap H^1_P(\Omega)$ in $W^{1,\infty}(\Omega)\cap H^1_P(\Omega)$ completes the proof of Theorem \ref{main-theorem}.
\section{Full Discretization}\label{sec:Full-Discrete}
As to fully discretizing the Hasegawa-Mima system, starting with   (\ref{HMC2-tau}) and to avoid any constraint of the  CFL type on the choice of $\tau$, the term $\int_t^{t+\tau}{\ltwoinner{\vec{V}(u(s)) \cdot \nabla v,w(s)}}$ is first discretized using an implicit right rectangular rule:
$$\int_t^{t+\tau} {\ltwoinner{\vec{V}(u(s)) \cdot \nabla v,w(s)}}={\tau}\ltwoinner{\vec{V}(u(t+\tau)) \cdot \nabla v,w(t+\tau)} +\epsilon_\tau,$$ leading to the following fully implicit {Computational Model}. Given $(u_N(t),w_N(t)) \in X_{N,P}\times X_{N,P}$, one seeks\linebreak $(u_N(t+\tau),w_N(t+\tau))\in  X_{N,P}\times X_{N,P}$, such that:
\begin{equation}\label{HMC-Comp}
\left\{\begin{array}{ll}
\ltwoinner{w_N(t+\tau)-w_N(t),v} = \tau\ltwoinner{\vec{V}(u_N(t+\tau)) \cdot \nabla v,w_N(t+\tau)} + {\tau}\ltwoinner{\hat{k}u_{N,y}(t+\tau),v}, &(1)\\
\honeinner{u_N(s),v}=\ltwoinner{w_N(s),v},   &(2)\\
\,\forall v\in  X_{N,P}, \quad \forall s\in\{t,t+\tau\}&
\end{array}\right.\end{equation} 
In matrix notations and using the expressions: $$w_N(t)=\sum_{I=1}^N{W_I(t)\varphi_I(x,y)},\, \mbox{ and } u_N(x,y,t)=\sum_{J=1}^N{U_J(t)\varphi_J(x,y)}, $$  where $W_I(t)=w_N(x_I,y_I,t)$, and $U_J(t)=w_N(x_J,y_J,t),$ then (\ref{HMC-Comp}) can be rewritten as follows:\\
Given $(U(t),W(t)) \in \mathbb{R}^{N}\times  \mathbb{R}^{N}$, seek $(U(t+\tau),W(t+\tau))\in  \mathbb{R}^{N}\times  \mathbb{R}^{N}$, such that:
\begin{equation}\label{HMC-Disc-Comp}
\left\{\begin{array}{ll}
(M+ \tau\,S(U(t+\tau))\;W(t+\tau)-\tau\, R\;U(t+\tau)=M\,W(t) &(1)\\
KU(s)=MW(s), \;\;\; \forall s \in \{t,t+\tau\} &(2)\\
\end{array}\right.
\end{equation}
In Section \ref{Sub-Full-Discrete}, using a fixed point approach we start by showing the existence of solution to \eqref{HMC-Disc-Comp}, i.e. \eqref{HMC-Comp}. Then we prove uniqueness of this solution in Section \ref{uniqueness}.{
\subsection{Existence of Solution to the Fully Discrete System}\label{Sub-Full-Discrete}
To prove the existence of a solution to the Fully Discrete System \eqref{HMC-Disc-Comp}, we start by transforming it into the fixed point problem \eqref{fixed3}. Then, we prove the existence of a solution to  \eqref{fixed3} using the Leray-Schauder fixed point Theorem \cite{mawhin}.\\

\noindent Using equation (\ref{HMC-Disc-Comp}.2) for $s = t+\tau$, one gets $U(t+\tau) = K^{-1}MW(t+\tau)$. Substituting $U(t+\tau)$ in  equation (\ref{HMC-Disc-Comp}.1), system \eqref{HMC-Disc-Comp} can be rewritten as follows
\begin{eqnarray}
\big{[}M +{\tau}\,S(K^{-1}MW(t+\tau))-{\tau}RK^{-1}M \big{]}W(t+\tau)=  M W(t), \label{fixed1}
\end{eqnarray} 
which is equivalent to: 
\begin{eqnarray}
\big{[}M -{\tau}RK^{-1}M \big{]}W(t+\tau)=  M W(t) -{\tau}\,S(K^{-1}MW(t+\tau))W(t+\tau), \label{fixed1p}
\end{eqnarray} 
Let $B = M -{\tau}RK^{-1}M $, $Z = W(t)$ and $Y = W(t+\tau)$, then we get the first fixed point 
\begin{eqnarray}
BY=  M Z -{\tau}\,S(K^{-1}MY)Y, \label{fixed2}
\end{eqnarray} 
\begin{theorem}\label{bilin} Define the bilinear form on $H_1(\Omega)\times H_1(\Omega)$: $$a_{\tau}(w,v) =  < w, v>_1 - \tau <\hat{k}w_y,v >.$$ Then,  for $\tau \leq \dfrac{1}{2||\hat{k}||_{\infty}}$ this bilinear form  is coercive in the sense that 
$$|a_{\tau}(w,w)| \geq \dfrac{1}{2}||w||_1^2, \qquad \forall w \in H_1(\Omega)$$
\end{theorem}
\begin{proof}
This simply results from: 
$$|a_{\tau}(w,w)| = |< w, w>_1 - \tau <\hat{k}w_y,w >| \geq \left(1 -\tau  ||\hat{k}||_{\infty}\right)||w||_1^2 \geq \dfrac{1}{2}||w||_1^2, \qquad for \;\; \tau \leq \dfrac{1}{2||\hat{k}||_{\infty}}$$
\end{proof}
\noindent As a consequence, one gets the following corollary.
\begin{corollary} The matrix $B$ is invertible for $\tau \leq \dfrac{1}{2||\hat{k}||_{\infty}}$. 
\end{corollary}
\begin{proof}
\begin{equation}
B\alpha = \beta \iff (I  -{\tau}RK^{-1})M  \alpha = \beta \iff (I  -{\tau}RK^{-1})\alpha^{(1)} = \beta
\end{equation}
 where $\alpha^{(1)} = M \alpha$.
Let $\alpha^{(2)} =  K^{-1} \alpha^{(1)}$, i.e. $\alpha^{(1)} = K \alpha^{(2)}$, and $\beta^{(1)} = M^{-1}\beta$. Hence,  
$$ (K  -{\tau}R)\alpha^{(2)} = M\beta^{(1)}$$
In variational form, this is equivalent to 
\begin{equation}\label{varFull} 
 < \alpha_N^{(2)}, v>_1 - \tau <\hat{k}\alpha_{N,y}^{(2)},v >= <\beta_N^{(1)},v>  
\end{equation}
for all $v \in X_{N,P}$, with $ \beta_N^{(1)} = \sum\limits_{i=1}^N {\beta_i^{(1)} \varphi_i}$,  $\alpha_N^{(2)}= \sum\limits_{i=1}^N {\alpha_{i}^{(2)} \varphi_i}$, and  $\alpha_{N,y}^{(2)}= \sum\limits_{i=1}^N {\alpha_{i,y}^{(2)} \varphi_i}$.

\noindent Since from theorem \ref{bilin}, the bilinear form $a_{\tau}(w,v) =  < w, v>_1 - \tau <\hat{k}w_y,v >$ is coercive,
then by applying Lax-Milgram on \eqref{varFull}, one completes the proof of the corollary. \\
\end{proof}
\noindent Given that $B$ is invertible for $ \tau \leq \dfrac{1}{2||\hat{k}||_{\infty}}$, the fixed point format \eqref{fixed2} becomes
\begin{equation} Y = G(Y) = B^{-1}\big[ M Z -{\tau}\,S(K^{-1}MY)Y\big]. \label{fixed3}\end{equation}
Using the Leray-Schauder fixed point theorem in $\mathbb{R}^N$, we prove the existence of a solution to the fixed point problem \eqref{fixed3}.

\begin{theorem}(Leray-Schauder Theorem \cite{mawhin} )
Let $\mathcal{T}$ be a continuous and compact mapping of a Banach space $X$ into itself, such that the set
$${ \{x\in X:x=\lambda \mathcal{T}x{\mbox{ for some }}0\leq \lambda \leq 1\}} $$
is bounded. Then $\mathcal{T}$ has a fixed point.
\end{theorem} 
\begin{theorem} Let $0\leq s \leq 1$ and consider the fixed point problem 
\begin{equation} B Y_s = s\big[ M Z -{\tau}\,S(K^{-1}MY_s)Y_s\big]\label{LS1}
\end{equation}  then for $ \tau \leq \dfrac{1}{2 ||\hat{k}||_{\infty}}$, the fixed point problem $Y = G(Y)$ admits a solution.
\end{theorem}
\begin{proof}
Multiplying equation \eqref{LS1} by $Y_s^T$ we get:
\begin{eqnarray}
Y_s^T B Y_s &=& sY_s^T M Z -{\tau}sY_s^T\,S(K^{-1}MY_s)Y_s \nonumber
\end{eqnarray}
Using skew-symmetry of the matrix $S(K^{-1}MY_s)$ we get:
\begin{eqnarray}
Y_s^T B Y_s &=& sY_s^T M Z \leq s||Y_s||_M ||Z||_M \\
Y_s^T (M - \tau RK^{-1}M) Y_s &=& ||Y_s||_M^2 - \tau Y_s^TRK^{-1}MY_s   \leq s ||Y_s||_M ||Z||_M \\
||Y_s||_M^2 &\leq & s ||Y_s||_M ||Z||_M +\tau Y_s^TRK^{-1}MY_s = s ||Y_s||_M ||Z||_M +\tau Y_s^TRW_s\label{LS2}
\end{eqnarray}
where $KW_s =MY_s$. Let $w_{s,N} = \sum\limits_{i=1}^N {W_{s,i}\varphi_i}$ and $y_{s,N} = \sum\limits_{i=1}^N {Y_{s,i}\varphi_i}$. Then, $$||w_{s,N}||_2 = ||W_s||_M\qquad and \qquad||y_{s,N}||_2 = ||Y_s||_M.$$ Also, \mbox{from} $KW_s =MY_s$ one gets: \begin{equation}
 <w_{s,N}, v>_1 = <y_{s,N},v>  \quad \implies ||w_{s,N}||_1 \leq ||y_{s,N}||_2\label{LS3}\end{equation} and  $$Y_s^T RW_s = <\hat{k}(w_{s,N})_y, y_{s,N}>.$$
Hence, using the last two equations, \eqref{LS2} becomes 
\begin{eqnarray}
||Y_s||_M^2 &\leq & s ||Y_s||_M ||Z||_M +\tau Y_s^TRW_s \leq ||Y_s||_M ||Z||_M + \tau ||\hat{k}||_{\infty} ||y_{s,N}||_2^2\\
&\leq & ||Y_s||_M ||Z||_M    + \tau ||\hat{k}||_{\infty} ||Y_s||_M^2\\
\implies (1 -\tau ||\hat{k}||_{\infty}) ||Y_s||_M^2 &\leq & ||Y_s||_M ||Z||_M 
\end{eqnarray}
Then for $\tau \leq \dfrac{1}{2||\hat{k}||_{\infty}}$ we get  $||Y_s||_M \leq  2 ||Z||_M$. Thus all solutions $Y_s$ are uniformly bounded, leading to the existence of a fixed point for all $0\leq s \leq 1$, using the  Leray-Schauder fixed point theorem. 
\end{proof}
\subsection{Uniqueness of Solution of the Fully Discrete System}\label{uniqueness}
\noindent Now we turn to uniqueness and consider the ball \begin{equation}\label{ball}
\mathcal{B}_N = \{Y \in \mathbb{R}^N \;\;|\;\; ||Y||_{M}\leq C_Z = 2 ||Z||_M\}
\end{equation}
we then prove the following theorem.
\begin{theorem}\label{fixProof} For $Y^{(k)} \in \mathcal{B}_N,\,k = 1,2$ and 
$\tau \leq \dfrac{1}{{\color{red} 8}}\min \left\{ \dfrac{4}{||\hat{k}||_{\infty}},\dfrac{h^2}{c C_Z}\right\} $, 
 one has:
  $$ ||G(Y^{(1)}) - G(Y^{(2)})||_M \leq \frac{1}{2}\;||Y^{(1)} - Y^{(2)}||_2$$
  and hence under such restriction on $\tau$ and $h$,  the fixed problem $Y=G(Y)$ admits a unique solution.
  \end{theorem}
\begin{proof}
For $ k=1,2$, let
\begin{eqnarray}
\Lambda^{(k)} &=& G(Y^{(k)})\\
K\Sigma^{(k)} &=& MY^{(k)} \label{62a}\\
 K\Delta^{(k)}   &=& M\Lambda^{(k)}\label{62b}
\end{eqnarray}
Then by applying this equation for $k = 1$ and $2$, followed by a subtraction, we get
 \begin{eqnarray}
 (M -{\tau}RK^{-1}M )( \Lambda^{(1)} - \Lambda^{(2)}) &=& {\tau}\,S(K^{-1}MY^{(2)})Y^{(2)} - {\tau}\,S(K^{-1}MY^{(1)})Y^{(1)} \nonumber \\
 &=& {\tau}\,S(\Sigma_2)Y^{(2)} - {\tau}\,S(\Sigma_1)Y^{(1)}\nonumber 
 \end{eqnarray}
Then, by using the newly introduced variables in \eqref{62a}, and \eqref{62b}, we obtain 
\begin{equation}
 (K - {\tau}R)(\Delta^{(1)} - \Delta^{(2)}) = {\tau}\,S(\Sigma^{(2)})Y^{(2)} - {\tau}\,S(\Sigma^{(1)})Y^{(1)} \label{61}
\end{equation}
Let $ Y^{(k)}_N = \sum\limits_{i=1}^N {Y_i^{(k)} \varphi_i}$ and $\epsilon_N =  Y^{(1)}_N -  Y^{(2)}_N$,  then $
||\epsilon_N ||^2 = ||Y^{(1)} - Y^{(2)}||_M^2$. \\
Let $ \Lambda^{(k)}_N = \sum\limits_{i=1}^N {\Lambda^{(k)} \varphi_i}$ and $\lambda_N = \Lambda^{(1)}_N - \Lambda^{(2)}_N$, then
$||\lambda_N ||^2 = ||G(Y^{(1)}) - G(Y^{(2)})||_M^2$.\\
In variational form, equations , \eqref{62a}, \eqref{62b} and \eqref{61} lead to:
\begin{eqnarray}
<\Delta_N^{(k)}, v>_1 - \tau <\hat{k}(\Delta_N^{(k)})_y,v >&=&<Z_{N},v> -\tau <\vec{V}(\Sigma^{(k)}_N).\nabla Y^{(k)}_N , v> \label{63k}\\
<\Delta_N^{(1)} - \Delta_N^{(2)}, v>_1 - \tau <\hat{k}(\Delta_N^{(1)} - \Delta_N^{(2)})_y,v >&=& -\tau <\vec{V}(\Sigma^{(1)}_N).\nabla Y^{(1)}_N - \vec{V}(\Sigma^{(2)}_N).\nabla Y^{(2)}_N, v> \label{63} \\
<\Sigma^{(k)}_N,v>_1 &=& <Y^{(k)}_N,v>\label{64} \\
<\Delta_N^{(k)},v>_1 &=& <\Lambda^{(k)}_N,v>   \label{65}
\end{eqnarray}
For all $v \in X_{N,P}$, where $ Z_{N} = \sum\limits_{i=1}^N {Z_i\varphi_i}$, $ \Sigma^{(k)}_N = \sum\limits_{i=1}^N {\Sigma_i^{(k)} \varphi_i}$, \;$\Delta_N^{(k)}= \sum\limits_{i=1}^N {\Delta_{i}^{(k)} \varphi_i}$,\;  and $\Delta_{N,y}^{(k)}= \sum\limits_{i=1}^N {\Delta_{i,y}^{(k)} \varphi_i}$,   for $k=1,2$.\\
 Let us define $\delta_N = \Delta_N^{(1)} - \Delta_N^{(2)}$, $\sigma_N = \Sigma^{(1)}_N - \Sigma^{(2)}_N$, and . Using again the bilinear form $a_{\tau}(.,.)$ defined in theorem \ref{bilin},
\begin{equation}
a_{\tau}(\delta_N, v) = -\tau <\vec{V}(\sigma_N).\nabla Y^{(1)}_N - \vec{V}(\Sigma^{(2)}_N).\nabla \epsilon_N, v>  
\end{equation} 
Then,  \eqref{63}, \eqref{64}, and \eqref{65} lead to:
\begin{eqnarray}
a_{\tau}(\delta_N, v)&=& -\tau <\vec{V}(\sigma_N).\nabla Y^{(1)}_N - \vec{V}(\Sigma^{(2)}_N).\nabla \epsilon_N, v>\label{66}  \\
<\sigma_N,v>_1 &=& <\epsilon_N,v>  \label{69}\\
  <\delta_N,v>_1 &=&<\lambda_N,v>  \label{68}
\end{eqnarray}
for all $v \in X_{N,P}$.\\ 
Now, let $v = \delta_N$ in \eqref{66}, then for $\tau \leq \dfrac{1}{2||\hat{k}||_{\infty}}$, we have:
\begin{eqnarray}
\dfrac{1}{2}||\delta_N||_1^2 &\leq & \tau |<\vec{V}(\sigma_N).\nabla Y^{(1)}_N,\delta_N >| + \tau |<\vec{V}(\Sigma^{(2)}_N).\nabla \epsilon_N, \delta_N>|\label{70}
\end{eqnarray}
Since, theorem \eqref{skewxnP} asserts that $$\forall \{v,z,\phi\}\in X_{N,P}\times X_{N,P}\times X_{N,P}:\,\ltwoinner{\vec{V}(v) \cdot \nabla z,\phi}= - \ltwoinner{\vec{V}(v) \cdot \nabla \phi,z},$$
then \eqref{70} can be rewritten as: 
\begin{eqnarray}
\dfrac{1}{2}||\delta_N||_1^2 &\leq & \tau |<\vec{V}(\sigma_N).\nabla \delta_N, Y^{(1)}_N>| + \tau |<\vec{V}(\Sigma^{(2)}_N).\nabla \delta_N, \epsilon_N>|\\
&\leq & \tau ||\vec{V}(\sigma_N).\nabla \delta_N|| \;||Y^{(1)}_N|| + \tau ||\vec{V}(\Sigma^{(2)}_N).\nabla  \delta_N||\; ||\epsilon_N||
\end{eqnarray}
Using Lemma \ref{lemma-skew}, this last inequality leads to: 
$$\dfrac{1}{2}||\delta_N||_1^2 \leq \frac{\tau}{h}||\delta_N||\left[C_Z||\sigma_N||_1+||\Sigma^{(2)}_N||_1.||\epsilon_N||\right].$$
Therefore, given that $||\delta_N||\le ||\delta_N||_1$, it results that:
\begin{eqnarray}\label{lastdelta}
\dfrac{1}{2}||\delta_N||_1 &\leq & \frac{\tau}{h} \left[C_Z||\sigma_N||_1+||\Sigma^{(2)}_N||_1.||\epsilon_N||\right].
\end{eqnarray}
Note that if in \eqref{69}, $v = \sigma_N$, one proves that:
$$||\sigma_N||_1\le ||\epsilon_N||, $$
and if in \eqref{64}, $k=2$ and $v=\Sigma^{(2)}_N$, we obtain:
$$||\Sigma^{(2)}_N||_1\le ||Y_N^{(2)}||\le C_Z. $$
Combining the last two inequalities with inequality \eqref{lastdelta}, we conclude the inequality:
\begin{eqnarray}\label{lastdelta2}
\dfrac{1}{2}||\delta_N||_1 &\leq& { 2}\frac{\tau C_Z}{h}||\epsilon_N||.
\end{eqnarray}
Finally, if in \eqref{68}, we let $v = \lambda_N$, then:
 $$||\lambda_N||^2 = |<\delta_N,\lambda_N>_1| \leq ||\delta_N||_1 \;||\lambda_N||_1 $$
{Using a result from Ciarlet (\cite{ciarlet}, Theorem 3.2.6), one has: $||\lambda_N||_1 \leq \dfrac{c}{h}||\lambda_N||$}, then 
\begin{equation} \label{eq:lambda}
||\lambda_N|| \leq \dfrac{c}{h} ||\delta_N||_1  
\end{equation}
When combining \eqref{lastdelta2} and \eqref{eq:lambda} we reach the final result:
$$ ||G(Y_1)-G(Y_2)||_M=||\Lambda_N^{(1)}-\Lambda^{(2)}_N||=||\lambda_N|| \leq \dfrac{c}{h} ||\delta_N||_1\leq {  4} \tau\; \frac{cC_Z}{h^2}||\epsilon_N||={  4}  \tau\; \frac{cC_Z}{h^2}||Y_1-Y_2||_M.$$
On the other hand, considering \eqref{63k}, one writes:
$$a_\tau(\Delta_N^{(k)},v)=<Z_{N},v> -\tau <\vec{V}(\Sigma^{(k)}_N).\nabla Y^{(k)}_N , v>,\,k=1,2 $$
Consequently, our Theorem is proved 
%
\end{proof}
}
\begin{corollary}\label{corol}
Under the conditions of theorem \ref{fixProof}, namely $\tau \leq \dfrac{1}{{ 8}}\min \left\{ \dfrac{4}{||\hat{k}||_{\infty}},\dfrac{h^2}{c C_Z}\right\}$, 
the iteration \\$Y^{(k+1)} = G(Y^{(k)})$ with $Y^{(0)}=Z$ in the ball $\mathcal{B}_N$ converges to the unique solution of $Y = G(Y)$.
\end{corollary}
\begin{proof}
Starting with 
  $$ ||Y^{(1)} - Y||_M = ||G(Z) - G(Y)||_M \leq \frac{1}{2}\;||Z - Y||_2$$
  then by induction we get that  
    $$ ||Y^{(k+1)} - Y||_M = ||G(Y^{(k)}) - G(Y)||_M \le \frac{1}{2}\;||Y^{(k)} - Y||_2  \leq \frac{1}{2^{k+1}}\;||Y^{(0)} - Y||_2$$
    where \quad $ 0 \leq \lim\limits_{k \rightarrow \infty} ||Y^{(k+1)} - Y||_M \leq \lim\limits_{k \rightarrow \infty} \dfrac{1}{2^{k+1}}\;||Y^{(0)} - Y||_2 = 0$.
\end{proof}
\section{Algorithm and Computer Simulations} \label{Sec:Algorithm}
Recall from  \eqref{HMC-Disc-Comp} that in matrix notations and using the expressions: $$w_N(t)=\sum_{I=1}^N{W_I(t)\varphi_I(x,y)},\, \mbox{ and } u_N(x,y,t)=\sum_{J=1}^N{U_J(t)\varphi_J(x,y)}, $$  where $W_I(t)=w_N(x_I,y_I,t)$, and $U_J(t)=w_N(x_J,y_J,t),$ then (\ref{HMC-Comp}) can be rewritten as follows:\\
Given $(U(t),W(t)) \in \mathbb{R}^{N}\times  \mathbb{R}^{N}$, seek $(U(t+\tau),W(t+\tau))\in  \mathbb{R}^{N}\times  \mathbb{R}^{N}$, such that:
$$
\left\{\begin{array}{ll}
(M+ \tau\,S(U(t+\tau))\;W(t+\tau)-\tau\, R\;U(t+\tau)=M\,W(t) &(1)\\
KU(s)=MW(s), \;\;\; \forall s \in \{t,t+\tau\} &(2)\\
\end{array}\right.
$$
\noindent To solve \eqref{HMC-Disc-Comp} we use in this paper a semi-linearized approach:\\
 \begin{equation}\label{semi-lin-var}
\left\{\begin{array}{ll}
\ltwoinner{w_N(t+\tau),v} -  \tau\ltwoinner{\vec{V}(u_N(t)) \cdot \nabla v,w_N(t+\tau)}=\ltwoinner{w_N(t),v} + {\tau}\ltwoinner{\hat{k}u_{N,y}(t),v}, &(1)\\
\honeinner{u_N(s),v}=\ltwoinner{w_N(s),v},\quad s\in\{t,t+\tau\}  &(2)\\
\end{array}\right.
\end{equation} 
 $\forall v\in  X_{N,P}$, which in matrix form is given by:
\begin{equation}\label{HMC-Disc-Comp-semi}
\left\{\begin{array}{ll}
(M+ \tau\,S(U(t))\;W(t+\tau)=M\,W(t)+\tau\, R\;U(t) &(1)\\
KU(t+\tau)=MW(t+\tau), \;\;\;  &(2)\\
\end{array}\right.
\end{equation}
\noindent where $M, K$, $S(U)$ and $R$ are $N\times N$ matrices defined in \eqref{HMC-Disc-V1-1}. 
However, by taking periodicity into account, the degrees of freedom are reduced from $N = n^2$ to $N_1 = (n-1)^2$. Note that $M$, is the well-known Mass matrix for periodic boundary conditions and $K = M+A$ where $A$ is the stiffness matrices for periodic boundary conditions. The nonlinearity of the problem originates from $S(U)$, which we derive its 
 corresponding local matrix, in addition to that of $R$, over each triangle for equally spaced nodes in Section \ref{sec:SR}. Then, deduce the block sparsity pattern of the global matrix, which is the same for $M$ and $K$.

\noindent Thus, to solve  \eqref{HMC-Disc-Comp-semi} these matrices should be generated for a given meshing. In Section \ref{sec:4.1}, we implement Algorithm \eqref{alg:HMC-Newton} using Freefem++ \cite{MR3043640}, a programming language and software focused on solving partial differential equations using the finite element method.\vspace{2mm}

\noindent Using different initial conditions, we test in Section \ref{Tests} our algorithm for the case when $p = \ln \dfrac{n_0}{\omega_{ci}}$ is a function of $x$, such that $\hat{k} = p_x$ is  a constant and $p_y = 0$,  and for the case when $p,p_x,$ and $p_y$ are functions of $(x,y)$. 
\subsection{Expressions of $S(U)$ and $R$}\label{sec:SR}
\begin{figure}[H]
\centering
\includegraphics[scale=0.17]{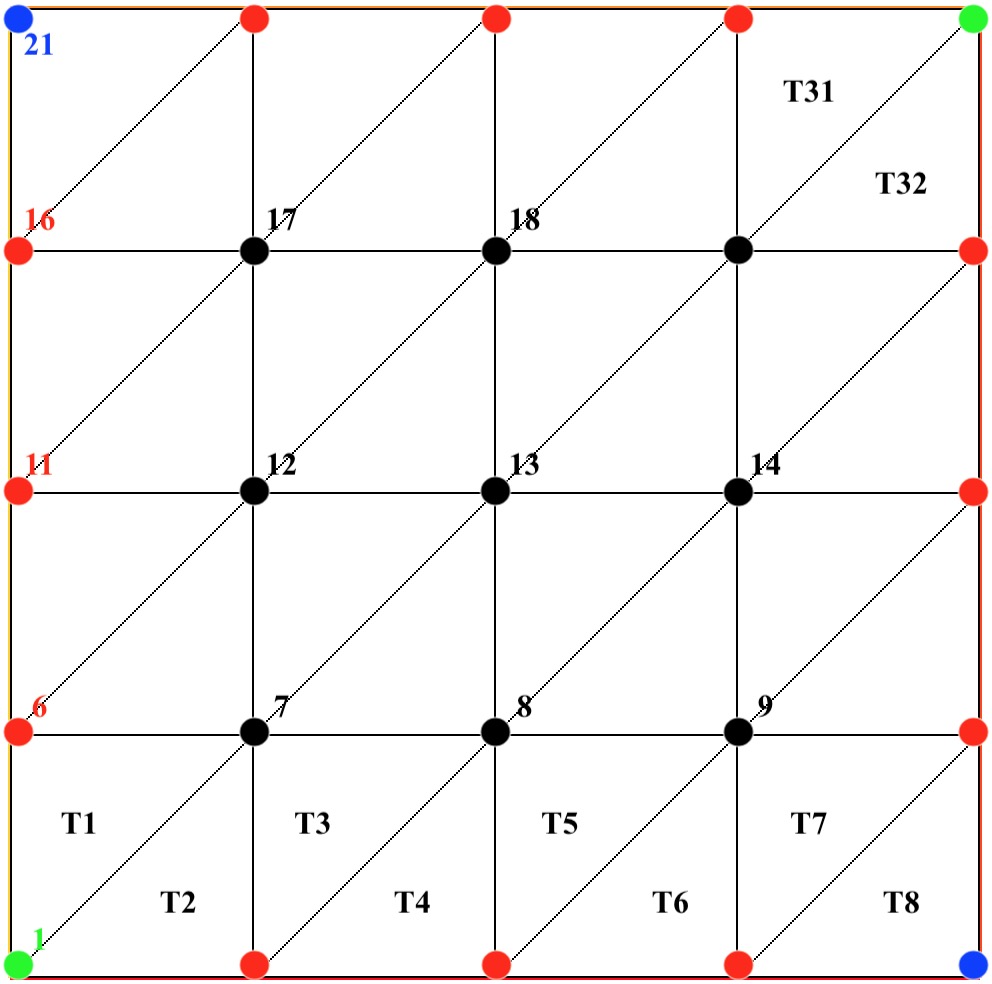}
\caption{\it\small A two-dimension meshing of $\Omega$ with the corresponding nodes and triangle indexing} \label{fig:mesh2}
\end{figure}

To compute the matrices $S(U)$ and $R$, the square domain $\overline{\Omega}$ is partitioned into $n$ equally-spaced nodes in each of the $x$ and $y$ direction, leading to a set of $n^2$ nodes.$$\mathcal{N}=\{P_I(x_i,y_j)|I=1,2,...,N=n^2\}=\p_x\times\p_y$$ The indexing of these nodes starts from left to right, and  bottom to top as shown in Figure \ref{fig:mesh2} for $n=5$. Moreover, the set of $M = 2(n-1)^2$ triangles covering $\overline{\Omega}$ are also indexed from  left to right, and bottom to top $$\mathcal{T}=\{T_J|J=1,2,...,M\},\,\,\overline{\Omega}=\cup_{J}{T_J}$$
The global indexing of the vertices of triangles $T_{2j-1}$ of type a  are $\{j+c,j+n+1+c,j+n+c\}$, whereas that of triangles $T_{2j}$ of type b  are $\{j+c,j+1+c, j+n+1+c\}$, for $j = 1,2,..,\dfrac{M}{2}=(n-1)^2$ and $c = \ceil{\dfrac{j}{n-1}}-1$.
\noindent The matrices $S(U)$ and $R$ are first computed locally on each triangle $T_J$ and then assembled globally. Triangle $T_J$ has three nodes with global indexing $\{\alpha, \beta, \gamma\}$ depending on its type, and local indexing $\{1,2,3\}$. On triangle $T_J$, the only non-zero basis functions are $\psi_{\alpha}$, $\psi_{\beta}$, and $\psi_{\gamma}$  are locally denoted by $\psi_{1}$, $\psi_{2}$, and $\psi_{3}$. Thus, the local $S(U)$  and $R$ matrices on triangle $T_J$ have at most $9$ nonzero entries (in rows and columns $ \alpha, \beta, \gamma$) that can be computed using a $3\times 3$ matrix, denoted by $\mathcal{S}^{(J)}(U^{(J)})$ and $R^{(J)}$, where $U^{(J)} = [U_{\alpha}(t), U_{\beta}(t), U_{\gamma}(t)]$ is a vector of length $3$. 

\noindent By assembling the global matrices and imposing periodic boundary conditions, the degrees of freedom are reduced from $N = n^2$ to $N_1 = (n-1)^2$.
Letting $\{i_1,i_2, \cdots, i_{N_1} \} \subset \{1,2, \cdots, N \}$, one defines the extracted vectors $\tilde{z}\in \mathbb{R}^{N_1}$ from $z\in \mathbb{R}^N$, and extracted matrices $\tilde{E}\in \mathbb{R}^{N_1\times N_1}$ from $E\in \mathbb{R}^{N\times N}$ as shown in the appendices \ref{sec:PSU} and \ref{sec:PR} for the matrices $S(U)$ and $R$. Note that in all following sections we drop the tilde notation, and the matrices $M, K$, $S(U)$ and $R$ are assumed to be of size $N_1 \times N_1$, and the vectors $U(t) \equiv U, W(t) \equiv W$ are of size $N_1$. Moreover, based on this extraction of the minimum number of  degrees of freedom, we consider 
$$w_N(x,y,t)=\sum_{I=1}^{N_1}{W_I(t)\psi_I(x,y)},\qquad \mbox{ and } \quad u_N(x,y,t)=\sum_{J=1}^{N_1}{U_J(t)\psi_J(x,y)} $$
where $\{\psi_J|J=i_1,...,i_{N_1}\}$ is the modified basis extracted from $\{\varphi_I|I=1,...,N\}$.

\noindent In what follows, we compute the local matrices $\mathcal{S}^{(J)}(U^{(J)})$ and $R^{(J)}$ and state the sparsity patterns of the resulting global matrices $S(U)$ and $R$.
\vskip -0.5cm
\subsubsection*{Local Matrix $\mathcal{S}^{(J)}(U^{(J)})$}
The 9 entries $\mathcal{S}^{(J)}_{i,j}$ of the local matrix $\mathcal{S}^{(J)}(U^{(J)})$ are defined as follows for $i,j = 1,2,3$
$$\mathcal{S}^{(J)}(U^{(J)}) =  \begin{bmatrix} \mathcal{S}^{(J)}_{1,1} & \mathcal{S}^{(J)}_{1,2} & \mathcal{S}^{(J)}_{1,3} \vspace{2mm}\\  
\mathcal{S}^{(J)}_{2,1} & \mathcal{S}^{(J)}_{2,2} & \mathcal{S}^{(J)}_{2,3} \vspace{2mm}\\  
\mathcal{S}^{(J)}_{3,1} & \mathcal{S}^{(J)}_{3,2} & \mathcal{S}^{(J)}_{3,3}\vspace{2mm}\\  
\end{bmatrix}
$$ 
$\begin{array}{ccl}
    \mathcal{S}^{(J)}_{i,j} &=& \int_{T_J} \vec{V}(u_n) \cdot \nabla\psi_i \; \psi_j \;\; dA \nonumber \\ &=&  - \int_{T_J} \dfrac{\partial \psi_i }{\partial x} \psi_j \sum_{k=1}^3 \para{U_k^{(J)} \; \dfrac{\partial\psi_k}{\partial y}} \; dA + \int_{T_J} \dfrac{\partial \psi_i }{\partial y} \psi_j \sum_{k=1}^3 \para{U_k^{(J)} \; \dfrac{\partial \psi_k}{\partial x}} \; dA
\end{array}$\vspace{1mm}
\noindent where $\psi_i(x,y) = a_i+b_ix+c_i y$ with 
$$\begin{array}{ccc}
a_1 = \dfrac{x_2 y_3 - x_3 y_2}{2Area(T_J)}, \quad & \quad b_1 = \dfrac{y_2  - y_3}{2Area(T_J)}, \quad & \quad c_1 = \dfrac{x_3 - x_2}{2Area(T_J)}\vspace{2mm}\\
a_2 = \dfrac{  x_3 y_1 - x_1 y_3}{2Area(T_J)}, \quad & \quad b_2 = \dfrac{y_3  - y_1}{2Area(T_J)}, \quad & \quad c_2 = \dfrac{x_1 - x_3}{2Area(T_J)}\vspace{2mm}\\
a_3 = \dfrac{  x_1 y_2-  x_2 y_1}{2Area(T_J)}, \quad & \quad b_3 = \dfrac{y_1  - y_2}{2Area(T_J)}, \quad & \quad c_3 = \dfrac{x_2 - x_1}{2Area(T_J)}\vspace{-1mm}
\end{array}$$
and $Area(T_J) = 0.5x_1(y_2 - y_3) + 0.5x_2(y_3 - y_1)+0.5x_3(y_1-y_2)$. Moreover, $\dfrac{\partial \psi_i }{\partial y} = c_i$, and $\dfrac{\partial \psi_i }{\partial x} = b_i$. Thus,  
\begin{eqnarray}
    \mathcal{S}^{(J)}_{i,j} &=& - \int_{T_J} b_i \psi_j \sum_{k=1}^3 \para{c_k \, U_k^{(J)} } \; dA + \int_{T_J} c_i \psi_j \sum_{k=1}^3 \para{ b_k \,U_k^{(J)} } \; dA\\
    &=& - b_i \para{\sum_{k=1}^3  c_k \, U_k^{(J)} } \int_{T_J}  \psi_j  \; dA + c_i \para{  \sum_{k=1}^3 b_k \,U_k^{(J)} } \int_{T_J}  \psi_j  \; dA
\end{eqnarray}
\noindent Let $b^{(J)} = [b_1, b_2, b_3]$ and $c^{(J)} = [c_1, c_2, c_3]$, then $\sum_{k=1}^3 b_k \,U_k^{(J)} = b^{(J)} \cdot U^{(J)} $ and  $\sum_{k=1}^3 c_k \, U_k^{(J)}  = c^{(J)} \cdot U^{(J)} $ are constants per triangle. Let $\eta_i = \int_{T_J}  \psi_i  \; dA$, and $\eta^{(J)} = [\eta_1,\; \eta_2, \;\eta_3]$ then,\\
\begin{eqnarray} \mathcal{S}^{(J)}(U^{(J)}) &=&  - (c^{(J)} \cdot U^{(J)}) \begin{bmatrix} b_{1} \eta_1   & b_{1} \eta_2 & b_{1} \eta_3 \vspace{2mm}\\  
b_{2}  \eta_1  & b_{2} \eta_2 & b_{2} \eta_3\vspace{2mm}\\  
b_{3}   \eta_1 & b_{3} \eta_2 & b_{3} \eta_3 \vspace{2mm}\\  
\end{bmatrix} + 
(b^{(J)} \cdot U^{(J)}) \begin{bmatrix} c_{1} \eta_1   & c_{1} \eta_2 & v_{1} \eta_3 \vspace{2mm}\\  
c_{2}  \eta_1  & c_{2} \eta_2 & c_{2} \eta_3\vspace{2mm}\\  
c_{3}   \eta_1 & c_{3} \eta_2 & c_{3} \eta_3 \vspace{2mm}\\  
\end{bmatrix}.\end{eqnarray}
\noindent Note that $\eta_i = \int_{T_J}  \psi_i  \; dA = \dfrac{1}{3} Area(T_J)$. Thus, $b_1 \eta_i = \dfrac{1}{6}(y_2 - y_3)$,  $c_1 \eta_i = \dfrac{1}{6}(x_3 - x_2)$, and 
\begin{eqnarray} \mathcal{S}^{(J)}(U^{(J)}) &=&  - \dfrac{1}{6}(c^{(J)} \cdot U^{(J)}) \begin{bmatrix} y_2 - y_3 \vspace{2mm}\\
y_3 - y_1\vspace{2mm}\\  
y_1 - y_2\vspace{2mm}
\end{bmatrix} \begin{bmatrix}1&1&1 \end{bmatrix} + 
\dfrac{1}{6}(b^{(J)} \cdot U^{(J)}) \begin{bmatrix} x_3 - x_2 \vspace{2mm}\\  
x_1 - x_3\vspace{2mm}\\  
x_2 - x_1\vspace{2mm} 
\end{bmatrix}\begin{bmatrix}1&1&1 \end{bmatrix}\\
&=& d_J \;\; \para{ \widehat{c}^{(J)} \cdot U^{(J)} \;\;\;\; \widehat{b}^{(J)} \begin{bmatrix}1&1&1 \end{bmatrix} \;\; - \;\; 
\widehat{b}^{(J)} \cdot U^{(J)}\;\;\;\;  \widehat{c}^{(J)}\begin{bmatrix}1&1&1 \end{bmatrix} }
\end{eqnarray}
where $\widehat{b}^{(j)} =  \begin{bmatrix} y_3 - y_2 \vspace{2mm}\\  
y_1 - y_3\vspace{2mm}\\  
y_2 - y_1\vspace{2mm} 
\end{bmatrix}$, \;\;\; $\widehat{c}^{(j)} =  \begin{bmatrix} x_3 - x_2 \vspace{2mm}\\  
x_1 - x_3\vspace{2mm}\\  
x_2 - x_1\vspace{2mm} 
\end{bmatrix}$, and $d_J =  \dfrac{1}{12Area(T_J)} $\vspace{2mm}

\noindent After computing $\mathcal{S}^{(J)}(U^{(J)})$ its rows and columns are mapped from local indexing $\{1,2,3\}$ to the global indexing $\{ \alpha, \beta, \gamma\}$ and added to the global matrix $S(U)$. Thus, for computing $S(U)$ two difference matrices $B, C$ of size $M\times 3$ have to be computed once and stored, where the $J^{th}$ row of $B$ is $\widehat{b}^{(J)}$ and the $J^{th}$ row of $C$ is $\widehat{c}^{(J)}$; in addition to an $M\times 1$ vector of    triangle areas. Note that the matrix $S(U)$ has to be assembled at each time iteration. \vspace{2mm}

\begin{figure}[H]
    \centering
%
\includegraphics[scale=0.3]{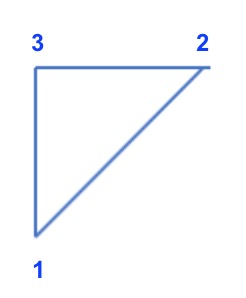} \qquad \qquad \qquad \qquad
\includegraphics[scale=0.3]{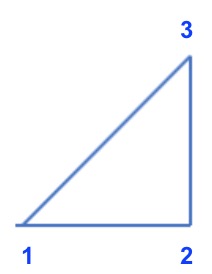}\vspace{-2mm}
    \caption{\it\small The 2 types of triangles with their local nodes' indexing, assuming that the set of nodes on $\Omega$ are equally spaced in the x and y directions. The triangle on the left is denoted by type a, whereas that on the right by type b.} \label{fig:triangles}
\end{figure}

\noindent Assuming that the set of nodes on $\Omega$ are equally spaced, i.e. $x_{i+1} - x_i = y_{i+1} - y_i = h,  \forall i = 0,1,..,n-1 $, then $\mathcal{S}^{(J)}(U^{(J)})$ can be further simplified. In this case, there are 2 types of triangles with the local nodes numbering as shown in Figure \ref{fig:triangles}. 
\vspace{-5mm}
\noindent For triangles of type a, \;$\widehat{b}^{(J)} = h\begin{bmatrix}0&-1&1 \end{bmatrix}^T$ \; and \; $\widehat{c}^{(J)} = h\begin{bmatrix}-1&0&1 \end{bmatrix}^T$. Whereas, for triangles of type b, \;$\widehat{b}^{(J)} = h\begin{bmatrix}1&-1&0 \end{bmatrix}^T$ \; and \; $\widehat{c}^{(J)} = h\begin{bmatrix}0&-1&1 \end{bmatrix}^T$. In both cases, $Area(T_J) = \dfrac{h^2}{2}$ and 
\begin{eqnarray} \mathcal{S}^{(J)}(U^{(J)})
&=& \dfrac{1}{6} \;\; \begin{bmatrix} U^{(J)}_3 - U^{(J)}_2 \vspace{2mm}\\ U^{(J)}_1 - U^{(J)}_3 \vspace{2mm}\\
U^{(J)}_2 - U^{(J)}_1 \vspace{2mm}\\
\end{bmatrix}  \begin{bmatrix}1&1&1 \end{bmatrix}
\end{eqnarray}

\noindent Thus, the computation of $S(U)$ in the case of equally spaced nodes reduces to taking differences of the $U$  vector entries, without the need to store any values.
In addition, $S(U)$ is a block tridiagonal matrix with 2 additional blocks in the upper right  and lower left corner. Moreover, it is a skew-symmetric matrix ($S(U)^T = -S(U)$) that is linear in $U$,  with  6 nonzero entries per row, 6 nonzero entries per column, and zeros on the diagonal assuming the meshing of $\Omega$ shown in Figure \ref{fig:mesh2}. 
$$ S(U) =\dfrac{1}{6} \left[ \begin{array}{cccccc} S_{1,1} & S_{1,2}& 0&\cdots&0&S_{1,k}\\
S_{2,1} & S_{2,2}& S_{2,3}&0&\cdots&0\\
0&\ddots&\ddots&\ddots&\ddots&\vdots\\
\vdots&\ddots&S_{j,l}&S_{j,j}&S_{j,i}&0\\
0&\cdots&0&S_{i,j}&S_{i,i}&S_{i,k}\\
A_{k,1}&0&\cdots&0&S_{k,i}&S_{k,k}\\
\end{array}\right]  \qquad \qquad with \quad S_{i,j} \equiv S_{i,j}(U)
$$
\noindent where $i = n-2, j = n-3, k= n-1, l=n-4$, and the $3(n-1)$ nonzero block matrices $S_{i,j}$ are of size $(n-1)\times (n-1)$ with $2(n-1)$ nonzero entries each, and the following sparsity patterns:
\begin{itemize}
\item $S_{i,i}$ for $i = 1,.., n-1$ are tridiagonal matrices with zero diagonal entries, and nonzero $S_{i,i}(1,n-1),$ and $ S_{i,i}(n-1,1)$. \vspace{-1mm}
\item $S_{1,n-1}$ and $S_{i+1,i}$ for $i = 1,2,3,.., n-2$ are lower bidiagonal matrices, with nonzero entry in first row and column $n-1$.\vspace{-1mm} 
\item  $S_{n-1,1}$ and $S_{i,i+1}$ for $i = 1,2,..,n-2$ are upper bidiagonal matrices with nonzero entry in first column and row $n-1$.
\end{itemize}
\noindent Thus, $S(U)$ has a total of $3(n-1)2(n-1) = 6N_1$ nonzero entries. As for the explicit expressions/values of the entries, refer to appendix \eqref{sec:PSU}.

\subsubsection*{Local Matrix $R^{(J)}$}
The 9 entries $R^{(J)}_{i,j}$ of the local matrix $R^{(J)}$ are defined as follows for $i,j = 1,2,3$
$$R^{(J)} =  \begin{bmatrix} R^{(J)}_{1,1} & R^{(J)}_{1,2} & R^{(J)}_{1,3} \vspace{2mm}\\  
R^{(J)}_{2,1} & R^{(J)}_{2,2} & R^{(J)}_{2,3} \vspace{2mm}\\  
R^{(J)}_{3,1} & R^{(J)}_{3,2} & R^{(J)}_{3,3}\vspace{2mm}\\  
\end{bmatrix}
$$ 
$  R^{(J)}_{i,j} = \int_{T_J}\hat{k}\;\psi_{i,y}\; \psi_j \;\; dA  =   \int_{T_J}\hat{k}\; \dfrac{\partial \psi_i }{\partial y} \psi_j dA$ where $\psi_i(x,y) = a_i+b_ix+c_i y,$ and $\dfrac{\partial \psi_i }{\partial y} = c_i$ with 
$$\begin{array}{ccc}
 c_1 = \dfrac{x_3 - x_2}{2Area(T_J)} \quad & \quad c_2 = \dfrac{x_1 - x_3}{2Area(T_J)}\quad & \quad c_3 = \dfrac{x_2 - x_1}{2Area(T_J)}\vspace{-1mm}
\end{array}$$
and $Area(T_J) = 0.5x_1(y_2 - y_3) + 0.5x_2(y_3 - y_1)+0.5x_3(y_1-y_2)$. Thus, assuming $\hat{k}$ is constant, then  
\begin{eqnarray}
    R^{(J)}_{i,j} &=& \hat{k}\; c_i\;\int_{T_J} \psi_j \; dA =  \hat{k}\; c_i\; \dfrac{1}{3} Area(T_J)\\
R^{(J)}&=& \dfrac{1}{6}\;\hat{k} \begin{bmatrix} x_3 - x_2  & x_3 - x_2  & x_3 - x_2 \vspace{2mm}\\  
x_1 - x_3 & x_1 - x_3 & x_1 - x_3\vspace{2mm}\\  
x_2 - x_1& x_2 - x_1 &x_2 - x_1 \vspace{2mm}\\  
\end{bmatrix} =  \dfrac{1}{6}\;\hat{k}\; \begin{bmatrix}
x_3 - x_2  \\
x_1 - x_3\\
x_2 - x_1
\end{bmatrix} \begin{bmatrix}
1 & 1& 1
\end{bmatrix} = \dfrac{1}{6}\;\hat{k}\; \widehat{c}^{(J)} \begin{bmatrix}
1 & 1& 1
\end{bmatrix} \end{eqnarray}

\noindent After computing $R^{(J)}$, its rows and columns are mapped from local indexing $\{1,2,3\}$ to the global indexing $\{ \alpha, \beta, \gamma\}$ and added to the global matrix $R$. Thus, for computing $R$ one difference matrix $C$ of size $M\times 3$ has to be computed once and stored, where the $J^{th}$ row of $C$ is $\widehat{c}^{(J)}$. Note that the matrix $R$ is computed once. 

\noindent Assuming that the set of nodes on $\Omega$ are equally spaced, i.e. $x_{i+1} - x_i = h,  \forall i = 0,1,..,n-1 $, then $R^{(J)}$ can be further simplified. 

\noindent For triangles of type a shown in Figure \ref{fig:triangles}, $\widehat{c}^{(J)} = h\begin{bmatrix}-1&0&1 \end{bmatrix}^T$ and $R^{(J)}_a =\dfrac{h}{6}\;\hat{k}\; \begin{bmatrix}   - 1 & - 1  &- 1 \vspace{2mm}\\  
0&0&0 \vspace{2mm}\\  
 1 &  1  & 1 \vspace{2mm}\\  
\end{bmatrix}$.\\
 Whereas, for triangles of type b, $\widehat{c}^{(J)} = h\begin{bmatrix}0&-1&1 \end{bmatrix}^T$ and $R^{(J)}_b =\dfrac{h}{6}\;\hat{k}\; \begin{bmatrix}0&0&0 \vspace{2mm}\\  
    - 1 & - 1  &- 1 \vspace{2mm}\\  
 1 &  1  & 1 \vspace{2mm}\\  
\end{bmatrix}$.\\

\noindent Given that the matrix $R$ is independent of $U$, its computation is straightforward. In addition, $R$ has the same block sparsity pattern as $S(U)$. For example, assuming that $\hat{k}$ is constant, then $R$ is a skew-symmetric matrix ($R^T = -R$) with zeros on the diagonal, and 6 nonzero entries per row and 6 nonzero entries per column, of the form $\dfrac{h}{6}\hat{k} \alpha$ where $\alpha = -2,-1,1,\; or \; 2$, with the sum of entries per row or column being zero. (refer to appendix \eqref{sec:PR}).
$$R =\dfrac{h}{6}\hat{k} \left[ \begin{array}{cccccc} R_{1,1} & R_{1,2}& 0&\cdots&0&R_{1,n-1}\\
R_{2,1} & R_{2,2}& R_{2,3}&0&\cdots&0\\
0&\ddots&\ddots&\ddots&\ddots&\vdots\\
\vdots&\ddots&R_{j,m}&R_{j,j}&R_{j,i}&0\\
0&\cdots&0&R_{i,j}&R_{i,i}&R_{i,l}\\
A_{l,1}&0&\cdots&0&R_{l,i}&R_{l,l}\\
\end{array}\right]
$$
\noindent where $i = n-2, j = n-3, l= n-1, m=n-4$, and the $3(n-1)$ nonzero block matrices $R_{i,j}$ are of size $(n-1)\times (n-1)$ with $2(n-1)$ nonzero entries each, and the following sparsity patterns:
\begin{itemize}
\item $R_{i,i}$ for $i = 1,2,.., n-1$ are such that: $R_{i,i}(j,j+1)=1$, $R_{i,i}(j+1,j)=-1$, for $j=1,2,..., n-2$, $R_{i,i}(1,n-1) = -1$, and $R_{i,i}(n-1,1) = 1$.
\item $R_{1,n-1} = R_{i+1,i}$ for $i = 1,2,.., n-2$ are lower bidiagonal matrices ($R_{i+1,i}(j,j) =2, R_{i+1,i}(j+1,j) =1$), with $R_{i+1,i}(1,n-1) = 1$.
\item  $R_{n-1,1} = R_{i,i+1}$ for $i = 1,.., n-2$ are upper bidiagonal matrices ($R_{i,i+1}(j,j) =-2, R_{i,i+1}(j,j+1) =-1$), with $R_{i,i+1}(n-1,1) = -1$.
\end{itemize}
\noindent Thus, $R$ has a total of $6N_1$ nonzero entries. \vspace{2mm}
\subsection{Solution of the Semi-Linear Scheme \eqref{HMC-Disc-Comp-semi}}\label{sec:4.1}
The semi linear scheme \eqref{HMC-Disc-Comp-semi} can be solved at each time step as shown in algorithm \eqref{alg:HMC-Newton}, where the matrices need to be generated as described previously, based on the meshing of the space domain given in Figure \ref{fig:mesh2}. For that purpose, we implemented Algorithm \ref{alg:HMC-Newton} using FreeFem++.
\begin{algorithm}[H]
\centering
\caption{ Numerical Hasegawa-Mima semi-linearized Finite Element Scheme}
{\renewcommand{\arraystretch}{1.3}
\begin{algorithmic}[1]
\Statex{\textbf{Input:} $M$:  mass matrix ; $K = M+A$, $A$: stiffness matrix; $S(U)$: algorithm that builds  $S(U)$ as in Appendix \ref{sec:SU};}
\Statex{ \qquad \quad $R$: matrix defined in Section \ref{sec:SR}; $U0, W0$: the discrete initial condition vectors; $\tau$: time step; $T$: end time;}
\State $U = U0$; $W = W0$; \vspace{1mm}
\For {$t = 0:\tau:T-\tau$\vspace{2mm}} 
\State $G=M*W$;\vspace{1mm}
\State Solve for $W$: $(M+ \tau\,S(U))\;W=\tau\, R\,U + G$;\vspace{1mm}
\State Solve for $U$: $K\,U= M\,W$;\vspace{1mm}
\EndFor
\end{algorithmic}}
\label{alg:HMC-Newton}
\end{algorithm}
\noindent We consider a square domain $[x_0,x_n]\times [y_0,y_n]$ with a uniform mesh in the $x$ and $y$ direction ($x_i - x_{i-1} = \dfrac{x_n - x_0}{n} = y_i - y_{i-1} = \dfrac{y_n - y_0}{n}$ for $i = 1,2,..,n$ and $n$ intervals in the $x$ and $y$ directions respectively) and the finite element $\mathbb{P}1$ space with periodic boundary conditions, using appropriate Freefem++ functions. 
The function $p$ of the initial Hasegawa-Mima PDE is defined, where in most simulations it is assumed that $p_y = 0$, and  $\hat{k} = p_x$ is a constant, unless stated otherwise. The initial conditions $u_0$ is  given as input. As for the initial condition $w_0 = u_0 - \Delta u_0$ it could be given as input if $u_0$ is a simple function. However, for any function $u_0$, we compute the vector $W_0 = W(0)$ by solving the linear system $${M}* W_0 = {K}*U_0$$ where the vectors $U_0 = U(0), W_0$,   the mass matrix ${M}$, and the matrix ${K}={M}+{A}$ are  defined in \eqref{HMC-Disc-V1-1}, with $A$ being the stiffness matrix. Note that the matrices ${M}$, ${A}$, ${R}$ and ${S}({U}^j)$ are generated in Freefem++ using the corresponding variational formulations:
\begin{eqnarray}
 a(u,v) &=& \int\limits_{Th} (u_x*v_x + u_y*v_y); \mbox{ where the matrix } {A} = a(Vh,Vh);\nonumber \\
  b(u,v) &=& \int\limits_{Th} (u*v); \mbox{ where the matrix }{M} = b(Vh,Vh);\nonumber\\
  c(u,v) &=& \int\limits_{Th} (p_x*u_y*v - p_y*u_x*v);  \mbox{ where the matrix }{R} = c(Vh,Vh); \nonumber\\
  d(w,v) &=& \int\limits_{Th}(u^j_x*w_y*v - u^j_y*w_x*v ) \mbox{ where the matrix } {S}(\tilde{U}^j) = d(Vh,Vh).
\end{eqnarray}

\noindent Another alternative in Freefem++ is to just define the variational formulation of \eqref{HMC-Disc-Comp-semi} as problems that are solved at each time iteration:  
\begin{eqnarray} 
 hypo(w^{j+1},v) &=& \int\limits_{Th} (w^{j+1}*v/\tau - w^j*v/\tau)+\int\limits_{Th}(u^j_x*w^{j+1}_y*v-u^j_y*w^{j+1}_x*v)  +\int\limits_{Th}(p_y*u^j_x*v -p_x*u^j_y*v )\nonumber\\
  ellip(u^{j+1},v) &=& \int\limits_{Th}(u^{j+1}_x*v_x + u^{j+1}_y*v_y)+\int\limits_{Th}(u^{j+1}*v)-\int\limits_{Th}(w^{j+1}*v) \nonumber
\end{eqnarray}
where $u^{j+1}, w^{j+1}$ refer to  the sought solutions at the current time iteration, i.e. $u(t+\tau), w(t+\tau)$, and $u^j, w^j$ refer to $u(t), w(t)$ the solutions at the previous time iteration. 
At each time iteration, Freefem++ will generate the corresponding matrices and vectors and solve the linear systems accordingly. However, this implies that the fixed matrices $M,K,$ and $R$ will be regenerated redundantly at each iteration. Thus, it is preferable timewise to generate the fixed matrices once in the algorithm and solve the matrix form of the problem (Algorithm \ ref{alg:HMC-Newton}). Table \ref{tab:time} validates this claim, where the execution time is reduced at least by half.

\noindent Note that the simulation is stopped once the maximum value of $u(t)$ at one of the mesh nodes is $0.3$, which corresponds to the maximum value attained physically.
\begin{table}[H]
\centering
\begin{tabular}{|c|r||r|r|}
\hline
\multirow{2}{*}{$\tau$}&\multirow{2}{*}{$T$}& FreeFem++& FreeFem++ \\
&& Variational Form& Matrix Form\\

\hline
1/8&34.2500&51.8857&25.6141\\
1/10&40.7000&78.3870&35.8261\\
1/16&61.3125&188.9680&85.7979\\
1/32&119.0000&676.0540&318.1010\\ 
1/64&235.0000&3245.3400&1272.4000\\ \hline
\end{tabular}
\caption{\small Execution time of the semi linear scheme \eqref{HMC-Disc-Comp-semi} in FreeFem++ using the Variational Form or the Matrix Form for $n=64$, $u_0 = 10^{-5}sin(3x)$, $\Omega = [0,\pi]\times[0,\pi]$, $\hat{k} = p_x = 12$, $p_y = 0$, and a given time step $\tau$. $T$ denotes the end time at which the simulation was stopped as $u(t)$ reached the maximum value of $0.3$ at one of the mesh nodes.} \label{tab:time}
\end{table}
\subsection{Testing}\label{Tests}
\begin{figure}[H]
\centering
\includegraphics[scale=0.17]{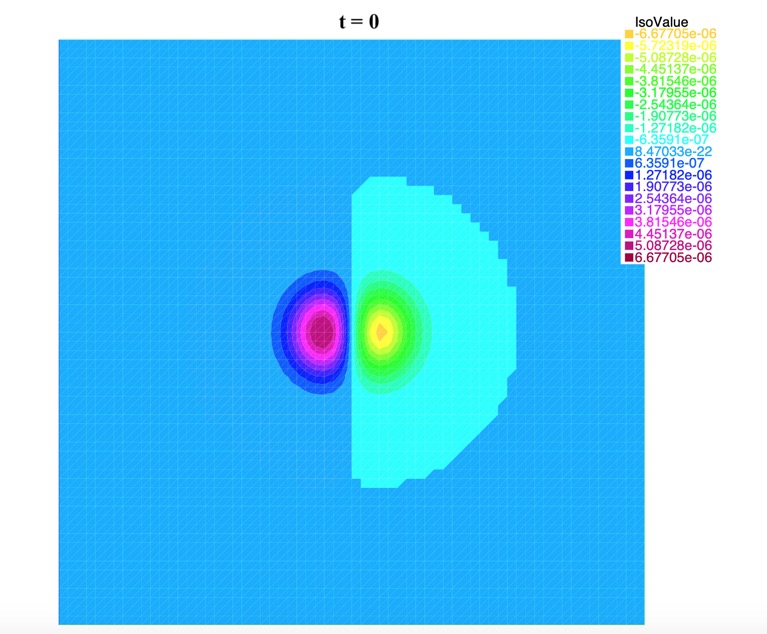}
\includegraphics[scale=0.17]{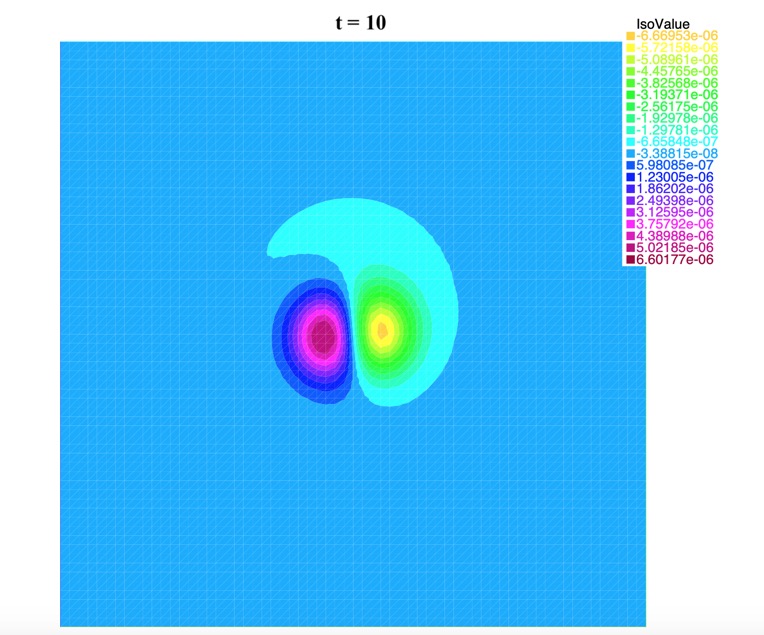}
\includegraphics[scale=0.17]{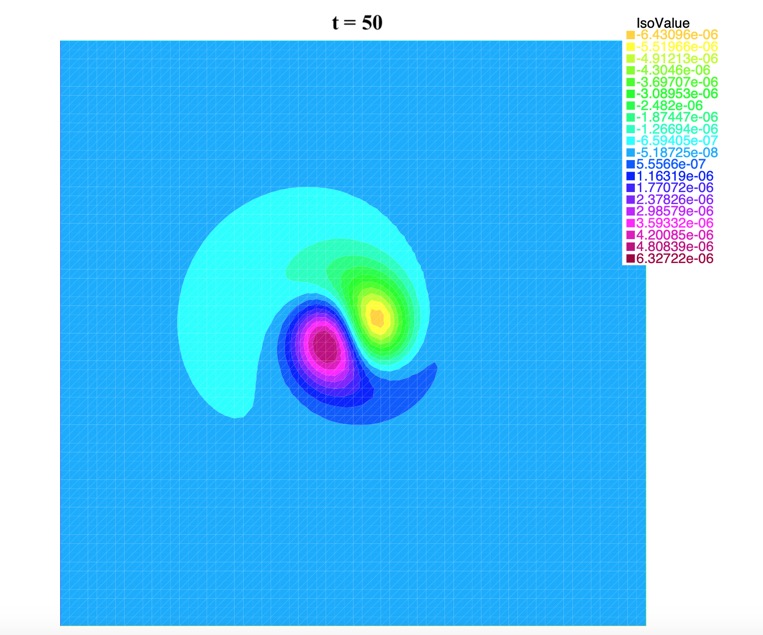}
\includegraphics[scale=0.17]{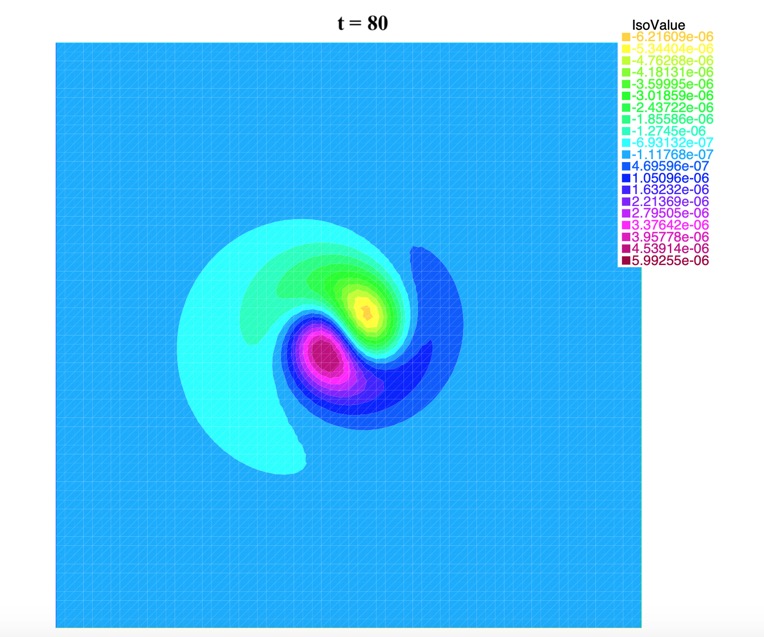}
\includegraphics[scale=0.17]{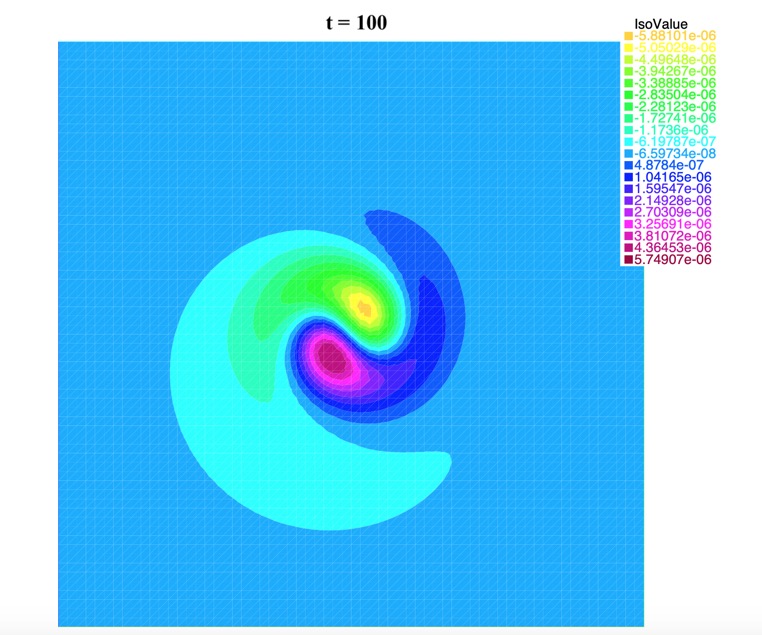}
\includegraphics[scale=0.17]{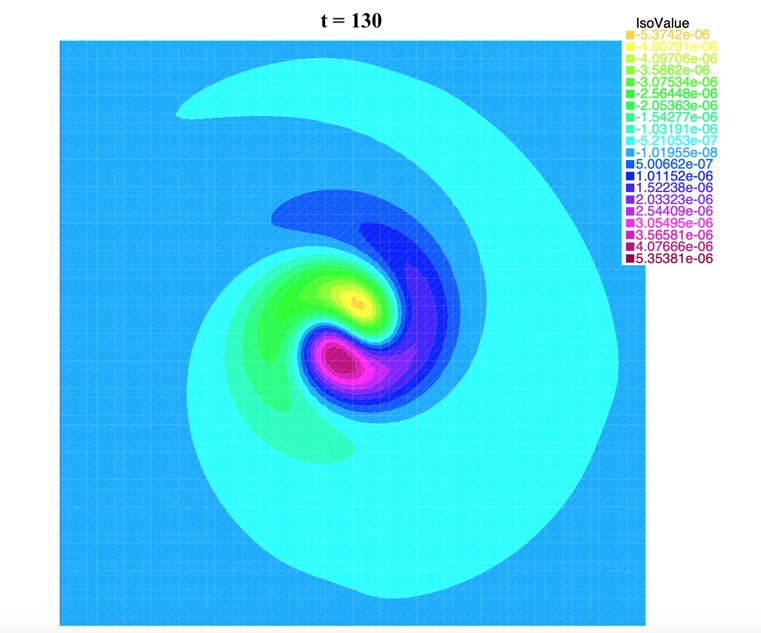}
\includegraphics[scale=0.17]{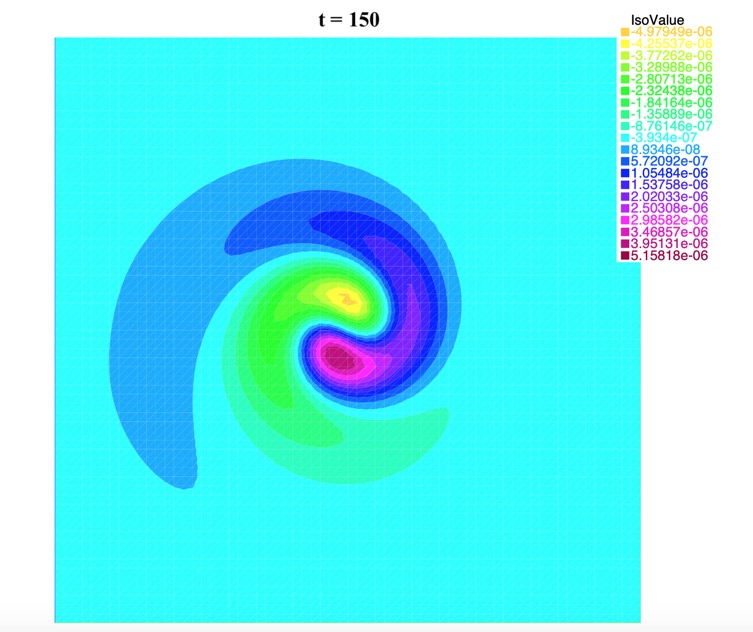}
\includegraphics[scale=0.17]{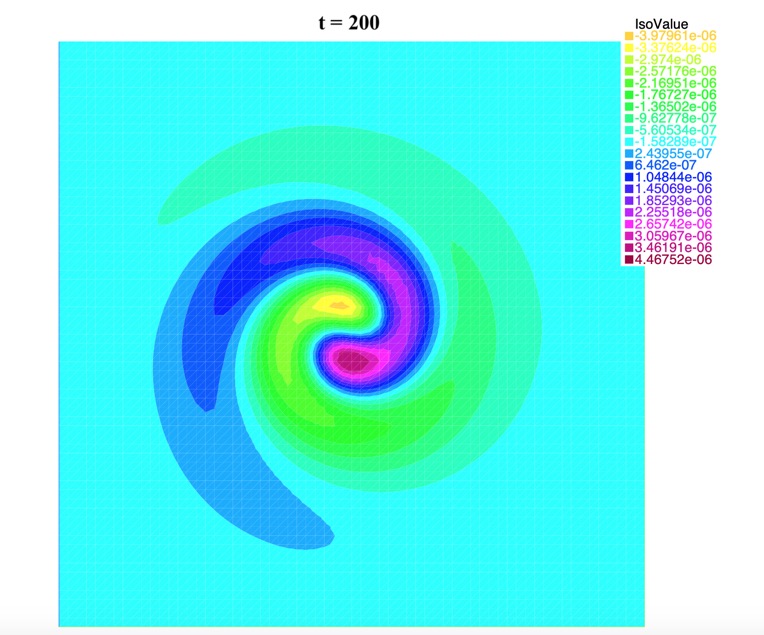}
\includegraphics[scale=0.17]{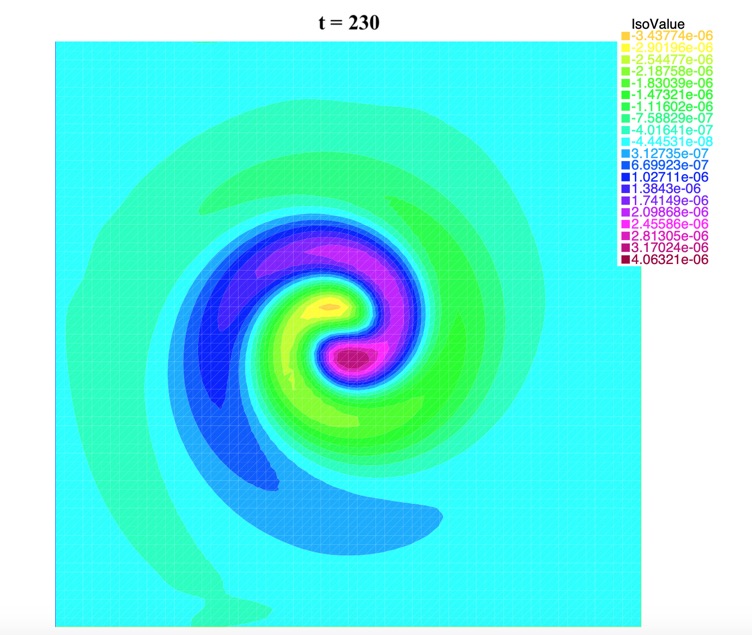}
\includegraphics[scale=0.17]{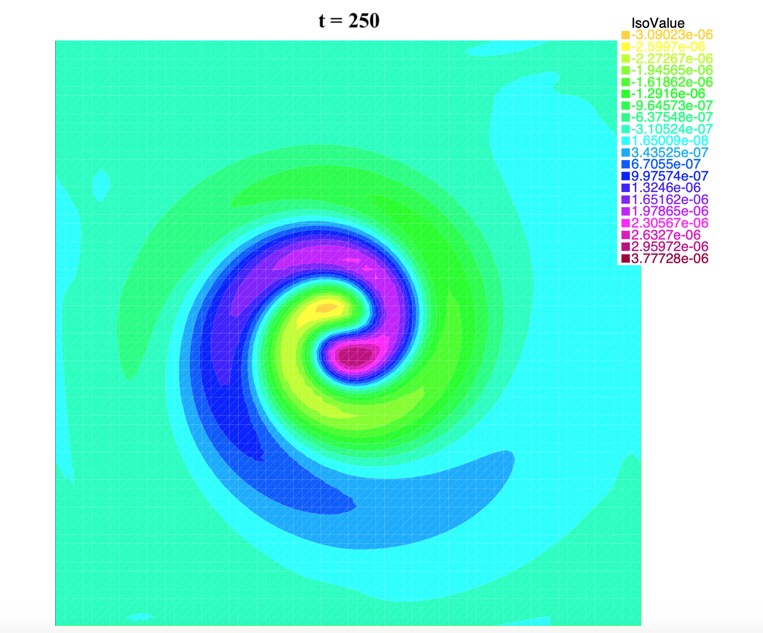} 
\includegraphics[scale=0.175]{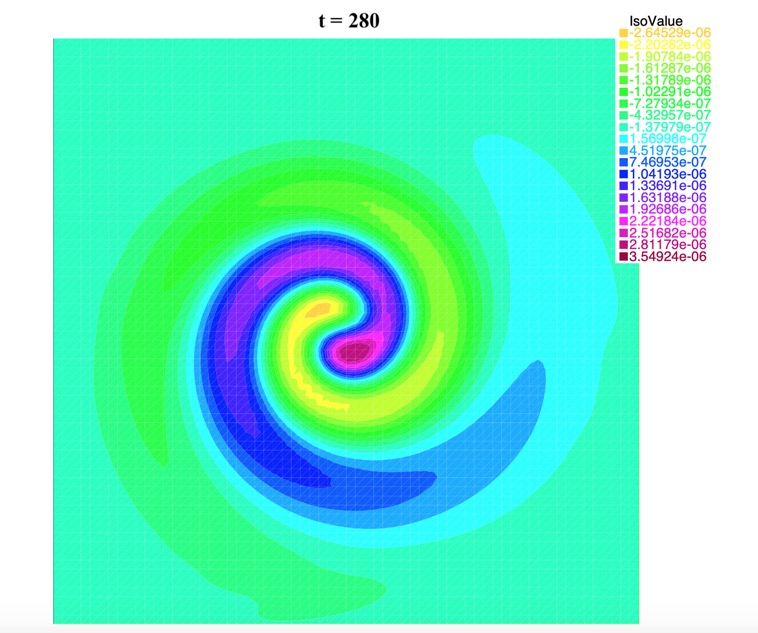} 
\includegraphics[scale=0.17]{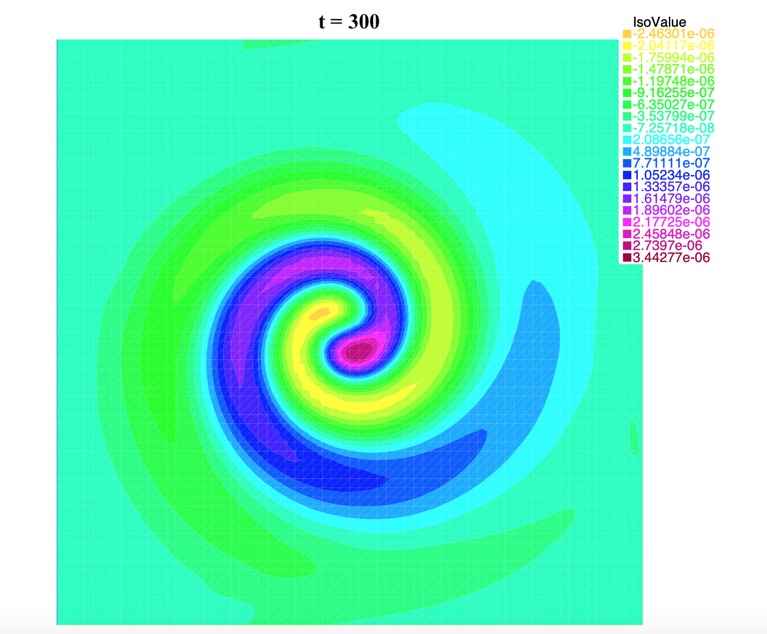} 
\includegraphics[scale=0.17]{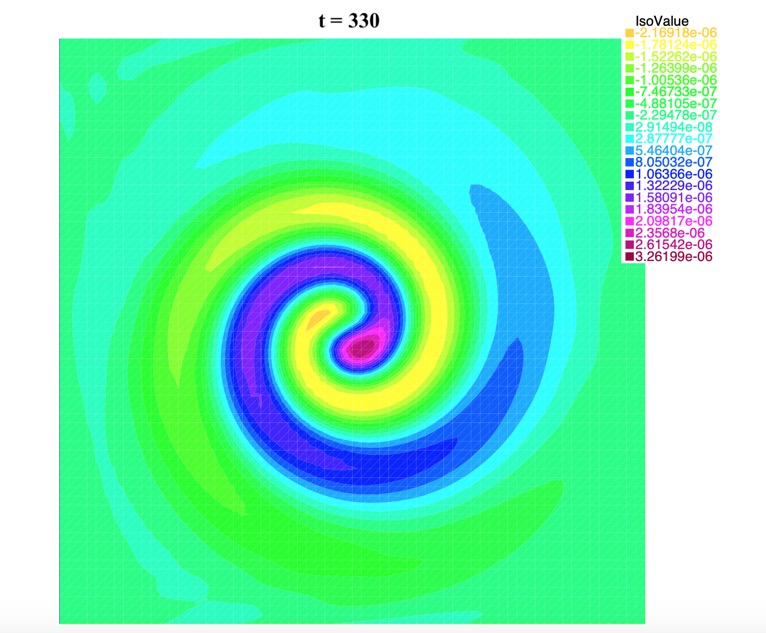} 
\includegraphics[scale=0.17]{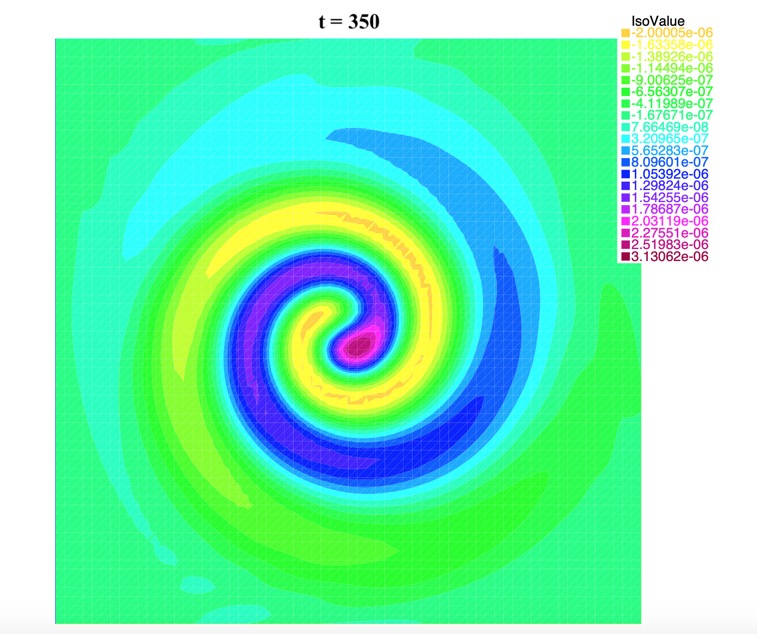} 
\includegraphics[scale=0.17]{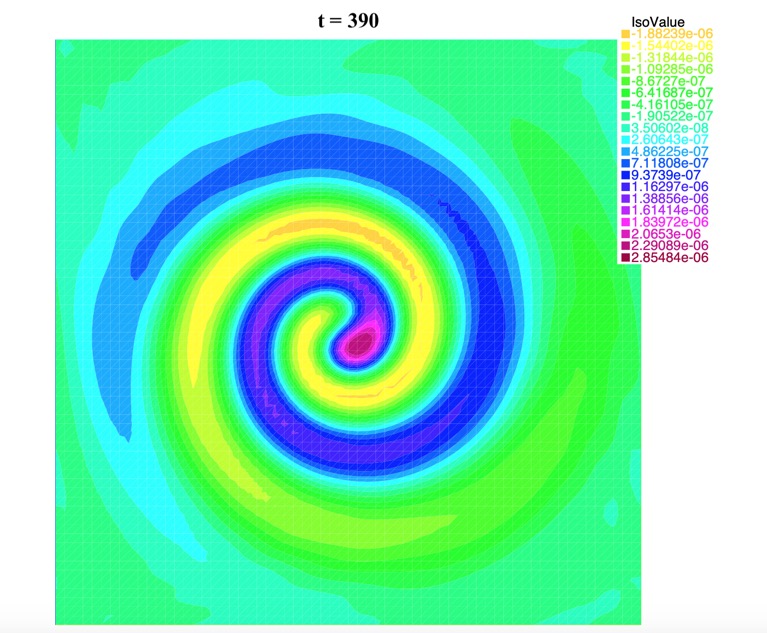} 
\caption{\it \small Time evolution of solution $u$ for $u_0= -10^{-5}(x-10)e^{ -0.5(x-10)^2 - 0.5(y-10)^2}$, $\tau = 0.1$, and a $65 \times 65$ grid on  $\Omega = [0,20] \times [0,20]$.}\label{fig:gauss-64}
\end{figure}

We start by testing Algorithm \ref{alg:HMC-Newton} for $\dfrac{n_0}{\omega_{ci}} = e^{Ax+B}$, i.e. $p = Ax+B$ where $\hat{k} = p_x = A$ and $p_y = 0$. As explained in \cite{hariri}, the solution is expected to be a traveling wave in the y-direction for a nonzero $A$. The speed of the motion and its direction depend on the magnitude and sign of $A$ respectively. We consider different initial conditions for the same exponential density profile $n_0$. We consider the following cases:
\begin{enumerate}
\item Domain $\Omega = [0,1] \times [0,1]$ with the number of intervals in  the $x$ and $y$ direction $n=64$ ($h = 1/64 \approx 0.015625 $), $A=12$, the initial condition $u_0(x,y) = 10^{-5}sin(10\pi y)$ and $\tau = 0.1$. Figure \ref{fig:sin9piy} shows the time evolution of the solution of the Semi-Linear Scheme \eqref{HMC-Disc-Comp-semi} at time $t = 0,5,10,15,20$. It is clear that the $sin$ function is moving in the y-direction as time proceeds without any perturbation in its initial form up till $t = 200$.  Afterwards, the solution grows with time to reach $||u||_{\infty} = 0.3$ at $t = 260.4$ when the algorithm is stopped.\\
 Note that even though $\tau > \dfrac{1}{2||\hat{k}||_{\infty}} = \dfrac{1}{24}$ does not satisfy the sufficient condition for proving the existence and uniqueness of a solution to  \eqref{HMC-Disc-Comp}, however the algorithm still converges and produces the expected behavior. \vspace{4mm}
\begin{figure}[H]
\begin{tabular}{cc}
\subfloat{\includegraphics[scale=0.15]{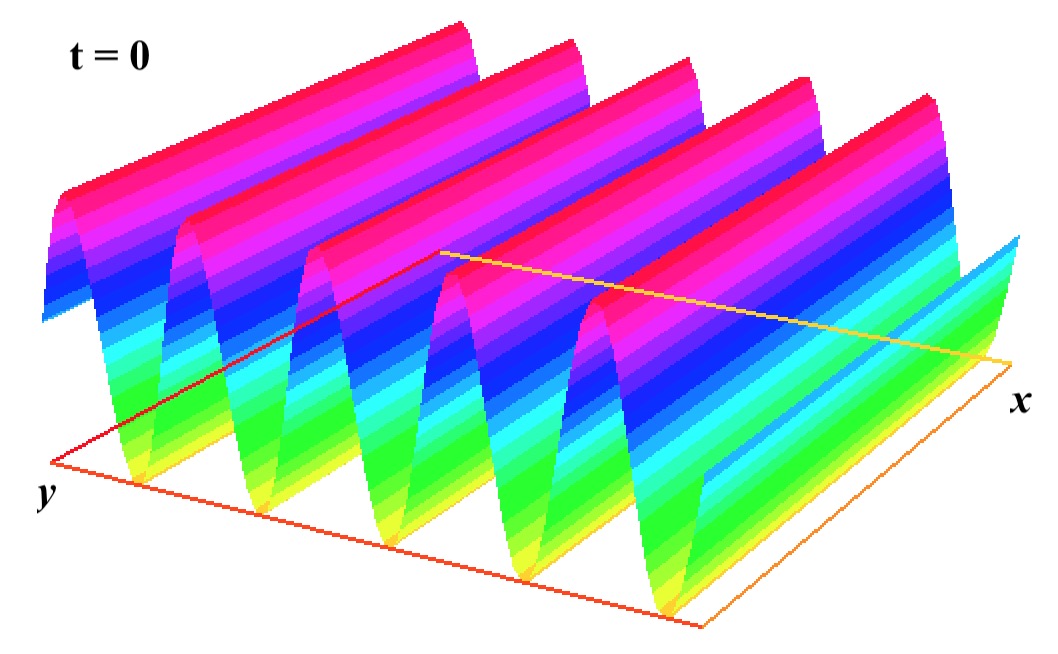}} &
\subfloat{\includegraphics[scale=0.15]{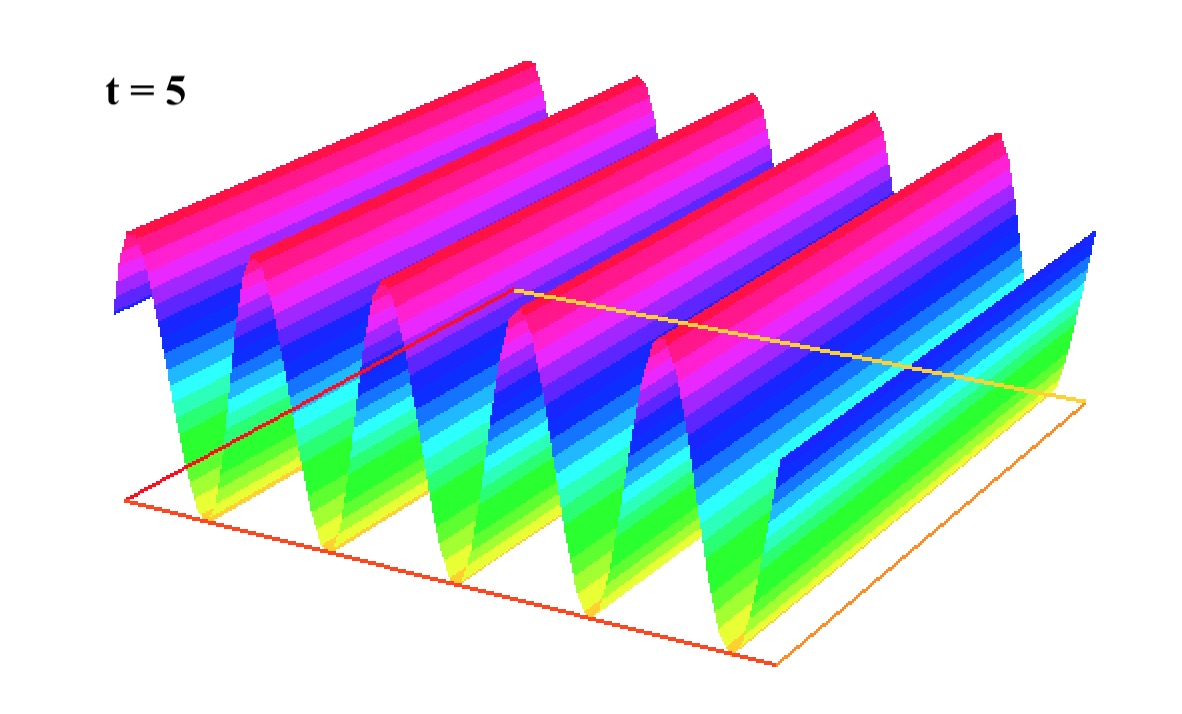}} \\
\subfloat{\includegraphics[scale=0.15]{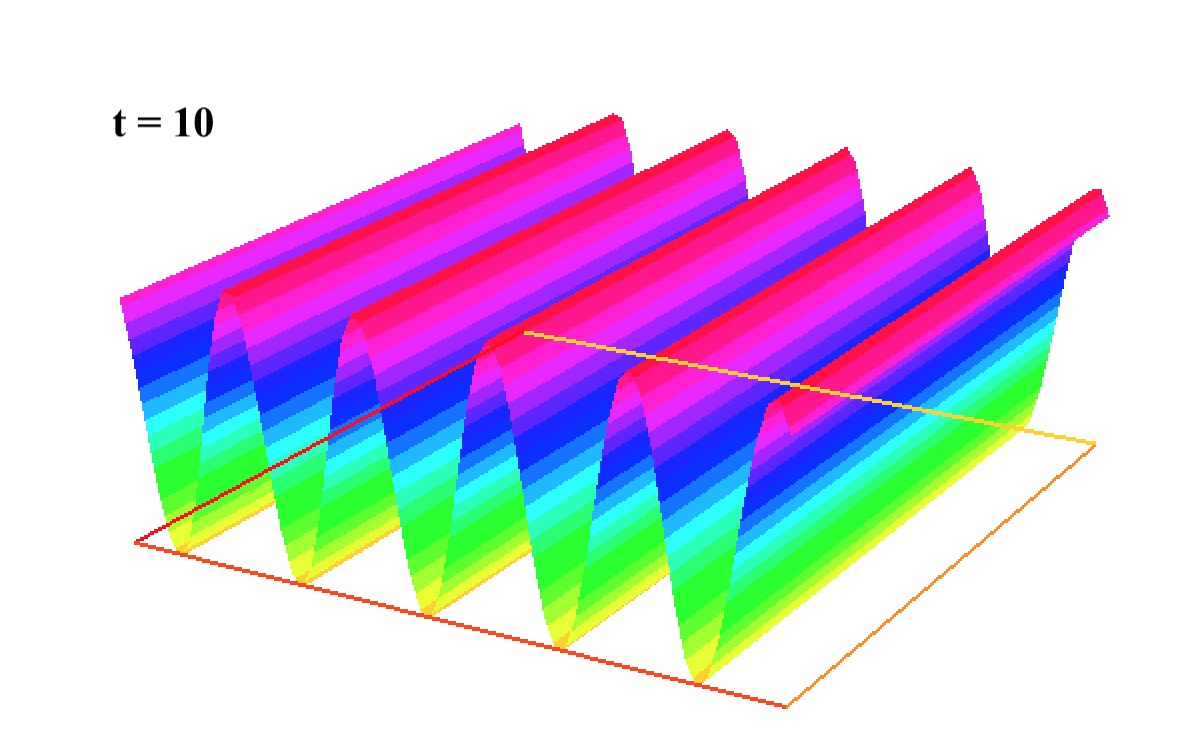}} &
\subfloat{\includegraphics[scale=0.15]{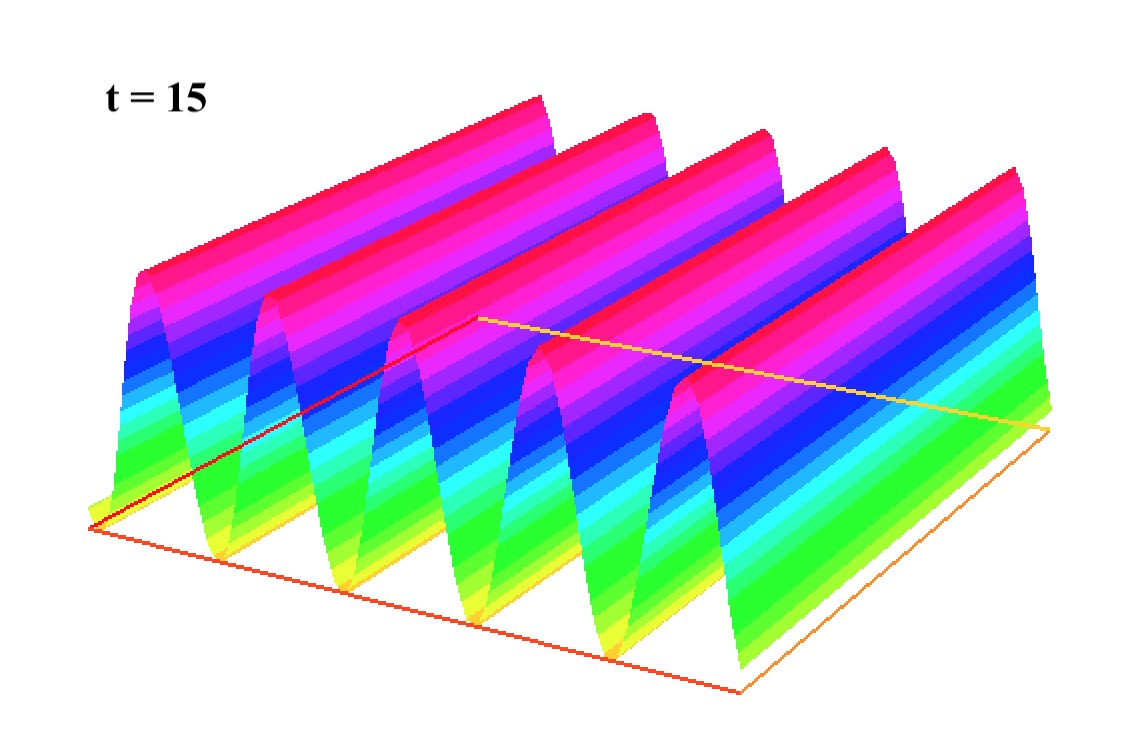}} \\
\subfloat{\includegraphics[scale=0.15]{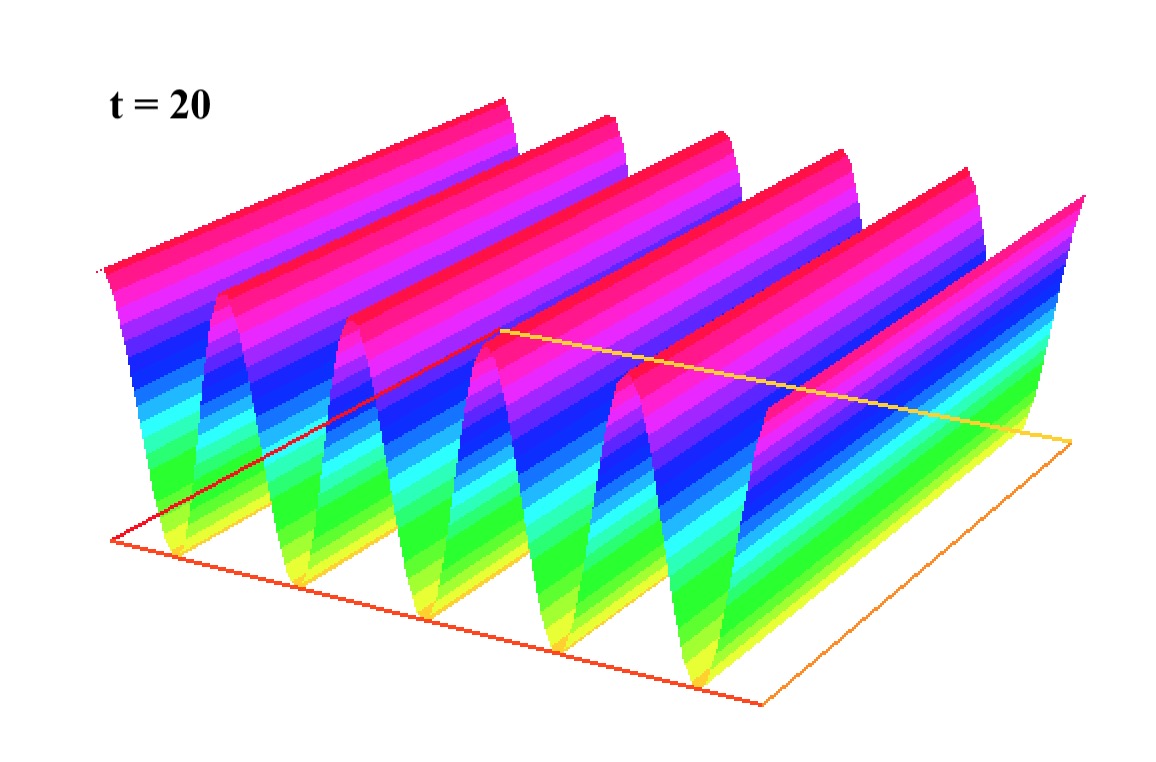}} &
\subfloat{\includegraphics[scale=0.2]{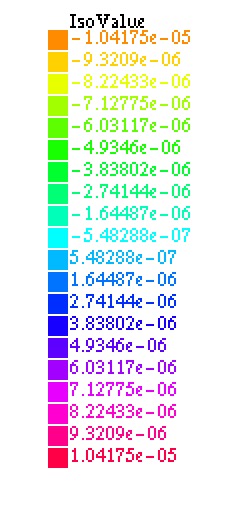}} \\
\end{tabular}
\centering
\caption{\it \small Time evolution of solution $u$ for $u_0= 10^{-5}sin(10\pi y)$, $\tau = 0.1$, and a $65 \times 65$ grid on  $\Omega = [0,1] \times [0,1]$.}\label{fig:sin9piy}
\end{figure}\newpage
\item Domain $\Omega = [0,\pi] \times [0,\pi]$ with $A=12$, the initial condition $u_0(x,y) = 10^{-5}sin(3 y)$, and $\tau = 0.1$. \\In Figure \ref{fig:sin3y32} we consider the number of intervals in  the $x$ and $y$ direction $n=32$ ($h = \pi/32 \approx 0.098175$), and in Figure \ref{fig:sin3y64} $n=64$ ($h = \pi/63 \approx 0.049087$). We notice the same behavior where the solution moves in y-direction and grows with time at a faster rate to reach $||u||_{\infty} = 0.3$ at $t = 9.6$. The larger $h$ value ($h \approx 0.1, n=32$) does not affect the solution with respect to that of $h \approx 0.05, n=64$, apart from the smoothness of the 3D surface. From this perspective, this shows the robustness of the algorithm for reasonable $h<1$ values.
\begin{figure}[H]
\centering
\includegraphics[scale=0.18]{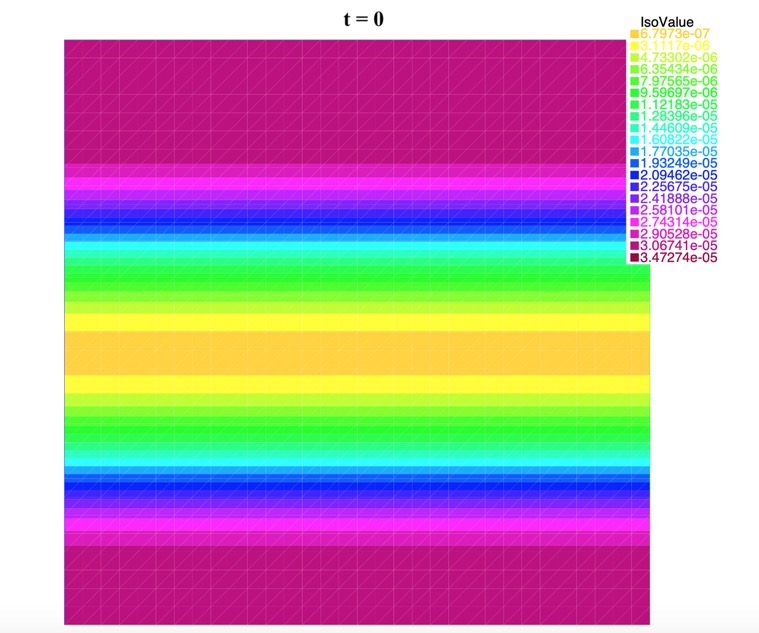}
\includegraphics[scale=0.18]{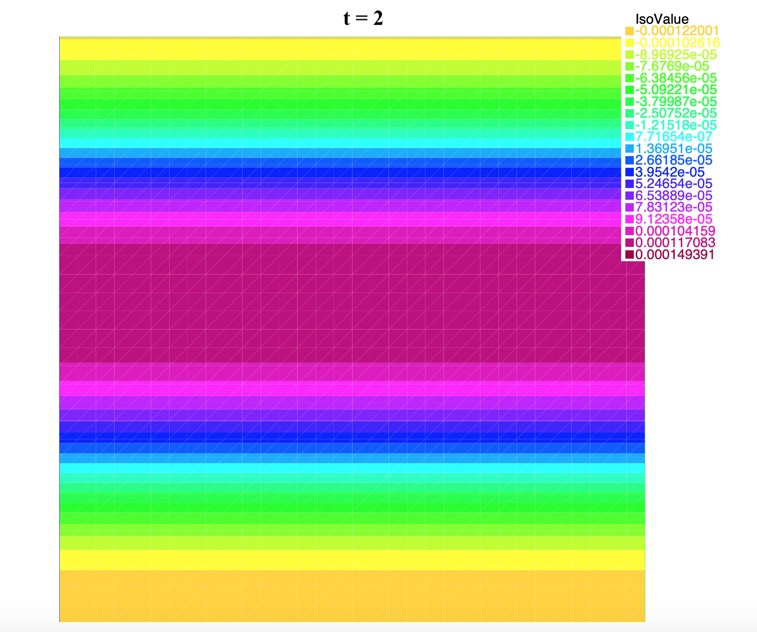}
\includegraphics[scale=0.18]{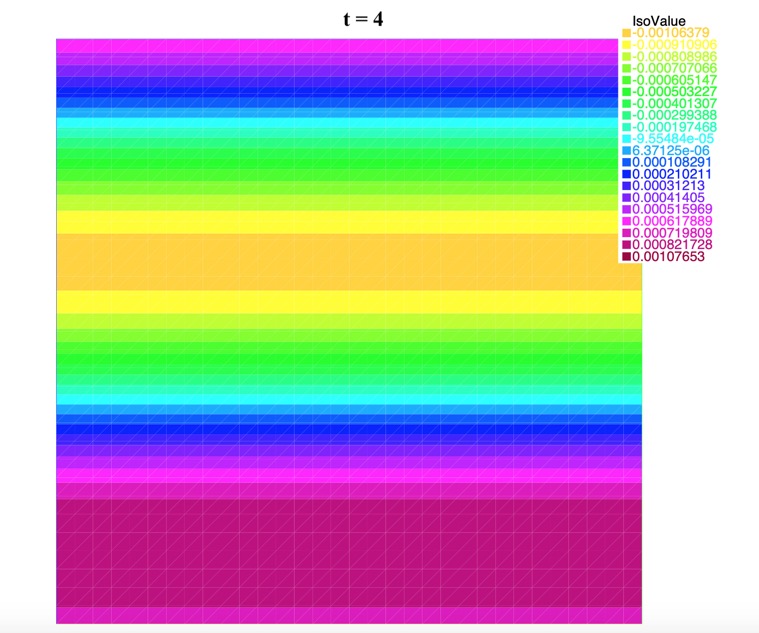}\\ \vspace{2mm}
\includegraphics[scale=0.18]{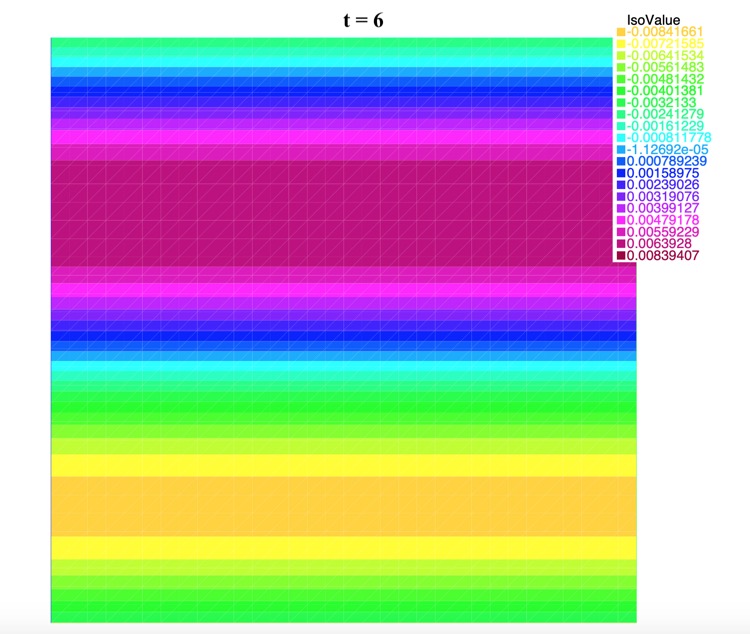}
\includegraphics[scale=0.185]{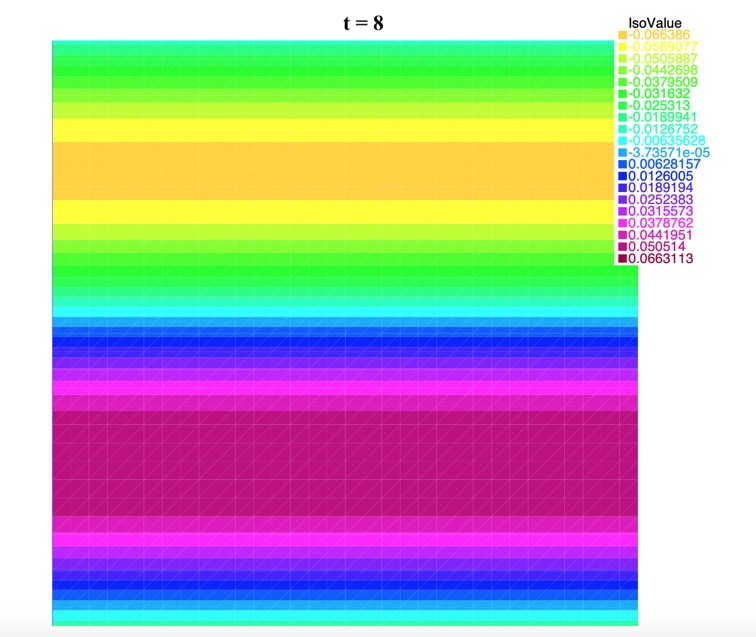}
\includegraphics[scale=0.18]{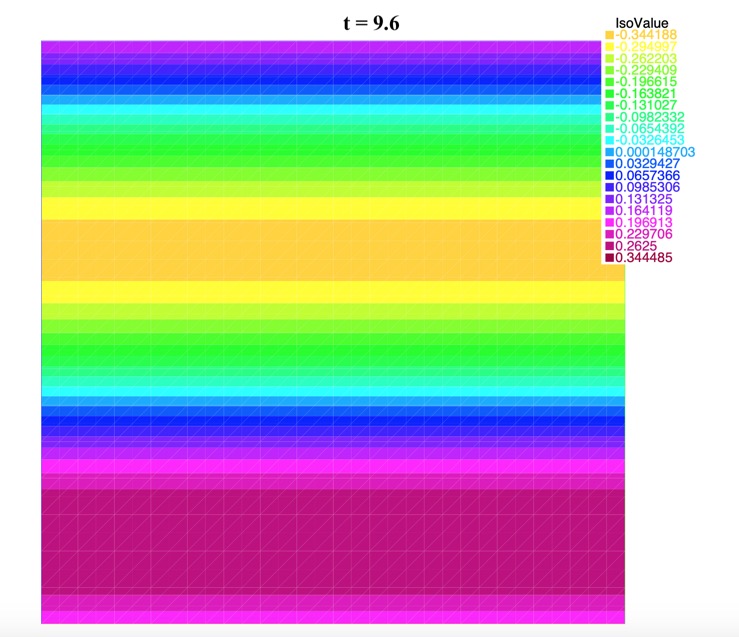} 
\caption{\it \small Time evolution of solution $u$ for $u_0= 10^{-5}sin(3 y)$, $\tau = 0.1$, and a $33 \times 33$ grid on  $\Omega = [0,\pi] \times [0,\pi]$.}\label{fig:sin3y32}
\end{figure}
 \begin{figure}[H]
\centering
\includegraphics[scale=0.18]{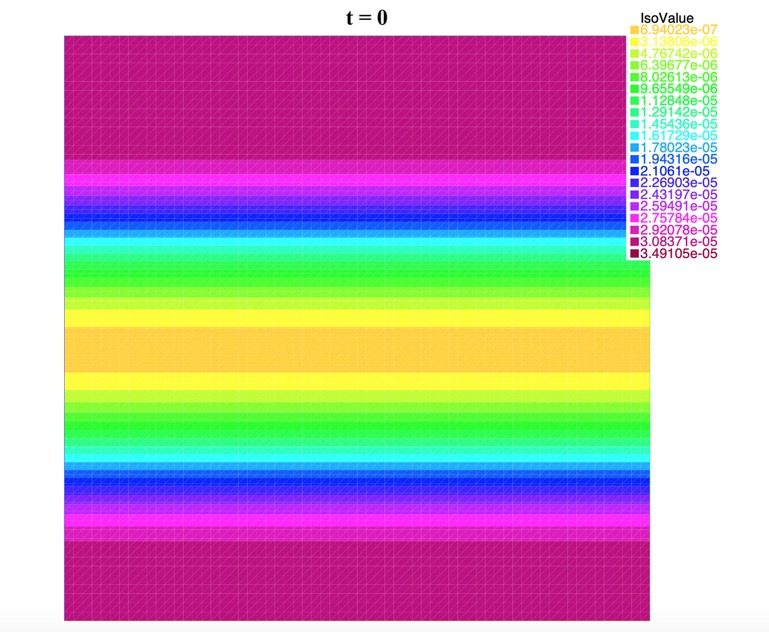}
\includegraphics[scale=0.18]{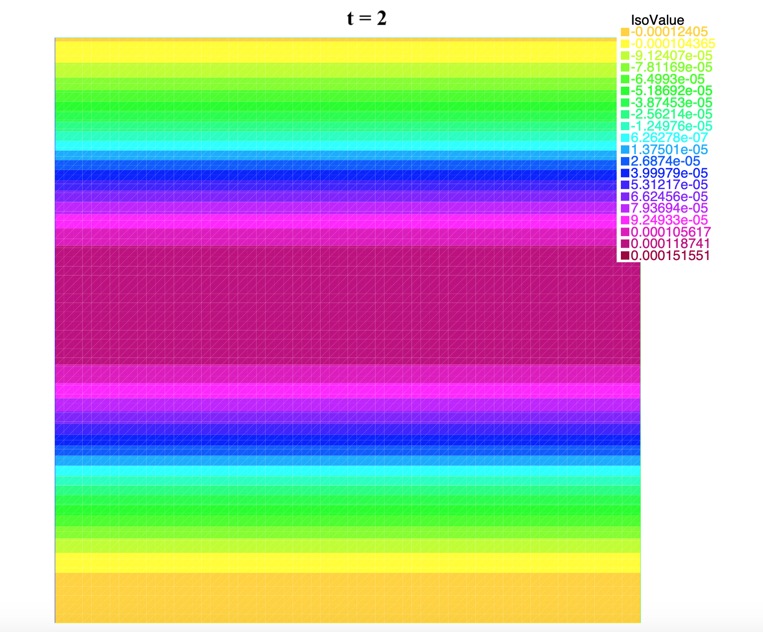}
\includegraphics[scale=0.18]{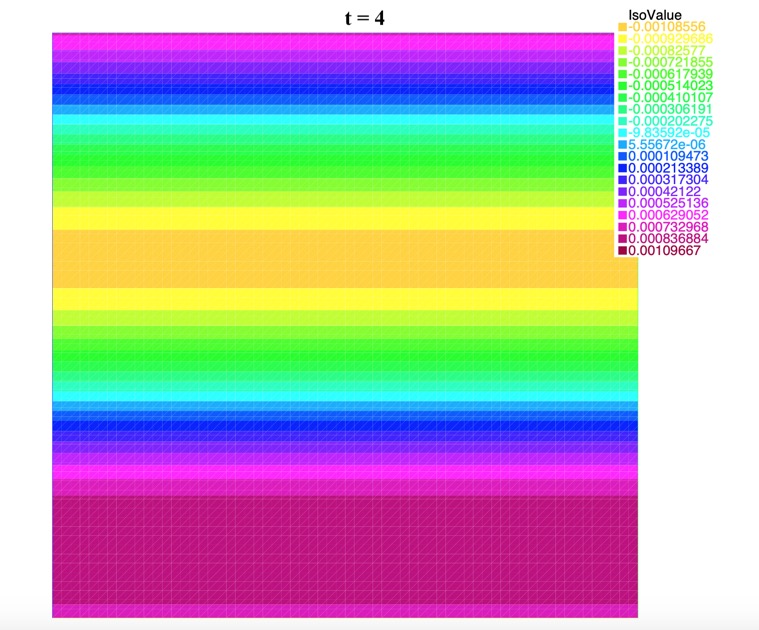}\\ \vspace{2mm}
\hspace{1mm}\includegraphics[scale=0.18]{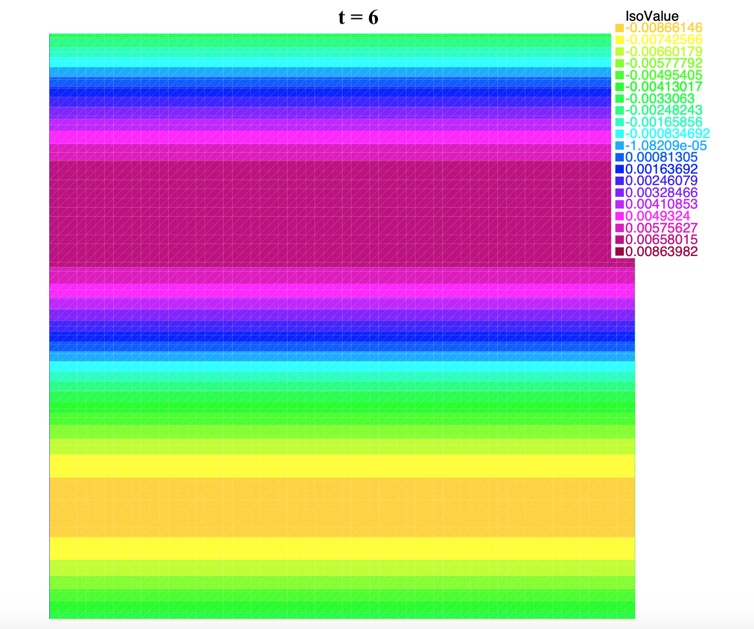}
\includegraphics[scale=0.18]{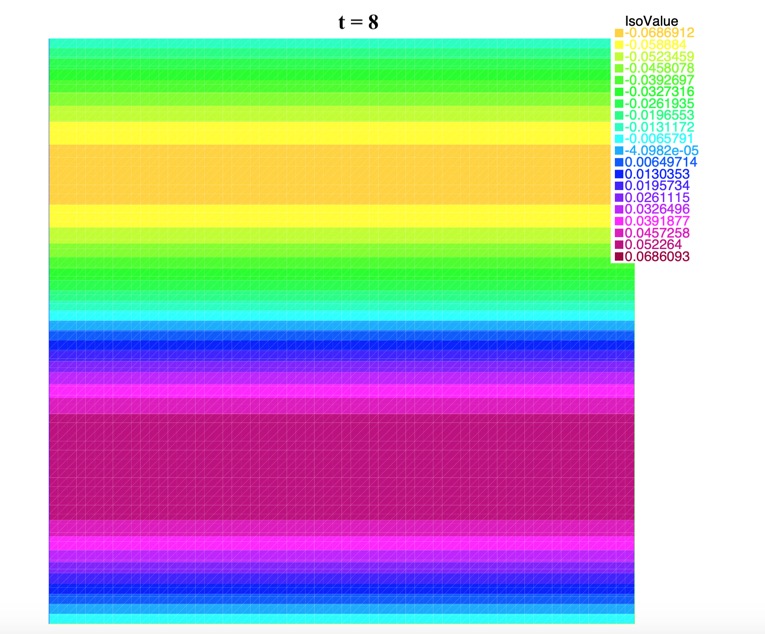}
\includegraphics[scale=0.18]{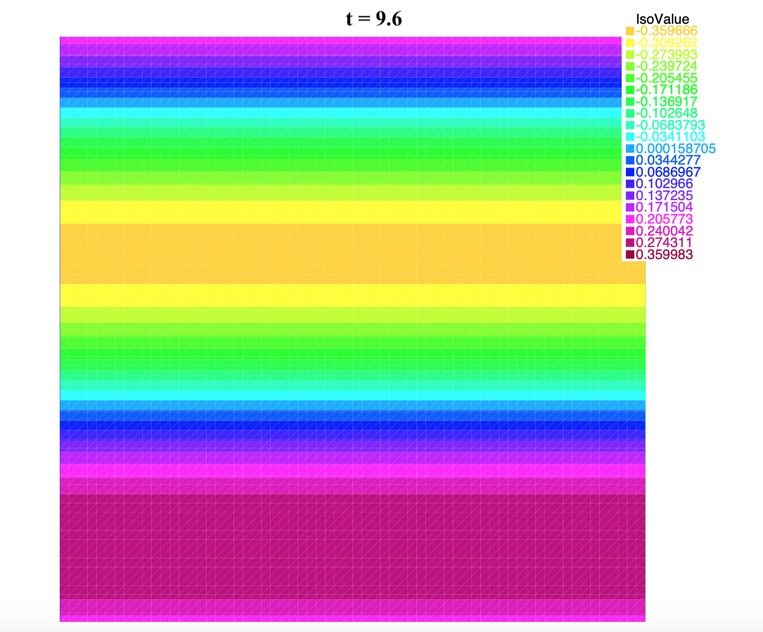} 
\caption{\it \small Time evolution of solution $u$ for $u_0= 10^{-5}sin(3 y)$, $\tau = 0.1$, and a $65 \times 65$ grid on  $\Omega = [0,\pi] \times [0,\pi]$.}\label{fig:sin3y64}
\end{figure}
 \item To test the fact that the solution will always converge to a sine function moving in the $y$ direction when $\dfrac{n_0}{\omega_{ci}} = e^{Ax+B}$, i.e. $p = Ax+B$, we start with $u_0(x,y) = 10^{-5}sin(3 x)$ for $\tau = 0.1$, $n=32$, with $A=12$. 
 The solution remains unchanged up till $t=26$, and after the transitional time where the solution shifts from a sine function in the x direction to a sine function in the y-direction, the same behavior is observed (Figure \ref{fig:sin3x-32}).
\begin{figure}[H]
\centering
\includegraphics[scale=0.17]{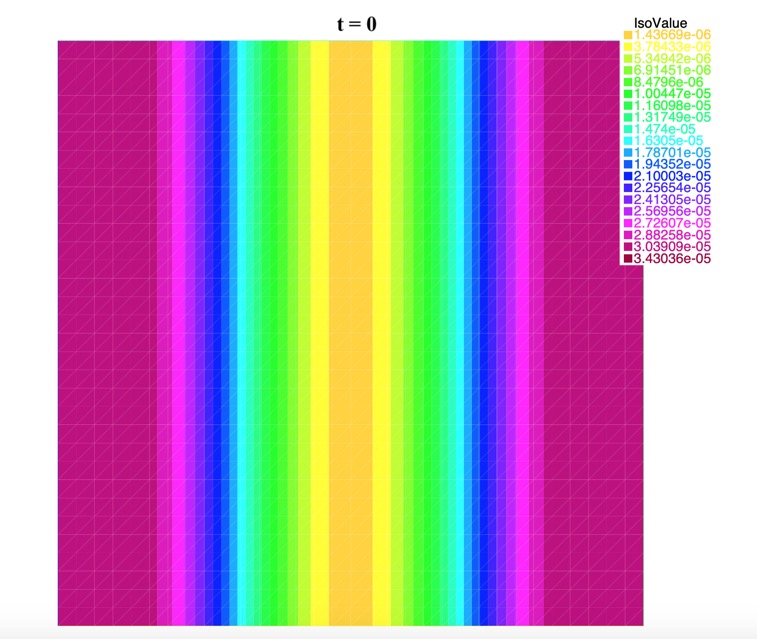}
\includegraphics[scale=0.17]{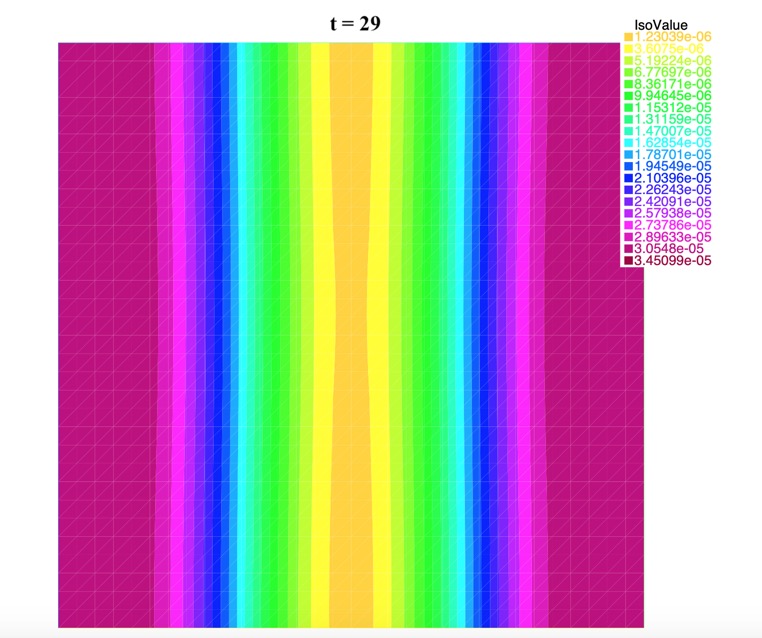}
\includegraphics[scale=0.17]{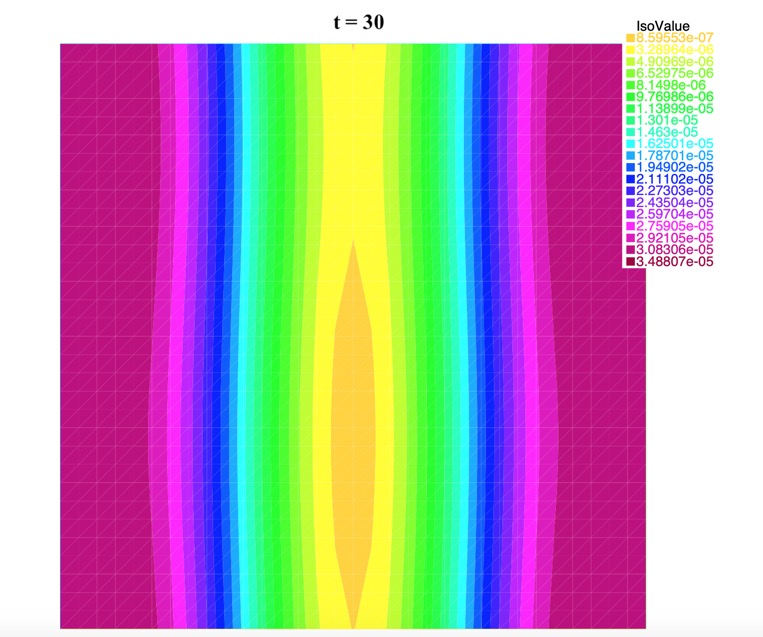}
\includegraphics[scale=0.17]{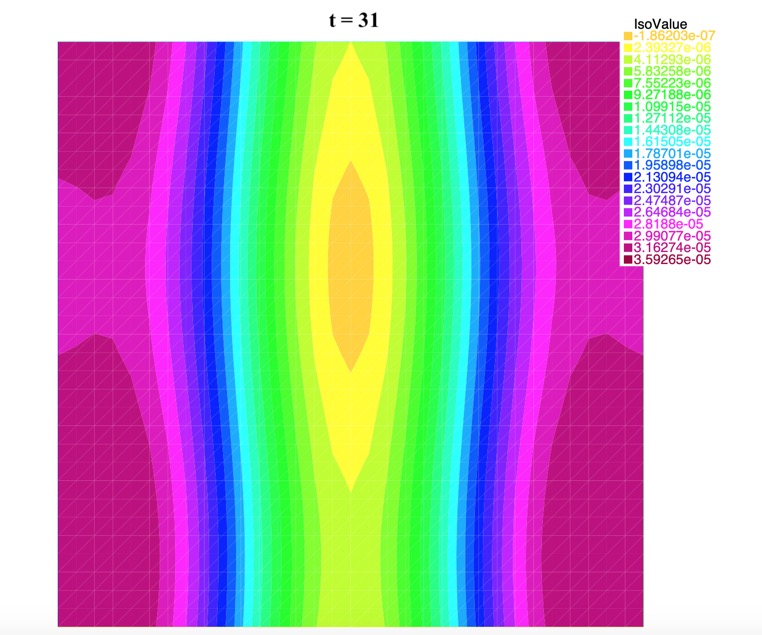}
\includegraphics[scale=0.17]{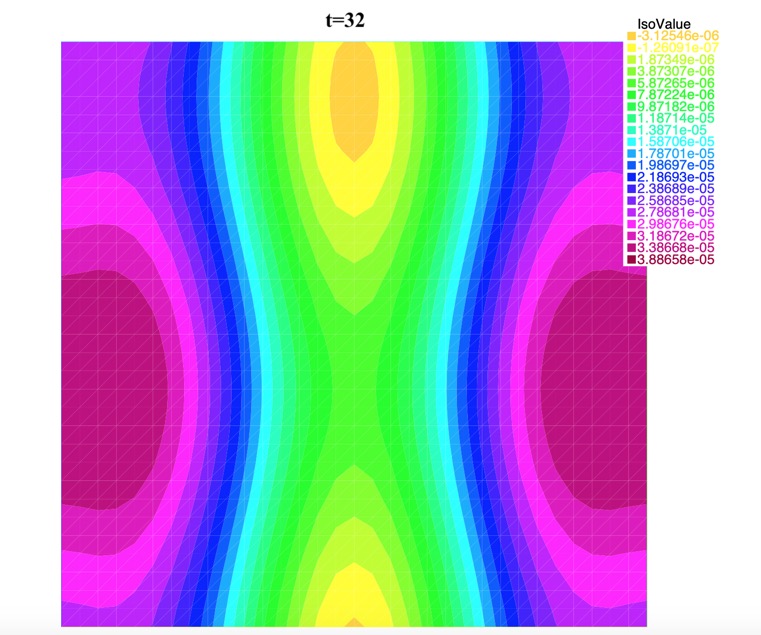}
\includegraphics[scale=0.17]{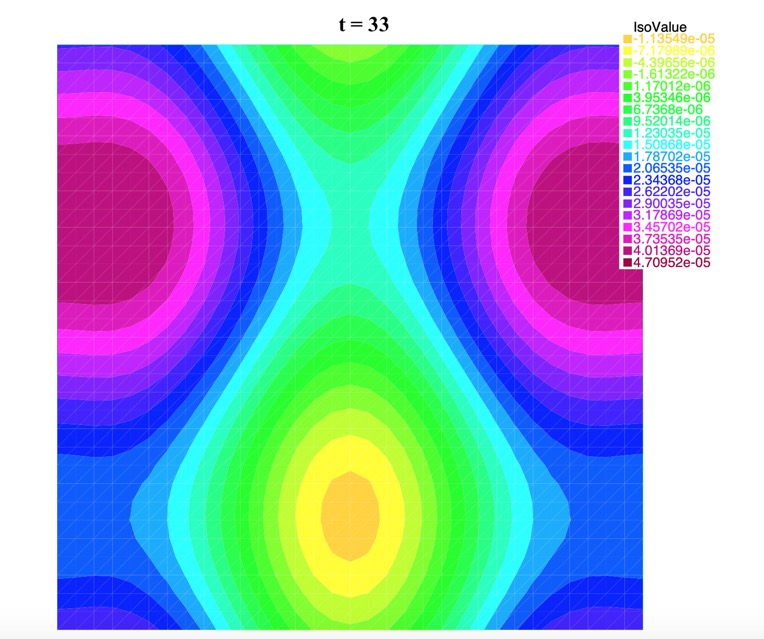}
\includegraphics[scale=0.17]{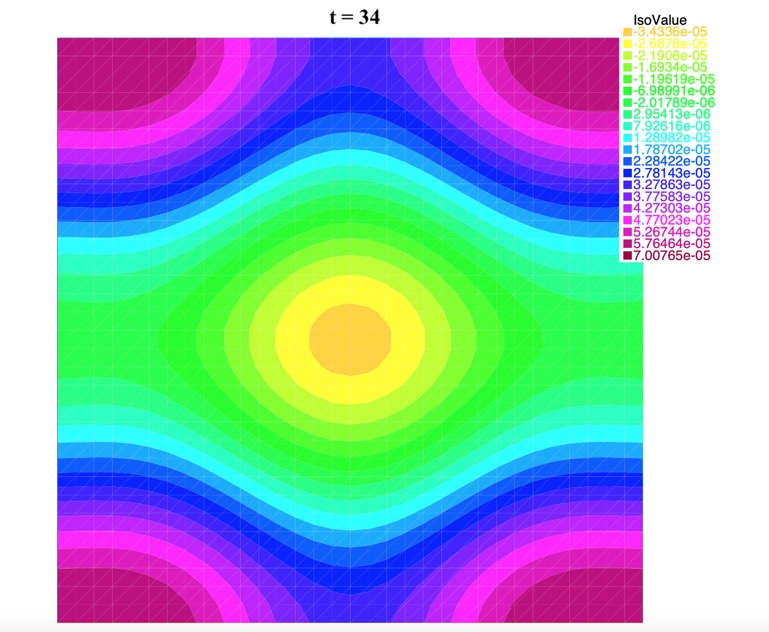}
\includegraphics[scale=0.17]{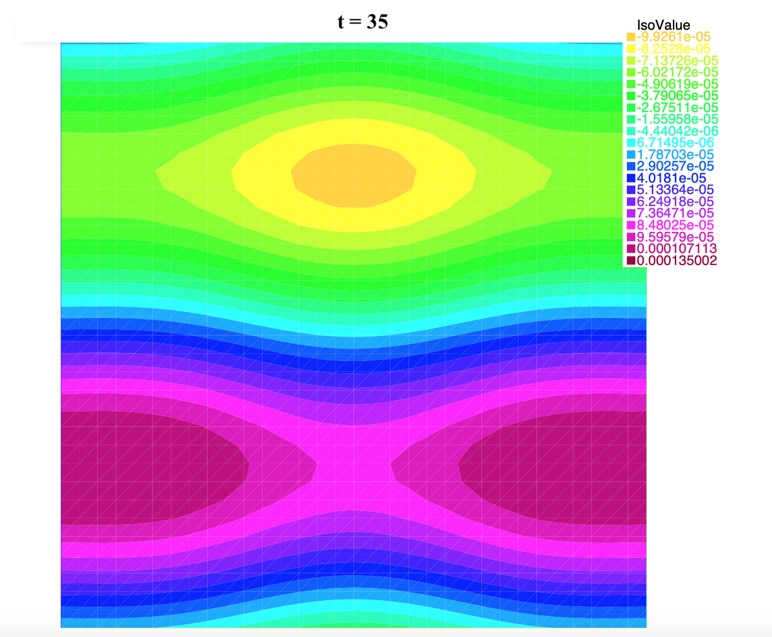}
\includegraphics[scale=0.17]{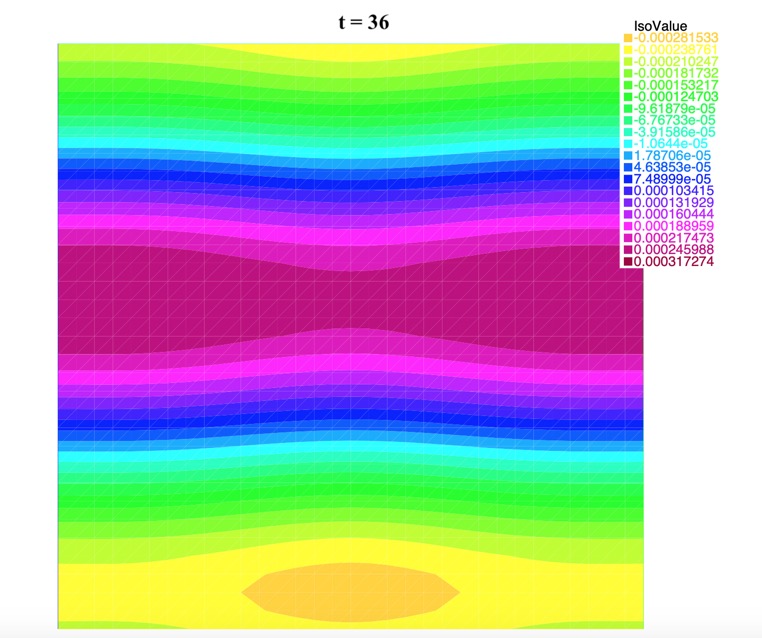}
\includegraphics[scale=0.17]{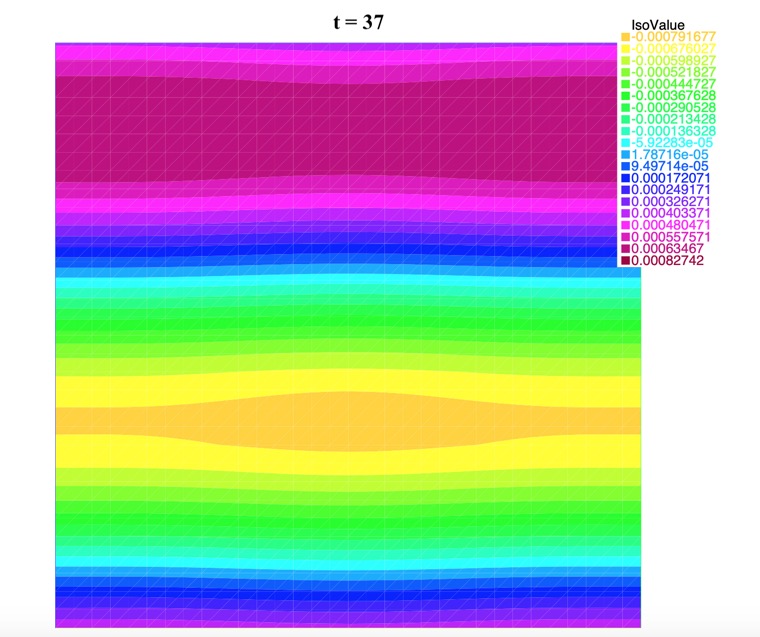} 
\includegraphics[scale=0.175]{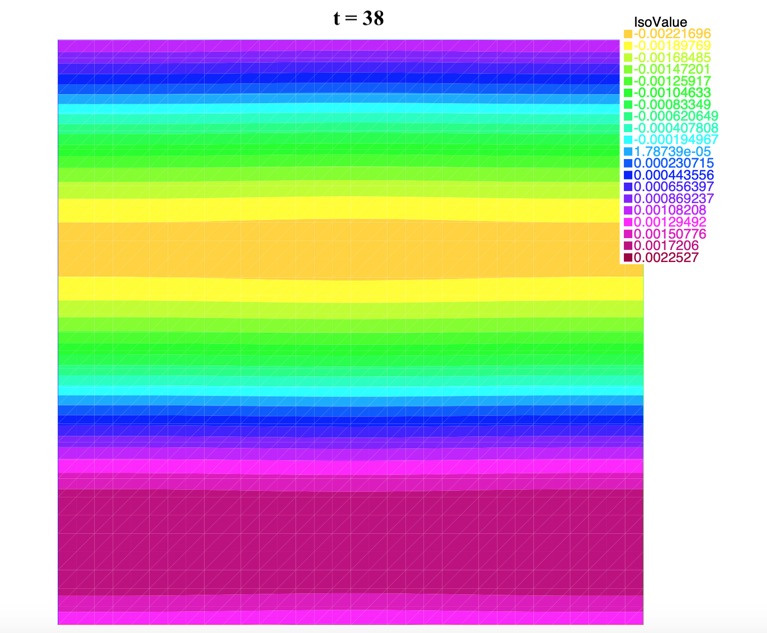} 
\includegraphics[scale=0.17]{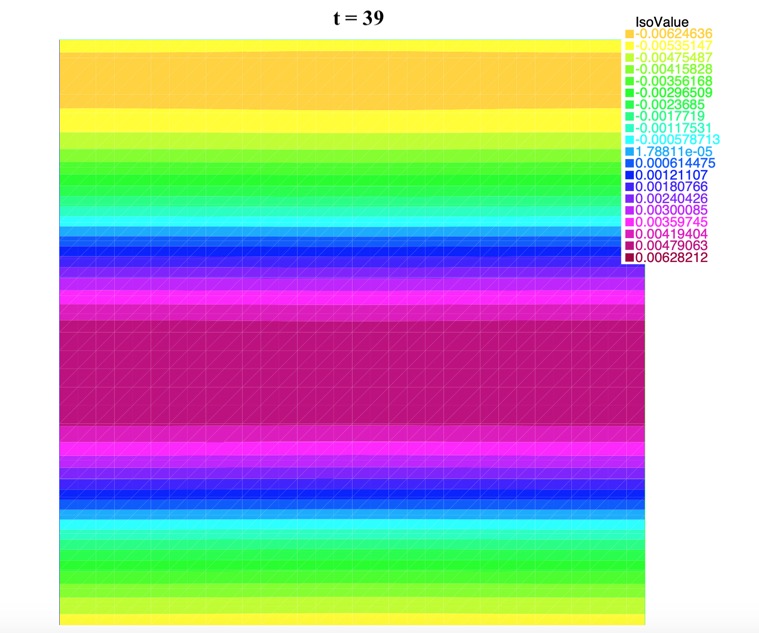} 
\includegraphics[scale=0.17]{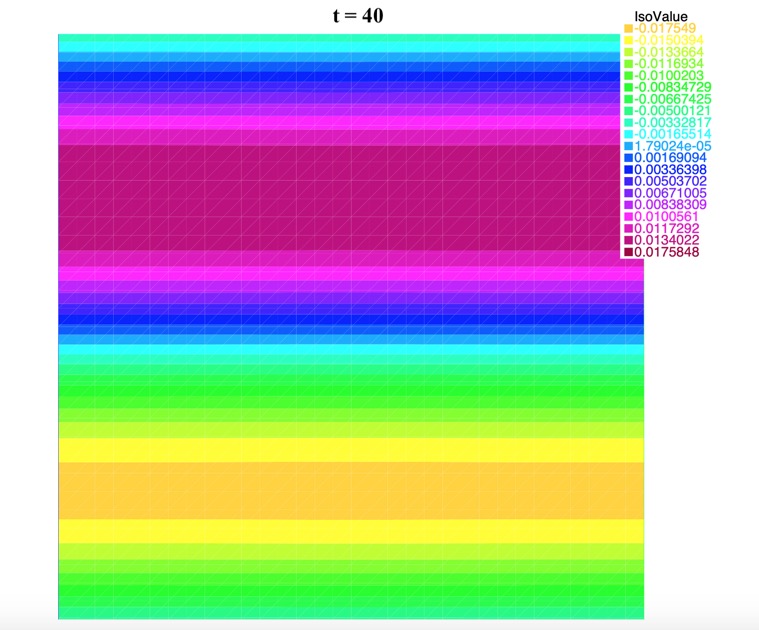} 
\includegraphics[scale=0.17]{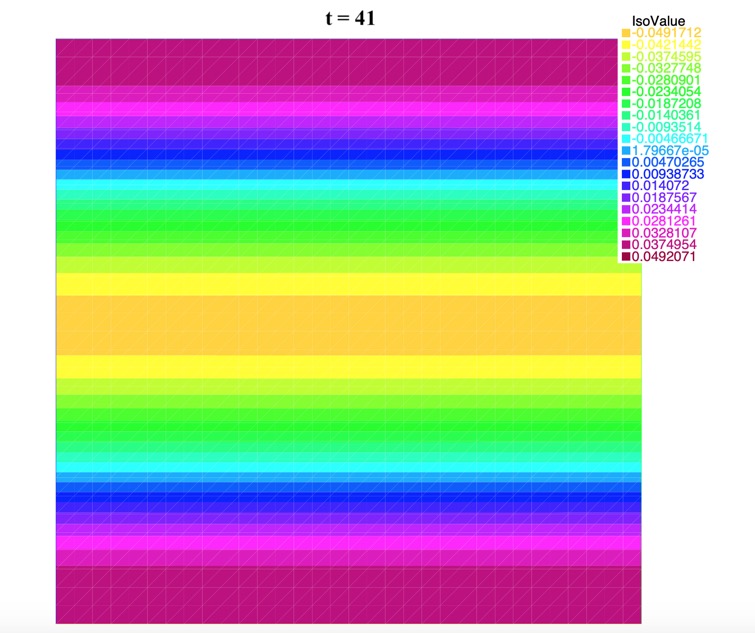} 
\includegraphics[scale=0.17]{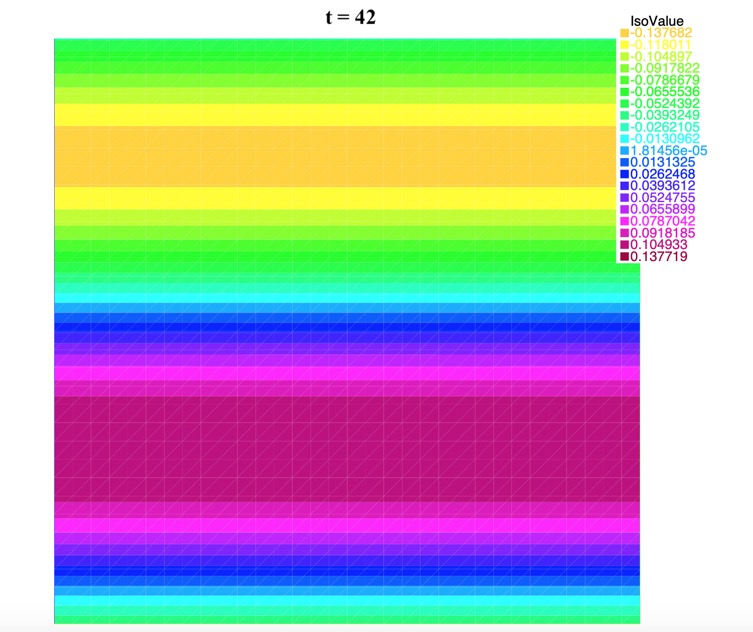} 
\caption{\it \small Time evolution of solution $u$ for $u_0= 10^{-5}sin(3x)$, $\tau = 0.1$, and a $33 \times 33$ grid on  $\Omega = [0,\pi] \times [0,\pi]$.}\label{fig:sin3x-32}
\end{figure}\newpage
\item It should be noted that even if we start with a random initial condition $u_0(x,y)$, after some time the solution converges to the same sine function moving in the $y$ direction. Figure \ref{fig:rand-32} shows the time evolution of the solution for $u_0(x,y) = 10^{-10}xy(x-2)sin(x)$, $\tau = 0.1$,  $\Omega = [0,\pi] \times [0,\pi]$, $n=32$ ($h = \pi/32 \approx 0.098175$) with $A=12$.
\begin{figure}[H]
\centering
\includegraphics[scale=0.17]{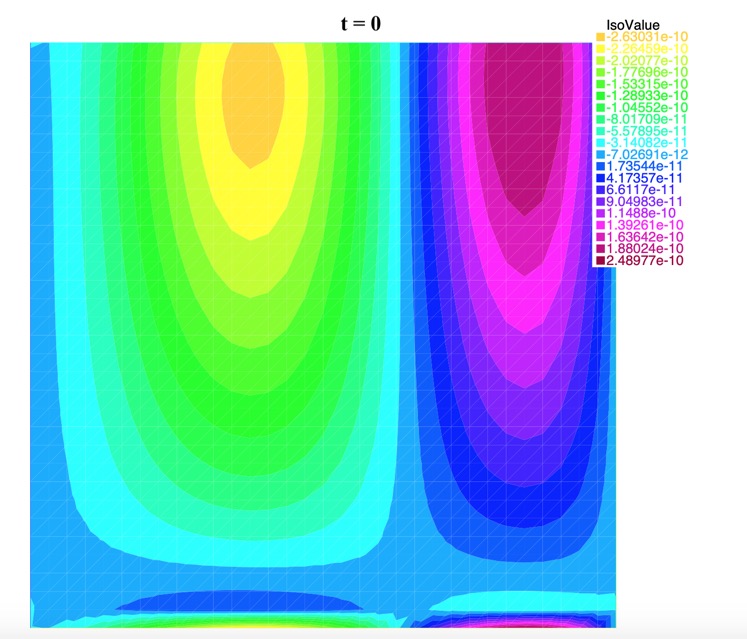}
\includegraphics[scale=0.17]{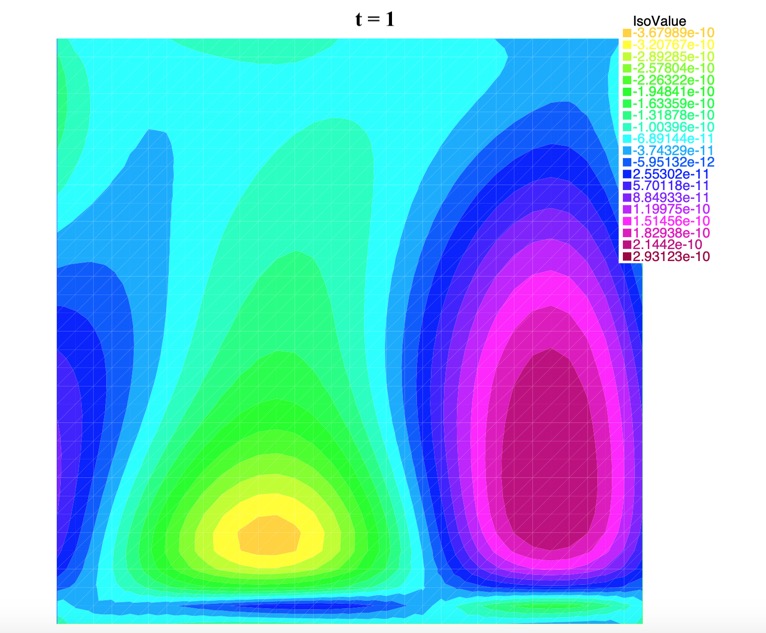}
\includegraphics[scale=0.17]{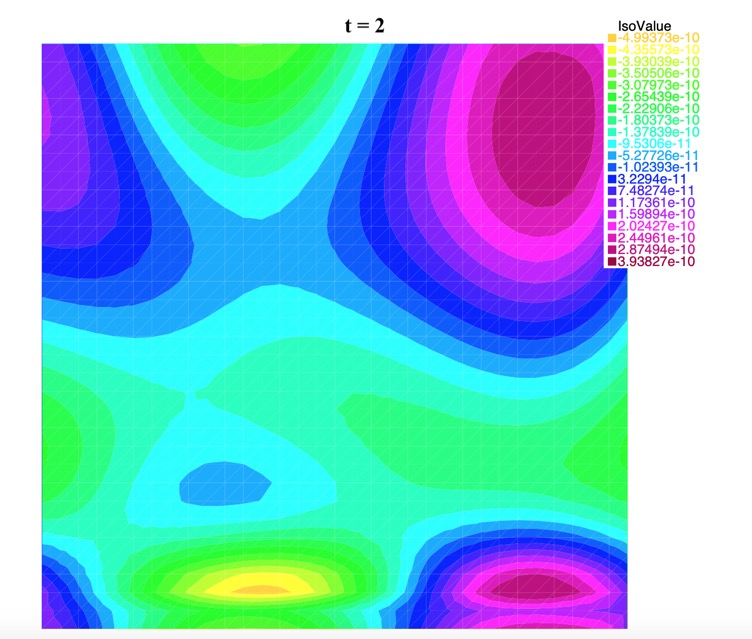}
\includegraphics[scale=0.17]{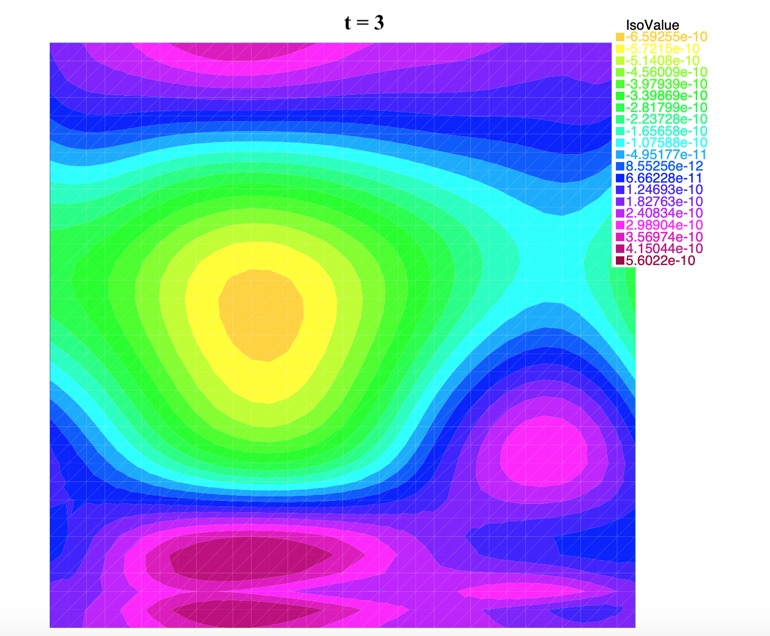}
\includegraphics[scale=0.17]{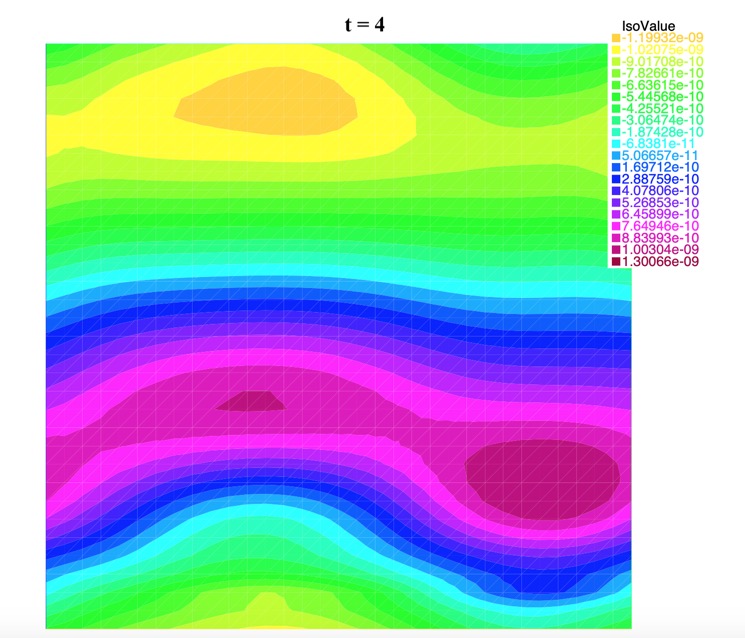}
\includegraphics[scale=0.17]{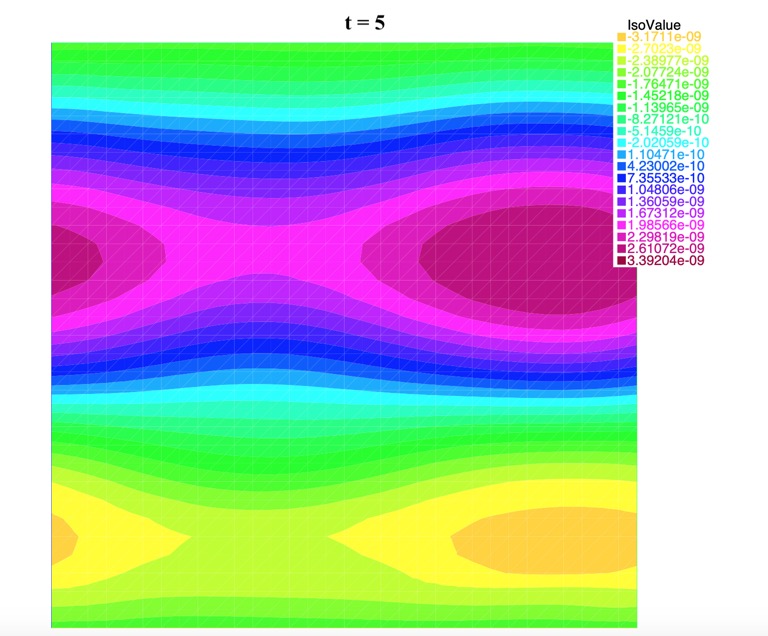}
\includegraphics[scale=0.17]{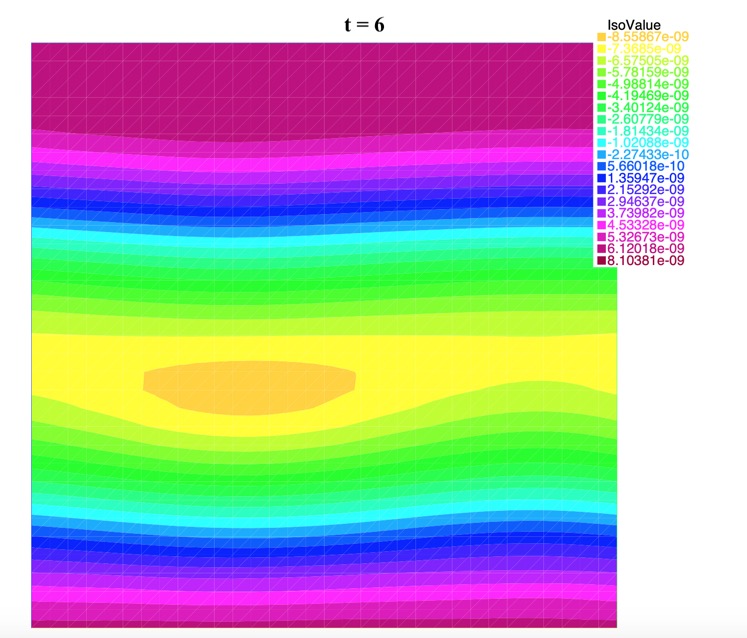}
\includegraphics[scale=0.17]{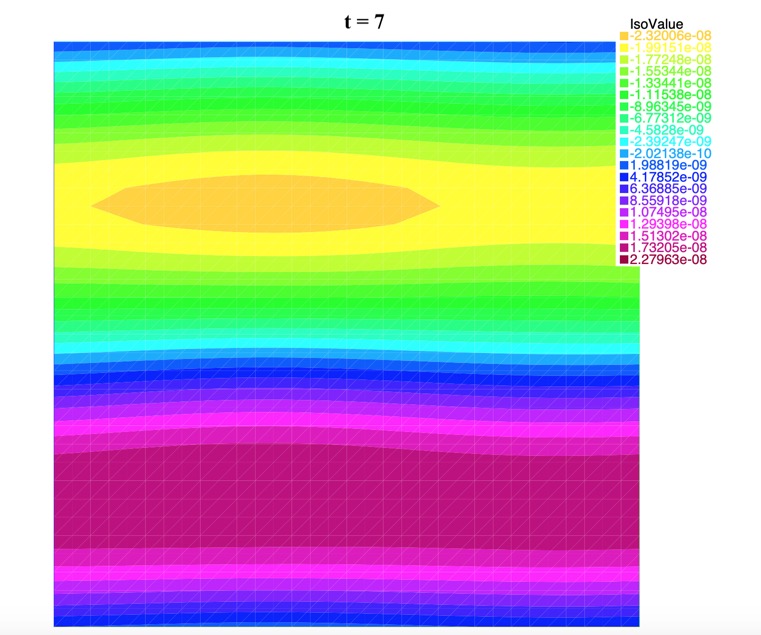}
\includegraphics[scale=0.17]{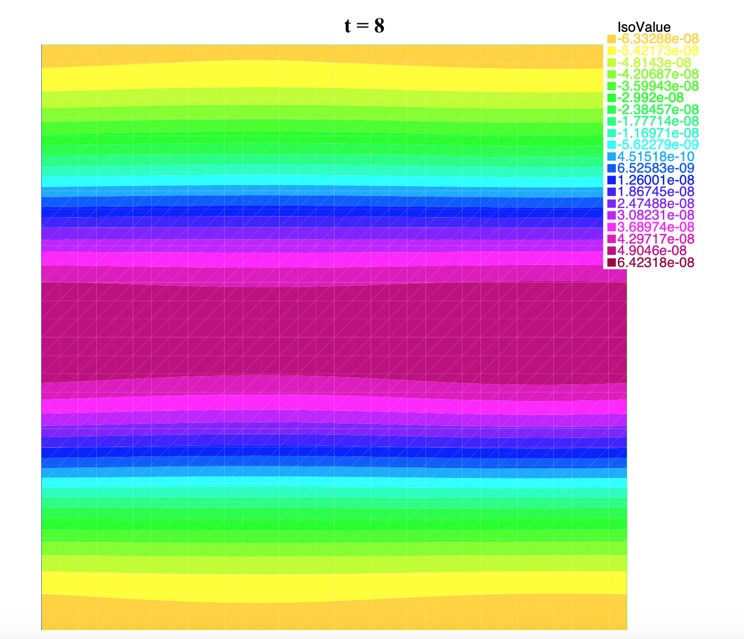}
\includegraphics[scale=0.17]{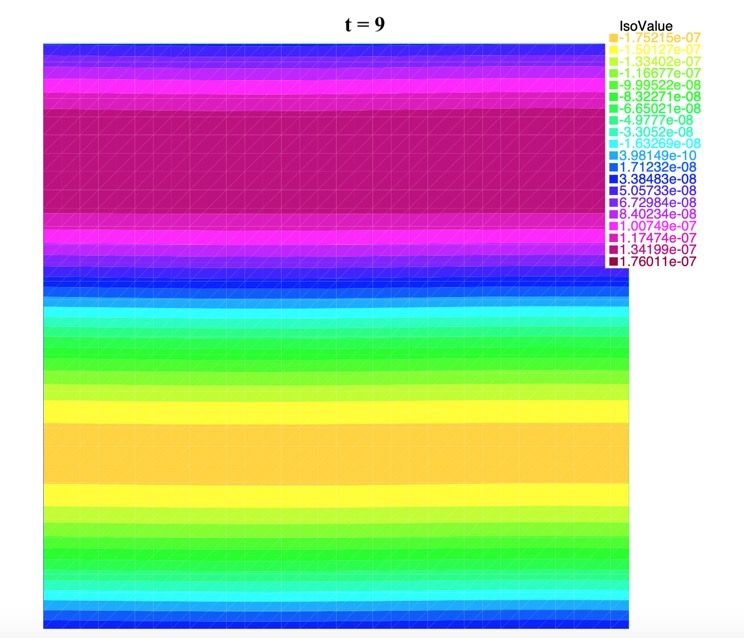} 
\includegraphics[scale=0.17]{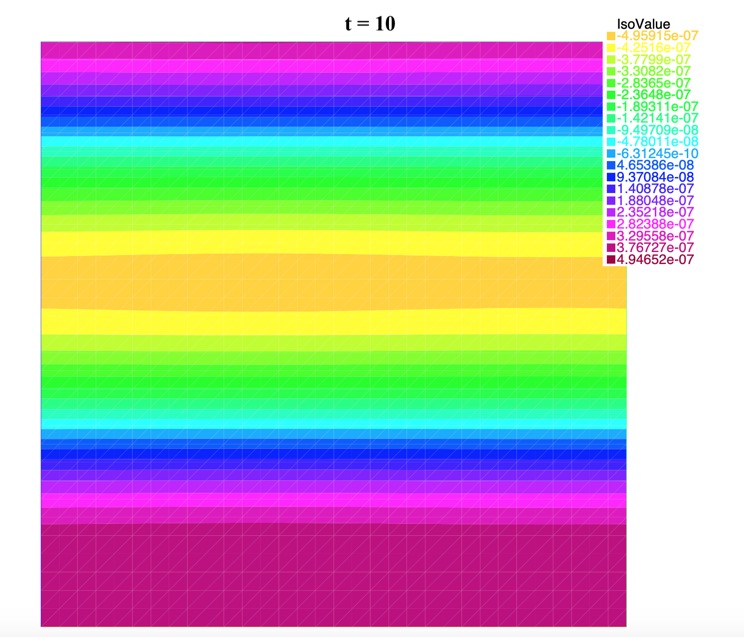} 
\includegraphics[scale=0.17]{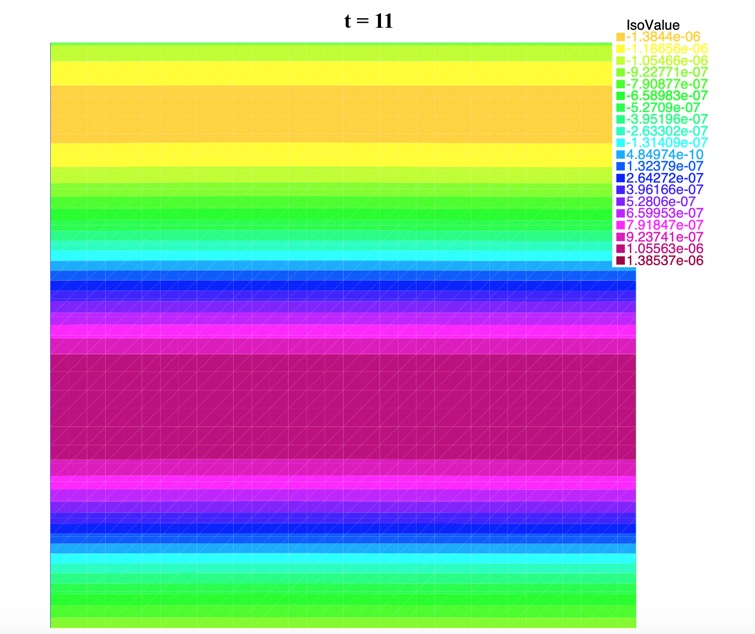} 
\includegraphics[scale=0.17]{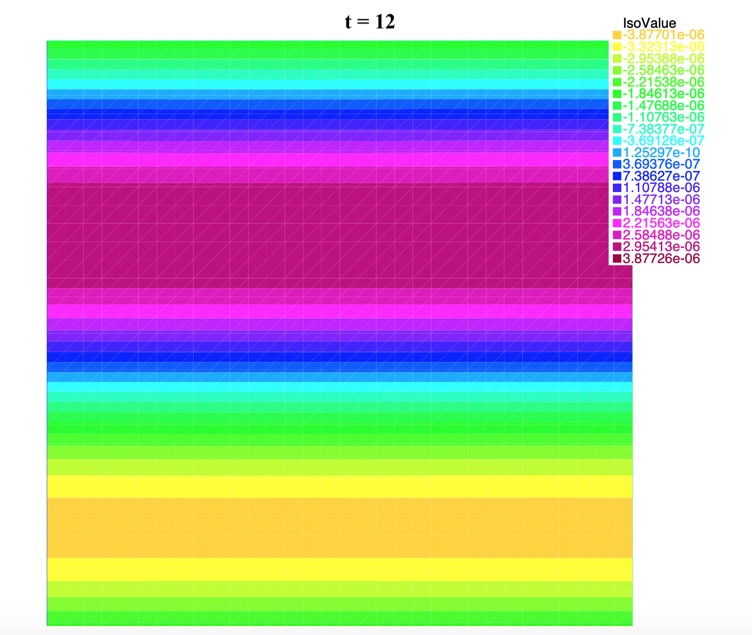} 
\includegraphics[scale=0.17]{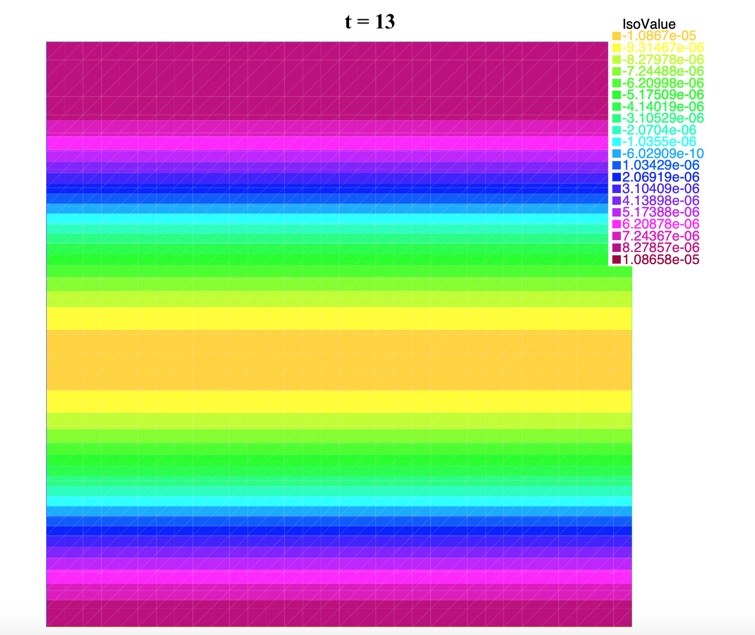} 
\includegraphics[scale=0.17]{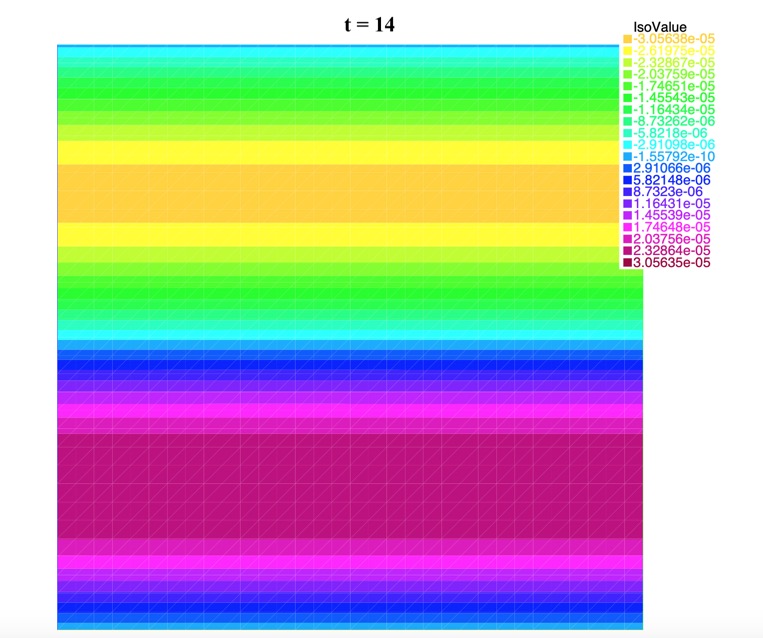} 
\caption{\it \small Time evolution of solution $u$ for $u_0= 10^{-10}xy(x-2)sin(x)$, $\tau = 0.1$, and a $33 \times 33$ grid on  $\Omega = [0,\pi] \times [0,\pi]$.}\label{fig:rand-32}
\end{figure}
 \end{enumerate}
Even though the theoretical study was done for the case where $p_y = 0$, however the algorithm works for any input function $p$. Note that if we set $p_x = 0$ and $p_y = 12$, then the solution will be moving in the x-direction in a similar manner as shown in the previous examples. We test the algorithm for the case where $n_0 = 10^{20} e^{-(x-10)^2/64-(y-10)^2/64}$, $w_{ci} = 10^7$, and $\Omega = [0,20]\times [0,20]$. Since $\nabla p = [-(x-10)/32, -(y-10)/32]$, the solution is expected to have a circular motion around the center of the domain $(10,10)$ as shown in Figure \ref{fig:gauss-64}. \vspace{-3mm}
\noindent  Note that if $\nabla p = [(x-10)/32, (y-10)/32]$ for the same initial conditions, then the solution will be moving in the opposite direction. 
\section{Concluding Remarks}\label{sec:Conclude}
To sum up the results of this paper, we have proven an existence theorem for the Hasegawa-Mima wave equation in $C(0,T;H^2\cap H^1_P)\times L^\infty(0,T;L^2)$. Uniqueness of the solution would require more regularity on the initial conditions as proven in \cite{kn}.
 On the other hand, we have also considered a full discretization scheme based on coupling the Finite-Element in space and a non-linear discretization of the time variable. The implementation of that scheme uses a semi-linear approach that provides a robust algorithm as revealed by early experiments.\\\\
Overall, future avenues of research include the following:
\begin{enumerate}
\item Proof of convergence of the solution to the nonlinear \eqref{HMC-Comp}, \eqref{HMC-Disc-Comp} schemes and the semi-linear \eqref{semi-lin-var} one as $\tau$ and $h$ go to zero.
\item Testing other alternatives to solve \eqref{HMC-Disc-Comp}. These include:
\begin{enumerate}
\item Using the fixed point approach discussed in Section \ref{Sub-Full-Discrete}, where the system can be written as follows:\\
Since $KU^j=MW^j$, and if we let $Z = MW(t)$, $Y = W(t+\tau)$, and $Y^1 = W(t)$,  then
$$Y^{k+1} = G(Y^k) \quad \iff \quad  (M-\tau RK^{-1}M)Y^{k+1} =  Z -{\tau}\,S(U^k)Y^k$$ 
which leads to the  
 following predictor-corrector scheme to move from time $t$ to $t+\tau$
\begin{equation}\label{HMC-Disc-Comp-it}
\left\{\begin{array}{ll}
U^1=U(t),\,Y^1=W(t),\,Z=MW(t)\,\mbox{error}=1,\,k=1& \\
\mbox{While error $>$ tol do:}&\\
\quad (M - \tau\, R\;K^{-1}\,M)\;Y^{k+1}  = Z  - \tau\,S(U^{k})Y^k &\\
\quad KU^{k+1}=MY^{k+1} &\\
 \quad \mbox{ error}=\max{\frac{||W^{k+1}-W^k||}{||W^{k+1}||},\frac{||U^{k+1}-U^k||}{||U^{k+1}||}}&\\
\quad k=k+1&\\
\mbox{end do}\\
U(t+\tau) = U^k \mbox{ \; and \;} W(t+\tau) = Y^k&
\end{array}\right.
\end{equation}
Following the results obtained in Section \ref{Sub-Full-Discrete} , specifically corollary \eqref{corol}, this algorithm is convergent. A first look at this approach indicates the necessity to deal with dense linear systems,  which matrix is \\$(M - \tau\, R\;K^{-1}\,M)$.
\noindent However, this difficulty can be lifted since $(M - \tau\, R\;K^{-1}\,M)Y = r$ is equivalent to $(K-\tau R)K^{-1}M Y = r$ where the solution is $Y = M^{-1}K(K-\tau R)^{-1}r$ can be obtained by solving two time-independent sparse systems
$(K-\tau R) \tilde{Y} = r \mbox{ followed by }MY = K\tilde{Y}$.
\item A second approach to handle \eqref{HMC-Disc-Comp} would be based on Newton's method.
\end{enumerate}
\item Another interesting problemhas to deal with the {Modon Traveling Waves Solutions} to \eqref{HMC-V1}. These solutions are obtained by  considering the pair of variables $(\xi,\eta)$ given by $\xi=x\,\,\eta=y-ct$, one  looks for solutions to (\ref{HMC2}) in the form $u(x,y,t)=\phi(\xi,\eta)=\phi(x,y-ct)\mbox{ and }w(x,y,t)= \psi(\xi,\eta)=\psi(x,y-ct)$. By defining $\forall t\in(0,T):\,\Omega_t=\{\xi,\eta\,| 0<\xi<L,\,-ct<\eta<L-ct\}, $ then in terms of $\phi$ and $\psi$, the system (\ref{HMC2}) reduces to be solved on $\Omega_0=\Omega$. Thus, with $\nabla=\nabla_{\xi,\eta}$, one seeks $\{\phi,\psi\}:\overline{\Omega}\to\mathbb{R}^2$, such that: 
\begin{equation}\label{HMCTravel}
\left\{\begin{array}{lll}
-c\psi_\eta + \vec{V}(\phi) \cdot \nabla\psi = k\phi_\eta & \mbox{on }  \Omega  &(1)\\
-\Delta \phi+\phi=\psi &\mbox{on }  \Omega &(2)\\
\mbox{PBC's on } \phi,\, \phi_\xi,\, \phi_\eta,\,\psi & \mbox{on }  \partial \Omega & (3)
\end{array}\right.
\end{equation}
Undergoing research is also being carried out on this problem. 
\end{enumerate}
\bibliographystyle{ieeetr}
\bibliography{references}

\newpage

\appendix
\section*{Appendices}
\addcontentsline{toc}{section}{Appendices}
\renewcommand{\thesubsection}{\Alph{subsection}}
\subsection{Assembling the Global matrix $S(U)$ }\label{sec:SU}
First, we derive the sparsity pattern of the global matrix $S(U)$ without any assumptions related to boundary conditions, i.e. define the general matrix $S(U)$, in Section  \ref{sec:GSU}. Then, we  derive the sparsity pattern of the matrix $S(U)$   assuming periodic boundary conditions in Section \ref{sec:PSU}.

\subsubsection{The General Matrix $S(U)$} \label{sec:GSU}
Given the meshing of Figure \ref{fig:mesh2} and without any assumptions on the boundary nodes, $S(U)$ is an $n^2\times n^2$ sparse matrix with at most 6 nonzero entries per row, as discussed below. To get those  nonzero entries, note that each node has at most 6 edges connecting it with its neighboring nodes, i.e. belongs to at most 6 triangles. We consider the 4 types of nodes shown in different colors in Figure \ref{fig:mesh2}: the $(n-2)^2$ black internal nodes that belong to 6 triangles, the $4(n-2)$ red boundary nodes that belong to 4 triangles, the 2 green corner nodes that belong to 2 triangles, and the 2 blue corner nodes that belong to 1 triangle.\vspace{5mm}

\noindent \textbf{1- Black Internal Nodes:}

Each of the $(n-2)^2$ black internal nodes with index $v = kn+i$ for $k = 1,2,...,n-2$ and $i=2,3,...,n-1$  belongs to triangles $T_{2j+1}, T_{2j+2}, T_{2j+3}, T_{2(j+n)}, T_{2(j+n)+1}, T_{2(j+n)+2}$,
where  $j = (n-1)(k-1)+(i-2)$. 

Thus, we first define the nonzero entries in row $v$ per local triangle and then add them up. \begin{itemize}
\item Triangle $T_{2j+1} = T_{2(j+1)-1}$ has global vertices 
$\begin{array}{lcl}
(j+1)+(k-1) &=& nk-n+i-1\\
 (j+1)+(n+1)+(k-1) &=& nk+i = v \\
 (j+1)+n+(k-1) &=& nk+i-1
 \end{array}$\\
 
 \noindent where node $v$ is the second local vertex of triangle $T_{2j+1}$. 
 Thus row $v$ of the matrix $S(U)$ has the entry \\$\dfrac{1}{6} U_{nk-n+i-1} - \dfrac{1}{6} U_{nk+i-1}$ in columns  $nk-n+i-1, nk+i-1$, and $v = nk+i$. 
 \item Triangle $T_{2j+2} = T_{2(j+1)}$ has global vertices 
$\begin{array}{lcl}
(j+1)+(k-1) &=& nk-n+i-1\\
 (j+1)+1+(k-1) &=& nk-n+i\\
 (j+1)+(n+1)+(k-1) &=& nk+i = v 
 \end{array}$\\
 
 \noindent where node $v$ is the third local vertex of triangle $T_{2j+2}$. 
 Thus row $v$ of the matrix $S(U)$ has the entry\\ $\dfrac{1}{6} U_{nk-n+i} - \dfrac{1}{6} U_{nk-n+i-1}$ in columns  $nk-n+i-1, nk-n+i$, and $v = nk+i$. 
  \item Triangle $T_{2j+3} = T_{2(j+2)-1}$ has global vertices 
$\begin{array}{lcl}
(j+2)+(k-1) &=& nk-n+i\\
 (j+2)+(n+1)+(k-1) &=& nk+i+1 \\
 (j+2)+n+(k-1) &=& nk+i = v 
 \end{array}$\\
 
 \noindent where node $v$ is the third local vertex of triangle $T_{2j+3}$. 
 Thus row $v$ of the matrix $S(U)$ has the entry\\ $\dfrac{1}{6}U_{nk+i+1} - \dfrac{1}{6}U_{nk-n+i}$ in columns  $nk-n+i, v = nk+i$, and $ nk+i+1$. 
  \item Triangle $T_{2(j+n)}$ has global vertices 
$\begin{array}{lcl}
(j+n)+(k) &=& nk+i-1\\
 (j+n)+1+(k) &=& nk+i = v\\
 (j+n)+(n+1)+(k) &=& nk+i+n  
 \end{array}$\\
 
\noindent where node $v$ is the second local vertex of triangle $T_{2(j+n)}$. 
 Thus row $v$ of the matrix $S(U)$ has the entry\\ $\dfrac{1}{6} U_{nk+i-1} - \dfrac{1}{6}U_{nk+i+n}$ in columns  $nk+i-1, v = nk+i$, and $nk+i+n $. 
 \item Triangle $T_{2(j+n)+1} = T_{2(j+n+1)-1}$ has global vertices 
$\begin{array}{lcl}
(j+n+1)+(k) &=& nk+i = v\\
 (j+n+1)+(n+1)+(k) &=& nk+i +n+1\\
 (j+n+1)+n+(k) &=& nk+i+n
 \end{array}$\\
 
 \noindent where node $v$ is the first local vertex of triangle $T_{2(j+n)+1}$. 
 Thus row $v$ of the matrix $S(U)$ has the entry\\ $\dfrac{1}{6} U_{nk+i+n}- \dfrac{1}{6} U_{ nk+i +n+1}$ in columns  $v = nk+i, nk+i+n$, and $nk+i+n+1$. 
  \item Triangle $T_{2(j+n)+2} = T_{2(j+n+1)} $ has global vertices 
$\begin{array}{lcl}
(j+n+1)+(k) &=& nk+i = v\\
 (j+n+1)+1+(k) &=& nk+i+1 \\
 (j+n+1)+(n+1)+(k) &=& nk+i+n +1 
 \end{array}$\\
 
 \noindent where node $v$ is the first local vertex of triangle $T_{2(j+n)+2}$. 
 Thus row $v$ of the matrix $S(U)$ has the entry\\ $\dfrac{1}{6} U_{ nk+i+n +1 } - \dfrac{1}{6}U_{nk+i+1}$ in columns  $v = nk+i, nk+i+1 $, and $nk+i+n+1 $. 
\end{itemize}
Thus the nonzero entries in row $v = nk+i$ of $S(U)$ for $k = 1,2,...,n-2$ and $i=2,3,...,n-1$ are in columns  
$$\begin{array}{lll}
nk-n+i-1 : & \dfrac{1}{6} U_{nk-n+i-1} - \dfrac{1}{6} U_{nk+i-1} + \dfrac{1}{6} U_{nk-n+i} - \dfrac{1}{6} U_{nk-n+i-1}  &=  \dfrac{1}{6} U_{nk-n+i} - \dfrac{1}{6} U_{nk+i-1} \nonumber \vspace{3mm}\\
nk-n+i : & \dfrac{1}{6} U_{nk-n+i} - \dfrac{1}{6} U_{nk-n+i-1}+ \dfrac{1}{6}U_{nk+i+1} - \dfrac{1}{6}U_{nk-n+i} &= \dfrac{1}{6}U_{nk+i+1}- \dfrac{1}{6} U_{nk-n+i-1}  \nonumber \vspace{3mm}\\
 nk+i-1 : & \dfrac{1}{6} U_{nk-n+i-1} - \dfrac{1}{6} U_{nk+i-1} + \dfrac{1}{6} U_{nk+i-1} - \dfrac{1}{6}U_{nk+i+n} &=\dfrac{1}{6} U_{nk-n+i-1} - \dfrac{1}{6}U_{nk+i+n}\nonumber \vspace{3mm}\\
 v = nk+i : & \dfrac{1}{6} U_{nk-n+i-1} - \dfrac{1}{6} U_{nk+i-1} + \dfrac{1}{6} U_{nk-n+i}- \dfrac{1}{6} U_{nk-n+i-1} &+ \dfrac{1}{6}U_{nk+i+1} - \dfrac{1}{6}U_{nk-n+i}\nonumber \vspace{3mm}\\
 &+\dfrac{1}{6} U_{nk+i-1}  - \dfrac{1}{6}U_{nk+i+n}+ \dfrac{1}{6} U_{nk+i+n} - \dfrac{1}{6} U_{ nk+i +n+1}&+\dfrac{1}{6} U_{ nk+i+n +1 } - \dfrac{1}{6}U_{nk+i+1}=0\nonumber \vspace{3mm}\\
  nk+i+1 : &  \dfrac{1}{6}U_{nk+i+1} - \dfrac{1}{6}U_{nk-n+i}+\dfrac{1}{6} U_{ nk+i+n +1 } - \dfrac{1}{6}U_{nk+i+1}& = \dfrac{1}{6} U_{ nk+i+n +1 }- \dfrac{1}{6}U_{nk-n+i}\nonumber \vspace{3mm}\\
  nk+i+n : & \dfrac{1}{6} U_{nk+i-1} - \dfrac{1}{6}U_{nk+i+n}+\dfrac{1}{6} U_{nk+i+n} - \dfrac{1}{6} U_{ nk+i +n+1}& = \dfrac{1}{6} U_{nk+i-1}- \dfrac{1}{6} U_{ nk+i +n+1} \nonumber \vspace{3mm}\\
   nk+i+n+1 : &\dfrac{1}{6} U_{nk+i+n} - \dfrac{1}{6} U_{ nk+i +n+1}+\dfrac{1}{6} U_{ nk+i+n +1 } - \dfrac{1}{6}U_{nk+i+1} & = \dfrac{1}{6} U_{nk+i+n}- \dfrac{1}{6}U_{nk+i+1}\nonumber
\end{array}$$

\noindent \textbf{2- Red Boundary Nodes:}

The red boundary nodes are of 4 types:
\begin{enumerate}
\item[\bf a)] the left boundary with index $v = kn+1$ and $k = 1,2,...,n-2$, that belong to triangles\\ $T_{2(k-1)(n-1)+1}, T_{2k(n-1)+1}$, and $T_{2k(n-1)+2}$.\\
Thus, we first define the nonzero entries in row $v=kn+1$ per local triangle and then add them up. 
\begin{itemize}
\item Triangle $T_{2(k-1)(n-1)+1} = T_{2{(}nk-n-k+2{)}-1}$ has global vertices 

$\begin{array}{lcl}
nk-n-k+2+(k-1) &=& nk-n+1 \\
 nk-n-k+2+(n+1)+(k-1) &=& nk+2\\
nk-n-k+2+n+(k-1) &=& nk+1 = v
 \end{array}$\\
 
 \noindent where node $v$ is the third local vertex of triangle $T_{2(k-1)(n-1)+1}$. 
 Thus row $v$ of the matrix $S(U)$ has the entry $\dfrac{1}{6} U_{ nk+2} - \dfrac{1}{6} U_{ nk-n+1}$ in columns  $nk-n+1, v = nk+1$, and $nk+2$. 
 \item Triangle $T_{2k(n-1)+1} = T_{2\big{(}k(n-1)+1\big{)}-1}$ has global vertices 

$\begin{array}{lcl}
k(n-1)+1+(k) &=& nk+1 = v \\
 k(n-1)+1+(n+1)+(k) &=& nk+n+2\\
k(n-1)+1+n+(k) &=& nk+n+1 
 \end{array}$\\
 
 \noindent where node $v$ is the first local vertex of triangle $T_{2k(n-1)+1}$. 
 Thus row $v$ of the matrix $S(U)$ has the entry $\dfrac{1}{6} U_{ nk+n+1} - \dfrac{1}{6} U_{ nk+n+2}$ in columns  $v = nk+1, nk+n+1, $, and $nk+n+2$. 
 \item Triangle $T_{2k(n-1)+2} = T_{2(k(n-1)+1)} $ has global vertices 

$\begin{array}{lcl}
k(n-1)+1+(k) &=& nk+1 = v\\
 k(n-1)+1+1+(k) &=& nk+2 \\
 k(n-1)+1+(n+1)+(k) &=& nk+n +2 
 \end{array}$\\
 
 \noindent where node $v$ is the first local vertex of triangle $T_{2k(n-1)+2}$. 
 Thus row $v$ of the matrix $S(U)$ has the entry $\dfrac{1}{6} U_{ nk+n +2 }- \dfrac{1}{6}U_{nk+2}$ in columns  $v = nk+1, nk+2 $, and $nk+n+2 $. 
\end{itemize}
Thus the nonzero entries in row $v = nk+1$ of $S(U)$ for $k = 1,2,...,n-2$ are in columns  
$$\begin{array}{lll}
 nk-n+1 : & \dfrac{1}{6} U_{ nk+2} - \dfrac{1}{6} U_{ nk-n+1} &= \dfrac{1}{6} U_{ nk+2} - \dfrac{1}{6} U_{ nk-n+1}  \nonumber \vspace{3mm}\\
 v = nk+1 : & \dfrac{1}{6} U_{ nk+2} - \dfrac{1}{6} U_{ nk-n+1}+\dfrac{1}{6} U_{ nk+n+1} - \dfrac{1}{6} U_{ nk+n+2}&  \nonumber \vspace{3mm}\\
 &+\dfrac{1}{6} U_{ nk+n +2 } - \dfrac{1}{6}U_{nk+2}&= \dfrac{1}{6} U_{ nk+n+1} - \dfrac{1}{6} U_{ nk-n+1}\nonumber \vspace{3mm}\\
  nk+2 : &\dfrac{1}{6} U_{ nk+2} - \dfrac{1}{6} U_{ nk-n+1} +\dfrac{1}{6} U_{ nk+n +2 } - \dfrac{1}{6}U_{nk+2}& = \dfrac{1}{6} U_{ nk+n +2 } - \dfrac{1}{6} U_{ nk-n+1}  \nonumber \vspace{3mm}\\
  nk+n+1 : & \dfrac{1}{6} U_{ nk+n+1} - \dfrac{1}{6} U_{ nk+n+2}& = \dfrac{1}{6} U_{ nk+n+1} - \dfrac{1}{6} U_{ nk+n+2}\nonumber \vspace{3mm}\\
   nk+n+2 : &\dfrac{1}{6} U_{ nk+n+1} - \dfrac{1}{6} U_{ nk+n+2}+\dfrac{1}{6} U_{ nk+n +2 } - \dfrac{1}{6}U_{nk+2}& = \dfrac{1}{6} U_{ nk+n+1} - \dfrac{1}{6}U_{nk+2}\nonumber \vspace{-1mm}
\end{array}$$
\item[\bf b)] the right boundary with index $v = kn$ and $k = 2,3,...,n-1$, that belong to triangles\\ $T_{2(k-1)(n-1)}, T_{2(k-1)(n-1)-1}$, and $T_{2k(n-1)}$.\\
Thus, we first define the nonzero entries in row $v=kn$ per local triangle and then add them up. 
\begin{itemize}
\item Triangle $T_{2(k-1)(n-1)-1}$ has global vertices 
$\begin{array}{lcl}
(k-1)(n-1)+(k-2) &=& nk-n-1 \\
 (k-1)(n-1)+(n+1)+(k-2) &=& nk = v\\
(k-1)(n-1)+n+(k-2) &=& nk-1 
 \end{array}$\\
 
 \noindent where node $v$ is the second local vertex of triangle $T_{2(k-1)(n-1)-1}$. 
 Thus row $v$ of the matrix $S(U)$ has the entry $\dfrac{1}{6} U_{ nk-n-1} - \dfrac{1}{6} U_{ nk-1}$ in columns  $nk-n-1, nk-1$, and $ v = nk$. 
  \item Triangle $T_{2(k-1)(n-1)} $ has global vertices 
$\begin{array}{lcl}
(k-1)(n-1)+(k-2) &=& nk-n-1 \\
 (k-1)(n-1)+1+(k-2) &=& nk-n  \\
 (k-1)(n-1)+(n+1)+(k-2) &=& nk = v
 \end{array}$\\
 
 \noindent where node $v$ is the third local vertex of triangle $T_{2(k-1)(n-1)}$. 
 Thus row $v$ of the matrix $S(U)$ has the entry $\dfrac{1}{6} U_{ nk-n } - \dfrac{1}{6}U_{nk-n-1}$ in columns  $nk-n-1, nk-n $, and $v = nk $. 
  \item Triangle $T_{2k(n-1)} $ has global vertices 
$\begin{array}{lcl}
k(n-1)+(k-1) &=& nk-1 \\
 k(n-1)+1+(k-1) &=& nk = v \\
 k(n-1)+(n+1)+(k-1) &=& nk+n 
 \end{array}$\\
 
 \noindent where node $v$ is the second local vertex of triangle $T_{2k(n-1)}$. 
 Thus row $v$ of the matrix $S(U)$ has the entry $\dfrac{1}{6} U_{ nk-1 } - \dfrac{1}{6}U_{nk+n}$ in columns  $ nk-1, v = nk$, and $nk+n $. 
\end{itemize}
Thus the nonzero entries in row $v = nk$ of $S(U)$ for $k = 2,3,...,n-1$ are in columns  
$$\begin{array}{lll}
 nk-n-1 : &\dfrac{1}{6} U_{ nk-n-1} - \dfrac{1}{6} U_{ nk-1}+\dfrac{1}{6} U_{ nk-n} - \dfrac{1}{6}U_{nk-n-1} &=  \dfrac{1}{6} U_{ nk-n} - \dfrac{1}{6} U_{ nk-1} \nonumber \vspace{3mm}\\
 nk-n : &\dfrac{1}{6} U_{ nk-n }- \dfrac{1}{6}U_{nk-n-1} &=  \dfrac{1}{6} U_{ nk-n } - \dfrac{1}{6}U_{nk-n-1}  \nonumber \vspace{3mm}\\
 nk-1 : &\dfrac{1}{6} U_{ nk-n-1} - \dfrac{1}{6} U_{ nk-1}+\dfrac{1}{6} U_{ nk-1 } - \dfrac{1}{6}U_{nk+n} &= \dfrac{1}{6} U_{ nk-n-1}  - \dfrac{1}{6}U_{nk+n} \nonumber \vspace{3mm}\\
 v = nk : &\dfrac{1}{6} U_{ nk-n-1} - \dfrac{1}{6} U_{ nk-1}+ \dfrac{1}{6} U_{ nk-n } - \dfrac{1}{6}U_{nk-n-1}&  \nonumber \vspace{3mm}\\
 &+\dfrac{1}{6} U_{ nk-1 } - \dfrac{1}{6}U_{nk+n} &= \dfrac{1}{6} U_{ nk-n} - \dfrac{1}{6} U_{ nk+n} \nonumber \vspace{3mm}\\
  nk+n : &\dfrac{1}{6} U_{ nk-1 } - \dfrac{1}{6}U_{nk+n} & = \dfrac{1}{6} U_{ nk-1 } - \dfrac{1}{6}U_{nk+n}\nonumber
\end{array}$$
\item[\bf c)] the lower boundary with index $v = i$ and $i = 2,3,...,n-1$, that belong to triangles\\ $T_{2(i-1)}, T_{2(i-1)+1}$, and $T_{2(i-1)+2}$.\\
Thus, we first define the nonzero entries in row $v=i$ per local triangle and then add them up. 
\begin{itemize}
 \item Triangle $T_{2(i-1)} $ has global vertices 
$\begin{array}{lcl}
i-1 &=& i-1 \\
 (i-1)+1&=& i  = v \\
 (i-1)+(n+1) &=& i+n
 \end{array}$\\
 
 \noindent where node $v = i$ is the second local vertex of triangle $T_{2(i-1)}$. 
 Thus row $v = i$ of the matrix $S(U)$ has the entry $\dfrac{1}{6} U_{ i-1 } - \dfrac{1}{6}U_{i+n}$ in columns  $i-1, v=i $, and $i+ n $. 
\item Triangle $T_{2(i-1)+1} = T_{2i-1}$ has global vertices 
$\begin{array}{lcl}
i &=&  v \\
 i+(n+1) &=& i+n+1\\
i+n &=& i+n
 \end{array}$\\ 
 
 \noindent where node $v = i$ is the first local vertex of triangle $T_{2(i-1)+1}$. 
 Thus row $v = i$ of the matrix $S(U)$ has the entry $\dfrac{1}{6} U_{ i+n} - \dfrac{1}{6} U_{ i+n+1}$ in columns  $ v = i, i+n$, and $ i+n+1$. 
 \item Triangle $T_{2(i-1) +2} = T_{2i}$ has global vertices 
$\begin{array}{lcl}
i &=&  v \\
 i+1&=& i +1  \\
 i+(n+1) &=& i+n+1
 \end{array}$\\ 
 
 \noindent where node $v = i$ is the first local vertex of triangle $T_{2(i-1)}$. 
 Thus row $v = i$ of the matrix $S(U)$ has the entry $\dfrac{1}{6}U_{i+n+1} - \dfrac{1}{6}U_{i+1}$ in columns  $v=i, i+1 $, and $i+ n+1$. 
\end{itemize}
Thus the nonzero entries in row $v = i$ of $S(U)$ for $i = 2,3,...,n-1$ are in columns  
$$\begin{array}{lll}
 i-1 : &\dfrac{1}{6} U_{ i-1 } - \dfrac{1}{6}U_{i+n} &= \dfrac{1}{6} U_{ i-1 } - \dfrac{1}{6}U_{i+n}   \nonumber \vspace{3mm}\\
 v = i : &\dfrac{1}{6} U_{ i-1 } - \dfrac{1}{6}U_{i+n} +\dfrac{1}{6} U_{ i+n} - \dfrac{1}{6} U_{ i+n+1} + \dfrac{1}{6}U_{i+n+1} - \dfrac{1}{6}U_{i+1}&= \dfrac{1}{6} U_{ i-1 }- \dfrac{1}{6}U_{i+1}\nonumber \vspace{3mm}\\
  i+1 : &\dfrac{1}{6}U_{i+n+1} - \dfrac{1}{6}U_{i+1} & = \dfrac{1}{6}U_{i+n+1} - \dfrac{1}{6}U_{i+1} \nonumber \vspace{3mm}\\
     i+n : &\dfrac{1}{6} U_{ i-1 } - \dfrac{1}{6}U_{i+n}+\dfrac{1}{6} U_{ i+n} - \dfrac{1}{6} U_{ i+n+1} &= \dfrac{1}{6} U_{ i-1 }- \dfrac{1}{6} U_{ i+n+1} \nonumber \vspace{3mm}\\
   i+n+1 : &\dfrac{1}{6} U_{ i+n} - \dfrac{1}{6} U_{ i+n+1}+\dfrac{1}{6}U_{i+n+1} - \dfrac{1}{6}U_{i+1} &= \dfrac{1}{6} U_{ i+n} - \dfrac{1}{6}U_{i+1} \nonumber \vspace{3mm}\\
\end{array}$$
\item[\bf d)] the upper boundary with index $v = n(n-1)+i$ and $i = 2,3,...,n-1$, that belong to triangles\\ $T_{2(n-1)(n-2)+2(i-1)-1}, T_{2(n-1)(n-2)+2(i-1)}$, and $T_{2(n-1)(n-2)+2(i-1)+1}$.
\\
Thus, we first define the nonzero entries in row $v=n(n-1)+i$ per local triangle and then add them up. 
\begin{itemize}\item Triangle $T_{2(n-1)(n-2)+2(i-1)-1}$ has global vertices 

$\begin{array}{lcl}
(n-1)(n-2)+i-1+(n-2)  &=&  n(n-2)+i-1 \\
 (n-1)(n-2)+i-1+(n+1)+(n-2)  &=& n(n-1)+i = v\\
(n-1)(n-2)+i-1+n +(n-2) &=& n(n-1)+i-1 
 \end{array}$\\ 
 
 \noindent where node $v =  n(n-1)+i$ is the second local vertex of triangle $T_{2(n-1)(n-2)+2(i-1)-1}$. 
 Thus row $v = n(n-1)+i$ of the matrix $S(U)$ has the entry $\dfrac{1}{6} U_{ n(n-2)+i-1} - \dfrac{1}{6} U_{ n(n-1)+i-1}$ in columns  $ n(n-2)+i-1,  n(n-1)+i-1$, and $v = n(n-1)+i$. 
 \item Triangle $T_{2(n-1)(n-2)+2(i-1)} $ has global vertices 
 
$\begin{array}{lcl}
(n-1)(n-2)+ i-1+(n-2) &=& n(n-2)+i-1 \\
 (n-1)(n-2)+ (i-1)+1+(n-2)&=&  n(n-2)+i   \\
 (n-1)(n-2)+(i-1)+(n+1)+(n-2) &=& n(n-1)+i = v
 \end{array}$\\
 
 \noindent where node $v = n^2-n+i$ is the third local vertex of triangle $T_{2(i-1)}$. 
 Thus row $v = i$ of the matrix $S(U)$ has the entry $\dfrac{1}{6} U_{ n(n-2)+i } - \dfrac{1}{6}U_{n(n-2)+i-1 }$ in columns  $n(n-2)+i-1 ,n(n-2)+i $, and $v = n(n-1)+i$. 
\item Triangle $T_{2(n-1)(n-2)+2(i-1)+1} = T_{2(n-1)(n-2)+2i-1}$ has global vertices 

$\begin{array}{lcl}
(n-1)(n-2)+i+(n-2)  &=&  n(n-2)+i \\
 (n-1)(n-2)+i+(n+1)+(n-2)  &=& n(n-1)+i+1\\
(n-1)(n-2)+i+n +(n-2) &=& n(n-1)+i = v
 \end{array}$\\ 
 
 \noindent where node $v =  n(n-1)+i$ is the third local vertex of triangle $T_{2(n-1)(n-2)+2(i-1)+1}$. 
 Thus row $v = n(n-1)+i$ of the matrix $S(U)$ has the entry $\dfrac{1}{6} U_{n(n-1)+i+1}- \dfrac{1}{6} U_{ n(n-2)+i}$ in columns  $ n(n-2)+i, v = n(n-1)+i$, and $n(n-1)+i+1$. 
\end{itemize}
Thus the nonzero entries in row $v = n(n-1)+i$ of $S(U)$ for $i = 2,3,...,n-1$, are in columns  
$$\begin{array}{lll}
 n(n-2)+i-1 : &\dfrac{1}{6} U_{ n(n-2)+i-1} - \dfrac{1}{6} U_{ n(n-1)+i-1}+\dfrac{1}{6} U_{ n(n-2)+i } &   \nonumber \vspace{3mm}\\
 &- \dfrac{1}{6}U_{n(n-2)+i-1 } &=\dfrac{1}{6} U_{ n(n-2)+i } - \dfrac{1}{6} U_{ n(n-1)+i-1}   \nonumber \vspace{3mm}\\
 \end{array}$$
 $$\begin{array}{lll}
 n(n-2)+i : &\dfrac{1}{6} U_{ n(n-2)+i } - \dfrac{1}{6}U_{n(n-2)+i-1 } +\dfrac{1}{6} U_{n(n-1)+i+1}&   \nonumber \vspace{3mm}\\
 &- \dfrac{1}{6} U_{ n(n-2)+i}&=  \dfrac{1}{6} U_{n(n-1)+i+1} - \dfrac{1}{6}U_{n(n-2)+i-1 } \nonumber \vspace{3mm}\\
 n(n-1)+i-1 : & \dfrac{1}{6} U_{ n(n-2)+i-1} - \dfrac{1}{6} U_{ n(n-1)+i-1}&= \dfrac{1}{6} U_{ n(n-2)+i-1} - \dfrac{1}{6} U_{ n(n-1)+i-1} \nonumber \vspace{3mm}\\
 v =  n(n-1)+i : & \dfrac{1}{6} U_{ n(n-2)+i-1} - \dfrac{1}{6} U_{ n(n-1)+i-1}+\dfrac{1}{6} U_{ n(n-2)+i } & \nonumber \vspace{3mm}\\
 &- \dfrac{1}{6}U_{n(n-2)+i-1 }+\dfrac{1}{6} U_{n(n-1)+i+1} - \dfrac{1}{6} U_{ n(n-2)+i}&= \dfrac{1}{6} U_{n(n-1)+i+1}-\dfrac{1}{6} U_{ n(n-1)+i-1} \nonumber \vspace{3mm}\\
 n(n-1)+ i+1 : &\dfrac{1}{6} U_{n(n-1)+i+1} - \dfrac{1}{6} U_{ n(n-2)+i} & = \dfrac{1}{6} U_{n(n-1)+i+1} - \dfrac{1}{6} U_{ n(n-2)+i} \nonumber \vspace{3mm}\\
\end{array}$$
\end{enumerate} \newpage

\noindent \textbf{3- Green Corner Nodes:}

There are 2 green corners with index :
\begin{enumerate}
\item[\bf a)] $v = 1$ that belongs to triangles $T_{1}$, and $T_{2}$. First, we define the nonzero entries in row $v=1$ per local triangle.
\begin{itemize}\item Triangle $T_{1}$ has global vertices $v=1, n+2, n+1$
 where node $v =  1$ is the first local vertex of triangle $T_{1}$. 
 Thus row $v = 1$ of the matrix $S(U)$ has the entry $\dfrac{1}{6} U_{ n+1} - \dfrac{1}{6} U_{ n+2}$ in columns  $v = 1, n+1 $, and $n+2$. 
 \item Triangle $T_{2}$ has global vertices $v=1, 2, n+2$
 where node $v =  1$ is the first local vertex of triangle $T_{1}$. 
 Thus row $v = 1$ of the matrix $S(U)$ has the entry $\dfrac{1}{6} U_{ n+2} - \dfrac{1}{6} U_{ 2}$ in columns  $v = 1,  2$, and $n+2$. 
\end{itemize}
Thus the nonzero entries in row $v = 1$ of $S(U)$  are in columns  
$$\begin{array}{lll}
 v=1 : & \dfrac{1}{6} U_{ n+1} - \dfrac{1}{6} U_{ n+2}+\dfrac{1}{6} U_{ n+2} - \dfrac{1}{6} U_{ 2}&=  \dfrac{1}{6} U_{ n+1} -  \dfrac{1}{6} U_{ 2} \nonumber \vspace{3mm}\\
 2 : & \dfrac{1}{6} U_{ n+2} - \dfrac{1}{6} U_{ 2}& = \dfrac{1}{6} U_{ n+2} - \dfrac{1}{6} U_{ 2}  \nonumber \vspace{3mm}\\
 n+1 : &\dfrac{1}{6} U_{ n+1} - \dfrac{1}{6} U_{ n+2} &= \dfrac{1}{6} U_{ n+1} - \dfrac{1}{6} U_{ n+2}  \nonumber \vspace{3mm}\\
  n+2 : & \dfrac{1}{6} U_{ n+1} - \dfrac{1}{6} U_{ n+2}+\dfrac{1}{6} U_{ n+2} - \dfrac{1}{6} U_{ 2}& = \dfrac{1}{6} U_{ n+1} - \dfrac{1}{6} U_{ 2} \nonumber \vspace{3mm}\\
\end{array}$$
\item[\bf b)] $v = n^2$ that belongs to triangles $T_{2(n-1)(n-1)- 1}$, and $T_{2(n-1)(n-1)}$.\\ First, we define the nonzero entries in row $v=n^2$ per local triangle and then add them up. 
\begin{itemize}\item Triangle $T_{2(n-1)(n-1)-1}$ has global vertices $n(n-1)-1,v=n^2, n^2 - 1$
 where node $v =  n^2$ is the second local vertex of triangle $T_{2(n-1)(n-1)-1}$. 
 Thus row $v = n^2$ of the matrix $S(U)$ has the entry $\dfrac{1}{6} U_{ n(n-1)-1} - \dfrac{1}{6} U_{ n^2-1}$ in columns  $n(n-1)-1, n^2-1$, and $v = n^2$. 
 \item Triangle $T_{2(n-1)(n-1)-1}$ has global vertices $n(n-1)-1, n(n-1), v=n^2$,
 where node $v =  n^2$ is the third local vertex of triangle $T_{1}$. 
 Thus row $v = n^2$ of the matrix $S(U)$ has the entry $\dfrac{1}{6} U_{ n(n-1)} - \dfrac{1}{6} U_{ n(n-1)-1}$ in columns  $n(n-1)-1, n(n-1)$, and $v=n^2$. 
\end{itemize}
Thus the nonzero entries in row $v = n^2$ of $S(U)$  are in columns  
$$\begin{array}{lll}
 n(n-1)-1 : &\dfrac{1}{6} U_{ n(n-1)-1} - \dfrac{1}{6} U_{ n^2-1} +\dfrac{1}{6} U_{ n(n-1)} - \dfrac{1}{6} U_{ n(n-1)-1}&= \dfrac{1}{6} U_{ n(n-1)} - \dfrac{1}{6} U_{ n^2-1}\nonumber \vspace{3mm}\\
 n(n-1) : &\dfrac{1}{6} U_{ n(n-1)} - \dfrac{1}{6} U_{ n(n-1)-1} & = \dfrac{1}{6} U_{ n(n-1)} - \dfrac{1}{6} U_{ n(n-1)-1}  \nonumber \vspace{3mm}\\
 \end{array}$$
 $$\begin{array}{lll}
 n^2-1 : &\dfrac{1}{6} U_{ n(n-1)-1} - \dfrac{1}{6} U_{ n^2-1} &= \dfrac{1}{6} U_{ n(n-1)-1} - \dfrac{1}{6} U_{ n^2-1}  \nonumber \vspace{3mm}\\
  v=n^2 : &\dfrac{1}{6} U_{ n(n-1)-1} - \dfrac{1}{6} U_{ n^2-1}+\dfrac{1}{6} U_{ n(n-1)} - \dfrac{1}{6} U_{ n(n-1)-1} & = \dfrac{1}{6} U_{ n(n-1)} - \dfrac{1}{6} U_{ n^2-1}\nonumber
\end{array}$$
\end{enumerate}

\noindent \textbf{4- Blue Corner Nodes:}

There are 2 blue corners with index :
\begin{enumerate}
\item[\bf a)] $v = n$ that belongs to triangle $T_{2(n-1)}$ with global vertices $n-1,v = n,2n$,
 where node $v =  n$ is the second local vertex of triangle $T_{2(n-1)}$. 
 Thus row $v = n$ of the matrix $S(U)$ has the entry $\dfrac{1}{6} U_{ n-1} - \dfrac{1}{6} U_{ 2n}$ in columns  $n-1,v = n$, and $2n$. 
\item[\bf b)] $v = n(n-1)+1$ that belongs to triangle $T_{2(n-2)(n-1)+ 1}$ with global vertices $n(n-2)+1,n(n-1) + 2, \\v = n(n-1) + 1$
 where node $v =  n(n-1) + 1$ is the third local vertex of triangle $T_{2(n-1)(n-1)-1}$. 
 Thus row $v = n(n-1) + 1$ of the matrix $S(U)$ has the entry $\dfrac{1}{6} U_{ n(n-1) + 2} - \dfrac{1}{6} U_{n(n-2)+1}$ in columns  $n(n-2)+1,\\ v = n(n-1)+1$, and $n(n-1)+2$. 
\end{enumerate}
\noindent For $n=5$, the matrix $S(U)$ corresponding to the mesh in Figure \ref{fig:mesh2}, 
has the following sparsity pattern:
$$S(U) = \left[
\begin{array}{ccccc|ccccc|ccccc|ccccc|ccccc}
* & * & 0&0&0& * &*&0&0&0&0&0&0&0&0&0&0&0&0&0&0&0&0&0&0\\
*&*&*&0&0&0&*&*&0&0&0&0&0&0&0&0&0&0&0&0&0&0&0&0&0\\
0&*&*&*&0&0&0&*&*&0&0&0&0&0&0&0&0&0&0&0&0&0&0&0&0\\
0&0&*&*&*&0&0&0&*&*&0&0&0&0&0&0&0&0&0&0&0&0&0&0&0\\
0&0&0&*&*&0&0&0&0&*&0&0&0&0&0&0&0&0&0&0&0&0&0&0&0 \\ \hline
*&0&0&0&0&*&*&0&0&0&*&*&0&0&0&0&0&0&0&0&0&0&0&0&0\\
*&*&0&0&0&*&0&*&0&0&0&*&*&0&0&0&0&0&0&0&0&0&0&0&0\\
0&*&*&0&0&0&*&0&*&0&0&0&*&*&0&0&0&0&0&0&0&0&0&0&0\\
0&0&*&*&0&0&0&*&0&*&0&0&0&*&*&0&0&0&0&0&0&0&0&0&0\\
0&0&0&*&*&0&0&0&*&*&0&0&0&0&*&0&0&0&0&0&0&0&0&0&0\\  \hline
0&0&0&0&0&*&0&0&0&0&*&*&0&0&0&*&*&0&0&0&0&0&0&0&0\\
0&0&0&0&0&*&*&0&0&0&*&0&*&0&0&0&*&*&0&0&0&0&0&0&0\\
0&0&0&0&0&0&*&*&0&0&0&*&0&*&0&0&0&*&*&0&0&0&0&0&0\\
0&0&0&0&0&0&0&*&*&0&0&0&*&0&*&0&0&0&*&*&0&0&0&0&0\\
0&0&0&0&0&0&0&0&*&*&0&0&0&*&*&0&0&0&0&*&0&0&0&0&0\\  \hline
0&0&0&0&0&0&0&0&0&0&*&0&0&0&0&*&*&0&0&0&*&*&0&0&0\\
0&0&0&0&0&0&0&0&0&0&*&*&0&0&0&*&0&*&0&0&0&*&*&0&0\\
0&0&0&0&0&0&0&0&0&0&0&*&*&0&0&0&*&0&*&0&0&0&*&*&0\\
0&0&0&0&0&0&0&0&0&0&0&0&*&*&0&0&0&*&0&*&0&0&0&*&*\\
0&0&0&0&0&0&0&0&0&0&0&0&0&*&*&0&0&0&*&*&0&0&0&0&*\\ \hline 
0&0&0&0&0&0&0&0&0&0&0&0&0&0&0&*&0&0&0&0&*&*&0&0&0\\
0&0&0&0&0&0&0&0&0&0&0&0&0&0&0&*&*&0&0&0&*&*&*&0&0\\
0&0&0&0&0&0&0&0&0&0&0&0&0&0&0&0&*&*&0&0&0&*&*&*&0\\
0&0&0&0&0&0&0&0&0&0&0&0&0&0&0&0&0&*&*&0&0&0&*&*&*\\
0&0&0&0&0&0&0&0&0&0&0&0&0&0&0&0&0&0&*&*&0&0&0&*&*\\
\end{array} \right]
$$ 

\noindent In general, the $S(U)$ matrix is an $n^2 \times n^2$ block tridiagonal matrix with at most 6 nonzero entries per row and 6 per column.
$$S(U) =\dfrac{1}{6} \left[ \begin{array}{cccccc} A_{1,1} & A_{1,2}& 0&0&\cdots&0\\
A_{2,1} & A_{2,2}& A_{2,3}&0&\cdots&0\\
0&\ddots&\ddots&\ddots&\ddots&\vdots\\
\vdots&\ddots&A_{j,k}&A_{j,j}&A_{j,i}&0\\
0&\cdots&0&A_{i,j}&A_{i,i}&A_{i,n}\\
0&\cdots&0&0&A_{n,i}&A_{n,n}\\
\end{array}\right]
$$
where $i = n-1, j = n-2, k= n-3$, and the block matrices are of size $n\times n$ with the following sparsity patterns:
\begin{itemize}
\item $A_{1,1}$ and $A_{n,n}$ are tridiagonal matrices, each with $3n-2$ nonzero entries.
\item $A_{i,i}$ for $i = 2,3,.., n-1$ are tridiagonal matrices with zero diagonal entries except for the first and the last.
\item $A_{i+1,i}$ for $i = 1,2,3,.., n-1$ are lower bidiagonal matrices, with $2n-1$ nonzero entries.
\item $A_{i,i+1}$ for $i = 1,2,3,.., n-1$ are upper bidiagonal matrices, with $2n-1$ nonzero entries.
\end{itemize}
Thus, $S(U)$ has a total of $2(3n-2)+(n-2)(2n)+(2n-2)(2n-1) = 6n-4+2n^2-4n+4n^2-6n+2 = 6n^2 - 4n - 2$ nonzero entries. As for the nonzero entries in $S(U)$, each is of order $||U||_{\infty}$ since each nonzero entry is of the form $$\dfrac{1}{6}(U_i-U_j) \leq \dfrac{1}{6}|U_i-U_j| \leq \dfrac{1}{6}(|U_i|+|U_j|) \leq \dfrac{1}{3} ||U||_{\infty}$$
\subsubsection{The Matrix $S(U)$ assuming Periodic Boundary Conditions on $U$}\label{sec:PSU}
Given the meshing of Figure \ref{fig:mesh2} and assuming periodic boundary conditions, where the values of U are equal at the upper and lower red vertices, the left and right red boundary vertices, and the corner vertices i.e. 
\begin{eqnarray}
U_{kn+1}&=&U_{kn+n}, \qquad \qquad \qquad \qquad \quad \;\; for \;\; k = 1,2,..,n-2\\
U_{i}&=&U_{i+n(n-1)}, \qquad \qquad \qquad \qquad for \;\; i = 2, 3,..,n-1\\ 
U_1 &=& U_n \; = \; U_{1+n(n-1)} \; = \; U_{n^2}
\end{eqnarray} and assuming that our domain is torus shaped and that the vector U is of size $(n-1)^2$, $$U = [U_1,.., U_{n-1},U_{n+1},...,U_{2n-1}, U_{2n+1},...,U_{3n-1}, ......, U_{n(n-1)-1}]$$ 
then $S(U)$ is an $(n-1)^2\times (n-1)^2$ sparse matrix with at most 6 nonzero entries per row, as discussed below. This ``periodic" $S(U)$ matrix, can be obtained from the general one described in the previous section by merging/adding the rows corresponding to equal $U$ entries and also the columns, and using the periodicity of $U$.

 Note that since rows/columns $n, 2n, 3n , ..., (n-1)n$ of the general $S(U)$ matrix are merged with other rows/columns, then the indices have to be reindexed to get the corresponding rows/columns of the ``periodic" matrix $S(U)$ as such:
\begin{eqnarray}
\left[ 1, \, \cdots , \, {n-1}, \, {n+1}, \, \cdots, \, {2n-1}, \right.&{2n+1},&\left. \cdots, \, {3n-1}, \, \cdots \cdots,\, {(n-1)n-1} \right] \nonumber\\
& \downarrow &  \label{reindex}\\
\left[ 1, \, \cdots , \, {n-1}, \, {n}, \, \cdots, \, {2n-2}, \right.&{2n-1},&\left. \cdots, \, {3n-3}, \, \cdots \cdots,\, {(n-1)^2} \right] \nonumber
\end{eqnarray}
But  the indices of the $U$ vector are not reindexed in what follows. The entries in the matrix $S(U)$ that will be modified are the ones corresponding to the green left corner, lower and left red boundary nodes, and the upper and left boundary black nodes in Figure \ref{fig:mesh2}.\vspace{2mm}

\noindent \textbf{1- Green Left Node:}
Assuming that the vertices $1$, $n$, $n(n-1) + 1$ and $n^2$ coincide, then the first row of the ``periodic" $S(U)$ will be the sum of the entries in rows  $1$, $n$, $n(n-1) + 1$ and $n^2$  of the general $S(U)$:
\begin{enumerate}
\item[\bf a)] Row $v = n$ of the general matrix $S(U)$ has the entry $\frac{1}{6} U_{ n-1} - \frac{1}{6} U_{ 2n} = \frac{1}{6} U_{ n-1} - \frac{1}{6} U_{ n+1}$ in columns  $n-1,v = n = 1$, and $2n = n+1$. 
\item[\bf b)] Row $v = n(n-1) + 1$ of the general matrix $S(U)$ has the entry $\frac{1}{6} U_{ n(n-1) + 2} - \frac{1}{6} U_{n(n-2)+1} = \frac{1}{6} U_{2} - \frac{1}{6} U_{n(n-2)+1}$ in columns  $n(n-2)+1, v = n(n-1)+1 = 1$, and $n(n-1)+2 = 2$.     
\item[\bf c)] Row $v = 1$ of the general matrix $S(U)$ has nonzero entries in columns  \\
$\begin{array}{lll}
 v=1 : & \dfrac{1}{6} U_{ n+1} - \dfrac{1}{6} U_{ n+2}+\dfrac{1}{6} U_{ n+2} - \dfrac{1}{6} U_{ 2}&=  \dfrac{1}{6} U_{ n+1} -  \dfrac{1}{6} U_{ 2} \nonumber \vspace{3mm}\\
 2 : & \dfrac{1}{6} U_{ n+2} - \dfrac{1}{6} U_{ 2}& = \dfrac{1}{6} U_{ n+2} - \dfrac{1}{6} U_{ 2}  \nonumber \vspace{3mm}\\
 n+1 : &\dfrac{1}{6} U_{ n+1} - \dfrac{1}{6} U_{ n+2} &= \dfrac{1}{6} U_{ n+1} - \dfrac{1}{6} U_{ n+2}  \nonumber \vspace{3mm}\\
  n+2 : & \dfrac{1}{6} U_{ n+1} - \dfrac{1}{6} U_{ n+2}+\dfrac{1}{6} U_{ n+2} - \dfrac{1}{6} U_{ 2}& = \dfrac{1}{6} U_{ n+1} - \dfrac{1}{6} U_{ 2} \nonumber \vspace{1mm}
\end{array}$   
\item[\bf d)] Row $v = n^2$ of the general matrix $S(U)$ has nonzero entries in columns  \\
$\begin{array}{lll}
 n(n-1)-1 : &\dfrac{1}{6} U_{ n(n-1)} - \dfrac{1}{6} U_{ n^2-1}&= \dfrac{1}{6} U_{ n(n-2)+1} - \dfrac{1}{6} U_{ n-1}\nonumber \vspace{3mm}\\
 n(n-1) = n(n-2)+1 : &\dfrac{1}{6} U_{ n(n-1)} - \dfrac{1}{6} U_{ n(n-1)-1}& =   \dfrac{1}{6} U_{ n(n-2)+1} - \dfrac{1}{6} U_{ n(n-1)-1}  \nonumber \vspace{3mm}\\
 n^2-1 = n-1 : &\dfrac{1}{6} U_{ n(n-1)-1} - \dfrac{1}{6} U_{ n^2-1} & =  \dfrac{1}{6} U_{ n(n-1)-1} - \dfrac{1}{6} U_{ n-1}  \nonumber \vspace{3mm}\\
  v=n^2 = 1: & \dfrac{1}{6} U_{ n(n-1)} - \dfrac{1}{6} U_{ n^2-1}& = \dfrac{1}{6} U_{ n(n-2)+1} - \dfrac{1}{6} U_{ n-1}\nonumber
\end{array}$
\end{enumerate}
\noindent Adding up these rows, we get row $1$ of the ``periodic" $S(U)$ with the following nonzero entries:\vspace{1mm}\\
  $\begin{array}{lll}
    1: &  \dfrac{1}{6} U_{n-1} - \dfrac{1}{6} U_{ n+1}+ \dfrac{1}{6} U_{2} - \dfrac{1}{6} U_{n(n-2)+1}+ \dfrac{1}{6} U_{ n+1} -  \dfrac{1}{6} U_{ 2}+ \dfrac{1}{6} U_{ n(n-2)+1} - \dfrac{1}{6} U_{ n-1}&=0\nonumber \vspace{3mm}\\
    2: &\dfrac{1}{6} U_{2} - \dfrac{1}{6} U_{n(n-2)+1} + \dfrac{1}{6} U_{ n+2} - \dfrac{1}{6} U_{ 2} \quad = \quad \dfrac{1}{6} U_{ n+2}- \dfrac{1}{6} U_{n(n-2)+1} &\nonumber \vspace{3mm}\\
    n-1: &\dfrac{1}{6} U_{n-1} - \dfrac{1}{6} U_{n+1} + \dfrac{1}{6} U_{ n(n-1)-1} - \dfrac{1}{6} U_{ n-1} \quad = \quad \dfrac{1}{6} U_{ n(n-1)-1}- \dfrac{1}{6} U_{n+1}\nonumber\vspace{3mm}\\
    n+1: &\dfrac{1}{6} U_{n-1} - \dfrac{1}{6} U_{n+1}+\dfrac{1}{6} U_{ n+1} - \dfrac{1}{6} U_{ n+2} \quad=\quad \dfrac{1}{6} U_{n-1} - \dfrac{1}{6} U_{ n+2}&\nonumber \vspace{3mm}\\
    n+2: &\dfrac{1}{6} U_{ n+1} - \dfrac{1}{6} U_{ 2} &\nonumber \vspace{3mm}\\
    n(n-2)+1: &\dfrac{1}{6} U_{2} - \dfrac{1}{6} U_{n(n-2)+1} + \dfrac{1}{6} U_{ n(n-2)+1} - \dfrac{1}{6} U_{ n(n-1)-1} \quad = \quad \dfrac{1}{6} U_{2}- \dfrac{1}{6} U_{ n(n-1)-1}&\nonumber \vspace{3mm}\\
    n(n-1)-1: & \dfrac{1}{6} U_{ n(n-2)+1} - \dfrac{1}{6} U_{ n-1}&\nonumber \vspace{3mm}\\
  \end{array}$ 
  
  Recall that these column indices have to be reindex by \eqref{reindex}.\vspace{2mm}
  
  \noindent \textbf{2- Lower Red Boundary Nodes:}
  Assuming that the vertices $i$ and ${i+n(n-1)}$   coincide for  $i = 2, 3,..,n-1$, then the rows $i$ will be the sum of the entries in rows $i$ and ${i+n(n-1)}$ of the general $S(U)$:
  \begin{enumerate}
\item[\bf a)] Row $i$ of the general matrix $S(U)$ for $i = 2,3,...,n-1$ has nonzero entries in columns  \vspace{2mm}\\
\begin{minipage}{0.3\textwidth}
$\begin{array}{lll}
 i-1 : &(U_{ i-1 } - U_{i+n})/6 &  \nonumber \vspace{3mm}\\
 v = i : &(U_{ i-1 }- U_{i+1})/6& \nonumber \vspace{1mm}
 \end{array}$
\end{minipage}
\begin{minipage}{0.3\textwidth}
$\begin{array}{lll}
  i+1 : &(U_{i+n+1} - U_{i+1})/6 & \nonumber \vspace{3mm}\\
     i+n : & (U_{ i-1 }- U_{ i+n+1})/6& \nonumber \vspace{1mm}
\end{array}$
 \end{minipage}
 \begin{minipage}{0.3\textwidth}
 $\begin{array}{lll}
   i+n+1 : & \dfrac{1}{6} U_{ i+n} - \dfrac{1}{6}U_{i+1}& \nonumber 
\end{array}$
 \end{minipage}

\item[\bf b)] Row $v = n(n-1)+i$ of the general matrix $S(U)$ for $i = 2,3,...,n-1$, has nonzero entries in columns \\
$\begin{array}{lll}
 n(n-2)+i-1 : &\dfrac{1}{6} U_{ n(n-2)+i } - \dfrac{1}{6} U_{ n(n-1)+i-1} & = \dfrac{1}{6} U_{ n(n-2)+i } - \dfrac{1}{6} U_{i-1}  \nonumber \vspace{3mm}\\
 n(n-2)+i : & \dfrac{1}{6} U_{n(n-1)+i+1} - \dfrac{1}{6}U_{n(n-2)+i-1 } & = \dfrac{1}{6} U_{i+1} - \dfrac{1}{6}U_{n(n-2)+i-1 } \nonumber \vspace{3mm}\\
 n(n-1)+i-1 = i-1: &  \dfrac{1}{6} U_{ n(n-2)+i-1} - \dfrac{1}{6} U_{ n(n-1)+i-1}& = \dfrac{1}{6} U_{ n(n-2)+i-1} - \dfrac{1}{6} U_{i-1} \nonumber \vspace{3mm}\\
 v =  n(n-1)+i = i : & \dfrac{1}{6} U_{n(n-1)+i+1}-\dfrac{1}{6} U_{ n(n-1)+i-1}& = \dfrac{1}{6} U_{i+1}-\dfrac{1}{6} U_{ i-1} \nonumber \vspace{3mm}\\
 n(n-1)+ i+1 = i+1 : &\dfrac{1}{6} U_{n(n-1)+i+1} - \dfrac{1}{6} U_{ n(n-2)+i} & = \dfrac{1}{6} U_{i+1} - \dfrac{1}{6} U_{ n(n-2)+i} \nonumber \vspace{3mm}\\
\end{array}$
\end{enumerate}
Adding up these 2 rows, we get row $i$ of the ``periodic" $S(U)$ with the following nonzero entries for $i = 2,3,...,n-1$:\\

 $\begin{array}{lll}
 i-1:& \dfrac{1}{6} U_{ i-1 } - \dfrac{1}{6}U_{i+n} + \dfrac{1}{6} U_{ n(n-2)+i-1} - \dfrac{1}{6} U_{i-1}& =  \dfrac{1}{6} U_{ n(n-2)+i-1} - \dfrac{1}{6}U_{i+n} \nonumber \vspace{3mm}\\
 i:&\dfrac{1}{6} U_{ i-1 }- \dfrac{1}{6}U_{i+1} +\dfrac{1}{6} U_{i+1}-\dfrac{1}{6} U_{ i-1}& =0 \nonumber \vspace{3mm}\\
 i+1:&\dfrac{1}{6}U_{i+n+1} - \dfrac{1}{6}U_{i+1} +\dfrac{1}{6} U_{i+1} - \dfrac{1}{6} U_{ n(n-2)+i} & = \dfrac{1}{6}U_{i+n+1} - \dfrac{1}{6} U_{ n(n-2)+i} \nonumber \vspace{3mm}\\
 i+n:&\dfrac{1}{6} U_{ i-1 }- \dfrac{1}{6} U_{ i+n+1} & \nonumber \vspace{3mm}\\
 i+n+1:&\dfrac{1}{6} U_{ i+n} - \dfrac{1}{6}U_{i+1} & \nonumber \vspace{3mm}\\
\end{array}$
 
 $\begin{array}{lll} 
 n(n-2)+i-1 :&\dfrac{1}{6} U_{ n(n-2)+i } - \dfrac{1}{6} U_{i-1} & \nonumber \vspace{3mm} \\
 n(n-2)+i:&\dfrac{1}{6} U_{i+1} - \dfrac{1}{6}U_{n(n-2)+i-1} & \nonumber \vspace{3mm}\\
 \end{array}$
 
 Recall that these column indices have to be reindex by \eqref{reindex}.
 Note that for $i=n-1$ we get the following nonzero entries:\\
 
  $\begin{array}{lll}
 n-2:&  \dfrac{1}{6} U_{ n(n-2)+n-2} - \dfrac{1}{6}U_{2n-1} & \nonumber \vspace{3mm}\\
 n = 1:& \dfrac{1}{6}U_{2n} - \dfrac{1}{6} U_{ n(n-2)+n-1}& = \dfrac{1}{6}U_{n+1} - \dfrac{1}{6} U_{ n(n-2)+n-1}  \nonumber \vspace{3mm}\\
 2n-1:&\dfrac{1}{6} U_{ n-2 }- \dfrac{1}{6} U_{ 2n} &= \dfrac{1}{6} U_{ n-2 }- \dfrac{1}{6} U_{ n+1} \nonumber \vspace{3mm}\\
 2n = n+1:&\dfrac{1}{6} U_{ 2n-1} - \dfrac{1}{6}U_{n} & = \dfrac{1}{6} U_{ 2n-1} - \dfrac{1}{6}U_{1} \nonumber \vspace{3mm}\\
 n(n-2)+n-2 :&\dfrac{1}{6} U_{ n(n-2)+n-1 } - \dfrac{1}{6} U_{n-2} & \nonumber \vspace{3mm} \\
 n(n-2)+n-1:&\dfrac{1}{6} U_{n} - \dfrac{1}{6}U_{n(n-2)+n-2} & = \dfrac{1}{6} U_{1} - \dfrac{1}{6}U_{n(n-2)+n-2} \nonumber \vspace{3mm}\\
 \end{array}$
  
  \noindent \textbf{3- Left Red Boundary Nodes:}
  Assuming that the vertices $kn+1$ and $kn+n$ coincide for  $k = 1, 2,..,n-2$, then the corresponding rows $kn+1-k$ of the ``periodic" $S(U)$  will be the sum of the entries in rows $kn+1$ and $kn+n$ of the general $S(U)$:
  \begin{enumerate}
\item[\bf a)] Row $kn+1$ of the general matrix $S(U)$ for $k = 1, 2,..,n-2$,  has nonzero entries in columns  \\
$nk-n+1 :  \quad ( U_{ nk+2} - U_{ nk-n+1})/6$

\begin{minipage}{0.45\textwidth}
$\begin{array}{lll}
 nk+1 : & \dfrac{1}{6} U_{ nk+n+1} - \dfrac{1}{6} U_{ nk-n+1} &\nonumber \vspace{3mm}\\
  nk+2 : &  \dfrac{1}{6} U_{ nk+n +2 } - \dfrac{1}{6} U_{ nk-n+1}&  \nonumber 
\end{array}$
\end{minipage}
\begin{minipage}{0.45\textwidth}
$\begin{array}{lll}
  nk+n+1 : & \dfrac{1}{6} U_{ nk+n+1} - \dfrac{1}{6} U_{ nk+n+2}& \nonumber \vspace{3mm}\\
   nk+n+2 : & \dfrac{1}{6} U_{ nk+n+1} - \dfrac{1}{6}U_{nk+2} &\nonumber
\end{array}$
\end{minipage}

\item[\bf b)] Row $(k+1)n$ of the general matrix $S(U)$ for $k = 1,2,...,n-2$, has nonzero entries in columns \\
$\begin{array}{lll}
 n(k+1)-n-1 = nk-1 : & \dfrac{1}{6} U_{ n(k+1)-n} - \dfrac{1}{6} U_{ n(k+1)-1}& = \dfrac{1}{6} U_{ nk-n+1} - \dfrac{1}{6} U_{ nk+n-1}  \nonumber \vspace{3mm}\\
 n(k+1)-n = nk = nk-n+1 : &\dfrac{1}{6} U_{ n(k+1)-n }- \dfrac{1}{6}U_{n(k+1)-n-1} &=  \dfrac{1}{6} U_{ nk-n+1 } - \dfrac{1}{6}U_{nk-1}  \nonumber \vspace{3mm}\\
 n(k+1)-1 = nk+n-1 : & \dfrac{1}{6} U_{ n(k+1)-n-1}  - \dfrac{1}{6}U_{n(k+1)+n}& = \dfrac{1}{6} U_{ nk-1}  - \dfrac{1}{6}U_{nk+n+1} \nonumber \vspace{3mm}\\
 n(k+1) = nk+1: & \dfrac{1}{6} U_{ n(k+1)-n} - \dfrac{1}{6} U_{ n(k+1)+n} &= \dfrac{1}{6} U_{ nk-n+1} - \dfrac{1}{6} U_{ nk+n+1} \nonumber \vspace{3mm}\\
  n(k+1)+n = nk+n+1  : &\dfrac{1}{6} U_{ n(k+1)-1 } - \dfrac{1}{6}U_{n(k+1)+n} & = \dfrac{1}{6} U_{ nk+n-1 } - \dfrac{1}{6}U_{nk+n+1}\nonumber \vspace{3mm}\\
\end{array}$
\end{enumerate}
Adding up these 2 rows, we get row $nk+1-k$ of the ``periodic" $S(U)$ with the following nonzero entries for $k = 1,2,...,n-2$:\\

 $\begin{array}{lll}
 nk-n+1 :&\dfrac{1}{6} U_{ nk+2} - \dfrac{1}{6} U_{ nk-n+1} +\dfrac{1}{6} U_{ nk-n+1 } - \dfrac{1}{6}U_{nk-1} & = \dfrac{1}{6} U_{ nk+2}  - \dfrac{1}{6}U_{nk-1} \vspace{3mm}\\
 nk-1:&\dfrac{1}{6} U_{ nk-n+1} - \dfrac{1}{6} U_{ nk+n-1} &  \vspace{3mm}\\
 \end{array}$
 
 $\begin{array}{lll}
 nk+1:& \dfrac{1}{6} U_{ nk+n+1} - \dfrac{1}{6} U_{ nk-n+1} + \dfrac{1}{6} U_{ nk-n+1} - \dfrac{1}{6} U_{ nk+n+1}& = 0\vspace{3mm}\\
 nk+2:&  \dfrac{1}{6} U_{ nk+n +2 } - \dfrac{1}{6} U_{ nk-n+1} &  \vspace{3mm}\\
 nk+n-1:&\dfrac{1}{6} U_{ nk-1}  - \dfrac{1}{6}U_{nk+n+1} & \vspace{3mm}\\
 nk+n+1:&\dfrac{1}{6} U_{ nk+n+1} - \dfrac{1}{6} U_{ nk+n+2} + \dfrac{1}{6} U_{ nk+n-1 } - \dfrac{1}{6}U_{nk+n+1} & = \dfrac{1}{6} U_{ nk+n-1 } - \dfrac{1}{6} U_{ nk+n+2} \vspace{3mm}\\
 nk+n+2:&\dfrac{1}{6} U_{ nk+n+1} - \dfrac{1}{6}U_{nk+2}  &  \vspace{3mm}\\
 \end{array}$
 
 \noindent Recall that these column indices have to be reindex by \eqref{reindex}.
 Note that for $k=n-2$ we get the following nonzero entries:\\
  $\begin{array}{lll}
 n(n-3)+1 :&\dfrac{1}{6} U_{ n(n-2)+2}  - \dfrac{1}{6}U_{n(n-2)-1} \vspace{3mm}\\
 n(n-2)-1:&\dfrac{1}{6} U_{ n(n-3)+1} - \dfrac{1}{6} U_{ n(n-1)-1} &  \vspace{3mm}\\
 n(n-2)+2:&  \dfrac{1}{6} U_{ n(n-1) +2 } - \dfrac{1}{6} U_{ n(n-3)+1} & = \dfrac{1}{6} U_{2 } - \dfrac{1}{6} U_{ n(n-3)+1} \vspace{3mm}\\
 n(n-1)-1:&\dfrac{1}{6} U_{ n(n-2)-1}  - \dfrac{1}{6}U_{n(n-1)+1} & = \dfrac{1}{6} U_{ n(n-2)-1}  - \dfrac{1}{6}U_{1}  \vspace{3mm}\\
 n(n-1)+1 = 1:&  \dfrac{1}{6} U_{ n(n-1)-1 } - \dfrac{1}{6} U_{ n(n-1)+2} &= \dfrac{1}{6} U_{ n(n-1)-1 } - \dfrac{1}{6} U_{2} \vspace{3mm}\\
 n(n-1)+2 = 2:&\dfrac{1}{6} U_{ n(n-1)+1} - \dfrac{1}{6}U_{n(n-2)+2}  &= \dfrac{1}{6} U_{1} - \dfrac{1}{6}U_{n(n-2)+2}  \vspace{3mm}\\
 \end{array}$
 

  \noindent \textbf{4- Right Black Boundary Nodes:}
  Rows $(k+1)n-k-1$  of the ``periodic" $S(U)$ matrix correspond to rows $(k+1)n-1$ of the general $S(U)$ for  $k = 1,..,n-3$ with nonzero entries  in columns:\\  
$\begin{array}{lll}
nk-n+n-1-1 = nk-2 : &  \dfrac{1}{6} U_{nk-n+n-1} - \dfrac{1}{6} U_{nk+n-1-1}&= \dfrac{1}{6} U_{nk-1} - \dfrac{1}{6} U_{nk+n-2} \nonumber \vspace{3mm}\\
nk-n+n-1 = nk-1 : &\dfrac{1}{6}U_{nk+n-1+1}- \dfrac{1}{6} U_{nk-n+n-1-1}&=\dfrac{1}{6}U_{nk+1}- \dfrac{1}{6} U_{nk-2}  \nonumber \vspace{3mm}\\
 nk+n-1-1 = nk+n-2 : & \dfrac{1}{6} U_{nk-n+n-1-1} - \dfrac{1}{6}U_{nk+n-1+n}&= \dfrac{1}{6} U_{nk-2} - \dfrac{1}{6}U_{nk+2n-1}\nonumber \vspace{3mm}\\
  nk+n-1+1=nk+n = nk+1 : &  \dfrac{1}{6} U_{ nk+n-1+n +1 }- \dfrac{1}{6}U_{nk-n+n-1}& = \dfrac{1}{6} U_{ nk+n+1}- \dfrac{1}{6}U_{nk-1}\nonumber \vspace{3mm}\\
  nk+n-1+n=nk+2n-1 : &  \dfrac{1}{6} U_{nk+n-1-1}- \dfrac{1}{6} U_{ nk+n-1 +n+1}&=  \dfrac{1}{6} U_{nk+n-2}- \dfrac{1}{6} U_{ nk+n+1} \nonumber \vspace{3mm}\\
   nk+n-1+n+1 = nk+2n = nk+n+1: & \dfrac{1}{6} U_{nk+n-1+n}- \dfrac{1}{6}U_{nk+n-1+1}&= \dfrac{1}{6} U_{nk+2n-1}- \dfrac{1}{6}U_{nk+1}\nonumber
\end{array}$

Recall that these column indices have to be reindex by \eqref{reindex}.

   \noindent \textbf{5- Upper Black Boundary Nodes:}
  Rows $(n-1)(n-2)+i$  of the ``periodic" $S(U)$ matrix correspond to rows $n(n-2)+i$ of the general $S(U)$ for  $i = 2,3,..,n-1$ with nonzero entries  in columns:\\ 
  $\begin{array}{lll}
n(n-2)-n+i-1 = n(n-3)+i-1   : & \dfrac{1}{6} U_{n(n-2)-n+i} - \dfrac{1}{6} U_{n(n-2)+i-1}&=\dfrac{1}{6} U_{n(n-3)+i} - \dfrac{1}{6} U_{n(n-2)+i-1} \nonumber \vspace{3mm}\\
n(n-2)-n+i  = n(n-3)+i : & \dfrac{1}{6}U_{n(n-2)+i+1}- \dfrac{1}{6} U_{n(n-2)-n+i-1} &= \dfrac{1}{6}U_{n(n-2)+i+1}- \dfrac{1}{6} U_{n(n-3)+i-1}  \nonumber \vspace{3mm}\\
\end{array}$
 
 $\begin{array}{lll}
 n(n-2)+i-1 : &\dfrac{1}{6} U_{n(n-2)-n+i-1} - \dfrac{1}{6}U_{n(n-2)+i+n}& = \dfrac{1}{6} U_{n(n-3)+i-1} - \dfrac{1}{6}U_{i} \nonumber \vspace{3mm}\\ 
  n(n-2)+i+1 : &   \dfrac{1}{6} U_{ n(n-2)+i+n +1 }- \dfrac{1}{6}U_{n(n-2)-n+i} & =   \dfrac{1}{6} U_{ i +1 }- \dfrac{1}{6}U_{n(n-3)+i} 
  \nonumber \vspace{3mm}\\
  n(n-2)+i+n =  n(n-1)+i = i : &  \dfrac{1}{6} U_{n(n-2)+i-1}- \dfrac{1}{6} U_{ n(n-2)+i +n+1}  &=  \dfrac{1}{6} U_{n(n-2)+i-1}- \dfrac{1}{6} U_{i+1}
   \nonumber \vspace{3mm}\\
   n(n-2)+i+n+1 = n(n-1)+i+1 = i+1 : & \dfrac{1}{6} U_{n(n-2)+i+n}- \dfrac{1}{6}U_{n(n-2)+i+1}& =\dfrac{1}{6} U_{i}- \dfrac{1}{6}U_{n(n-2)+i+1}
   \nonumber
\end{array}$

  Recall that these column indices have to be reindex by \eqref{reindex}.
 Note that for $i=n-1$ we get the following nonzero entries:\\
 \noindent$\begin{array}{lll}
 n(n-3)+n-1-1 = n(n-2)-2: &\dfrac{1}{6} U_{n(n-3)+n-1} - \dfrac{1}{6} U_{n(n-2)+n-1-1} &=\dfrac{1}{6} U_{n(n-2)-1} - \dfrac{1}{6} U_{n(n-1)-2} \nonumber \vspace{3mm}\\
n(n-3)+n-1 = n(n-2)-1: & \dfrac{1}{6}U_{n(n-2)+n-1+1}- \dfrac{1}{6} U_{n(n-3)+n-1-1}   & =\dfrac{1}{6}U_{n(n-2)+1}- \dfrac{1}{6} U_{n(n-2)-2}  \nonumber \vspace{3mm}\\
 n(n-2)+n-1-1= n(n-1)-2 : & \dfrac{1}{6} U_{n(n-3)+n-1-1} - \dfrac{1}{6}U_{n-1} & = \dfrac{1}{6} U_{n(n-2)-2} - \dfrac{1}{6}U_{n-1} \nonumber \vspace{3mm}\\
  n(n-2)+n-1+1 = n(n-1) = n(n-2)+1  : &   \dfrac{1}{6} U_{ n-1 +1 }- \dfrac{1}{6}U_{n(n-3)+n-1}  &  = \dfrac{1}{6} U_{ 1 }- \dfrac{1}{6}U_{n(n-2)-1} 
  \nonumber \vspace{3mm}\\
   n-1: &  \dfrac{1}{6} U_{n(n-2)+n-1-1}- \dfrac{1}{6} U_{n-1+1} &  =\dfrac{1}{6} U_{n(n-1)-2}- \dfrac{1}{6} U_{1}
   \nonumber \vspace{3mm}\\
  n-1+1 = n = 1 : & \dfrac{1}{6} U_{n-1}- \dfrac{1}{6}U_{n(n-2)+n-1+1} & = \dfrac{1}{6} U_{n-1}- \dfrac{1}{6}U_{n(n-2)+1}
   \nonumber
\end{array}$

 \noindent For $n=5$, the matrix $S(U)$ corresponding to the mesh in Figure \ref{fig:mesh2}, 
 has the following sparsity pattern:
  $$S(U) = \left[
\begin{array}{cccc|cccc|cccc|cccc}
0&*&0&*&*&*&0&0&0&0&0&0&*&0&0&*\\
*&0&*&0&0&*&*&0&0&0&0&0&*&*&0&0\\
0&*&0&*&0&0&*&*&0&0&0&0&0&*&*&0\\
*&0&*&0&*&0&0&*&0&0&0&0&0&0&*&*\\
\hline
*&0&0&*&0&*&0&*&*&*&0&0&0&0&0&0\\
*&*&0&0&*&0&*&0&0&*&*&0&0&0&0&0\\
0&*&*&0&0&*&0&*&0&0&*&*&0&0&0&0\\
0&0&*&*&*&0&*&0&*&0&0&*&0&0&0&0\\
 \hline
0&0&0&0&*&0&0&*&0&*&0&*&*&*&0&0\\
0&0&0&0&*&*&0&0&*&0&*&0&0&*&*&0\\
0&0&0&0&0&*&*&0&0&*&0&*&0&0&*&*\\
0&0&0&0&0&0&*&*&*&0&*&0&*&0&0&*\\
\hline 
*&*&0&0&0&0&0&0&*&0&0&*&0&*&0&*\\
0&*&*&0&0&0&0&0&*&*&0&0&*&0&*&0\\
0&0&*&*&0&0&0&0&0&*&*&0&0&*&0&*\\
*&0&0&*&0&0&0&0&0&0&*&*&*&0&*&0\\

\end{array} \right]
$$ with $U = [U_1, U_2, U_3, U_4, U_6,U_7, U_8, U_9, U_{11}, U_{12}, U_{13}, U_{14}, U_{16},U_{17},U_{18}, U_{19}]$

 \noindent In general, the $S(U)$ matrix is an $(n-1)^2 \times (n-1)^2$ block tridiagonal matrix with 2 additional blocks in the upper right  and lower left corner. Moreover, it is a skew-symmetric matrix ($A^T = -A$) with  6 nonzero entries per row, 6 nonzero entries per column, and zeros on the diagonal.
 
$$S(U) =\dfrac{1}{6} \left[ \begin{array}{cccccc} A_{1,1} & A_{1,2}& 0&\cdots&0&A_{1,n-1}\\
A_{2,1} & A_{2,2}& A_{2,3}&0&\cdots&0\\
0&\ddots&\ddots&\ddots&\ddots&\vdots\\
\vdots&\ddots&A_{j,h}&A_{j,j}&A_{j,i}&0\\
0&\cdots&0&A_{i,j}&A_{i,i}&A_{i,k}\\
A_{k,1}&0&\cdots&0&A_{k,i}&A_{k,k}\\
\end{array}\right]
$$
where $i = n-2, j = n-3, k= n-1, h=n-4$, and the $3(n-1)$ nonzero block matrices are of size $(n-1)\times (n-1)$ with $2(n-1)$ nonzero entries each, and the following sparsity patterns:
\begin{itemize}
\item $A_{i,i}$ for $i = 1,2,.., n-1$ are tridiagonal matrices with zero diagonal entries, and nonzero $A_{i,i}(1,n-1), A_{i,i}(n-1,1)$.
\item $A_{1,n-1}$ and $A_{i+1,i}$ for $i = 1,2,3,.., n-2$ are lower bidiagonal matrices, with nonzero entry in first row, column $n-1$ . 
\item  $A_{n-1,1}$ and $A_{i,i+1}$ for $i = 1,3,.., n-2$ are upper bidiagonal matrices with nonzero entry in first column, row $n-1$.
\end{itemize}
Thus, $S(U)$ has a total of $3(n-1)2(n-1) = 6(n-1)^2$ nonzero entries.
 
 \noindent As for the nonzero entries in $S(U)$, each is of order $||U||_{\infty}$ since each nonzero entry is of the form $$\dfrac{1}{6}(U_i-U_j) \leq \dfrac{1}{6}|U_i-U_j| \leq \dfrac{1}{6}(|U_i|+|U_j|) \leq \dfrac{1}{3} ||U||_{\infty}$$

\subsection{Assembling the Global matrix $R$ }\label{sec:R}

First, we derive the sparsity pattern of the global matrix $R$ without any assumptions related to boundary conditions, i.e. define the general matrix $R$, in Section  \ref{sec:GR}. Then, we  derive the sparsity pattern of the matrix $R$  assuming periodic boundary conditions in Section \ref{sec:PR}.

\subsubsection{The General Matrix $R$} \label{sec:GR}
Given the meshing of Figure \ref{fig:mesh2} and without any assumptions on the boundary nodes, $R$ is an $n^2\times n^2$ sparse matrix with at most 6 nonzero entries per row, as discussed below. To get those  nonzero entries, note that each node has at most 6 edges connecting it with its neighboring nodes, i.e. belongs to at most 6 triangles. We consider the 4 types of nodes shown in different colors in Figure \ref{fig:mesh2}: the $(n-2)^2$ black internal nodes that belong to 6 triangles, the $4(n-2)$ red boundary nodes that belong to 4 triangles, the 2 green corner nodes that belong to 2 triangles, and the 2 blue corner nodes that belong to 1 triangle.\vspace{5mm}

\noindent \textbf{1- Black Internal Nodes:}

Each of the $(n-2)^2$ black internal nodes with index $v = ln+i$ for $l = 1,2,...,n-2$ and $i=2,3,...,n-1$  belongs to triangles $T_{2j+1}, T_{2j+2}, T_{2j+3}, T_{2(j+n)}, T_{2(j+n)+1}, T_{2(j+n)+2}$,
where  $j = (n-1)(l-1)+(i-2)$. 

Thus, we first define the nonzero entries in row $v$ per local triangle and then add them up. \begin{itemize}
\item Triangle $T_{2j+1}$ of type a, has node $v$ as the second local vertex.
 Thus row $v$ of the matrix $R$ has the entry $0$ in columns  $nl-n+i-1, nl+i-1$, and $v = nl+i$. 
 \item Triangle $T_{2j+2}$ of type b, has node $v$ as the third local vertex. 
 Thus row $v$ of the matrix $R$ has the entry $\dfrac{h}{6}\;\hat{k}$ in columns  $nl-n+i-1, nl-n+i$, and $v = nl+i$. 
  \item Triangle $T_{2j+3}$ of type a, has node $v$ as the third local vertex. 
 Thus row $v$ of the matrix $R$ has the entry $\dfrac{h}{6}\;\hat{k}$ in columns  $nl-n+i, v = nl+i$, and $ nl+i+1$. 
  \item Triangle $T_{2(j+n)}$ of type b, has node $v$ as the second local vertex . 
 Thus row $v$ of the matrix $R$ has the entry $-\dfrac{h}{6}\;\hat{k}$ in columns  $nl+i-1, v = nl+i$, and $nl+i+n $. 
 \item Triangle $T_{2(j+n)+1}$  of type a, has node $v$ as the first local vertex. 
 Thus row $v$ of the matrix $R$ has the entry $-\dfrac{h}{6}\;\hat{k}$ in columns  $v = nl+i, nl+i+n$, and $nl+i+n+1$. 
  \item Triangle $T_{2(j+n)+2}$  of type b, has node $v$ as the first local vertex. 
 Thus row $v$ of the matrix $R$ has the entry $0$ in columns  $v = nl+i, nl+i+1 $, and $nl+i+n+1 $. 
\end{itemize}
Thus the nonzero entries in row $v = nl+i$ of $R$ for $l = 1,2,...,n-2$ and $i=2,3,...,n-1$ are in columns  \\
$\begin{array}{lll}
nl-n+i-1 : & \dfrac{h}{6}\;\hat{k}  &=\dfrac{h}{6}\;\hat{k}\nonumber \vspace{3mm}\\
nl-n+i : & \dfrac{h}{6}\;\hat{k} + \dfrac{h}{6}\;\hat{k}&= \dfrac{h}{3}\;\hat{k} \nonumber \vspace{3mm}\\
 nl+i-1 : & -\dfrac{h}{6}\;\hat{k} &=-\dfrac{h}{6}\;\hat{k}\nonumber \vspace{3mm}\\
 v = nl+i : &\dfrac{h}{6}\;\hat{k}+ \dfrac{h}{6}\;\hat{k} - \dfrac{h}{6}\;\hat{k} -\dfrac{h}{6}\;\hat{k} &=0\nonumber \vspace{3mm}\\
  nl+i+1 : & \dfrac{h}{6}\;\hat{k}&=\dfrac{h}{6}\;\hat{k} \nonumber \vspace{3mm}\\
  nl+i+n : & -\dfrac{h}{6}\;\hat{k}-\dfrac{h}{6}\;\hat{k}& = -\dfrac{h}{3}\;\hat{k} \nonumber \vspace{3mm}\\
   nl+i+n+1 : &-\dfrac{h}{6}\;\hat{k} &=-\dfrac{h}{6}\;\hat{k} \nonumber
\end{array}$\vspace{2mm}

\noindent \textbf{2- Red Boundary Nodes:}

The red boundary nodes are of 4 types:
\begin{enumerate}
\item[\bf a)] the left boundary with index $v = ln+1$ and $l = 1,2,...,n-2$, that belong to triangles\\ $T_{2(l-1)(n-1)+1}, T_{2l(n-1)+1}$, and $T_{2l(n-1)+2}$.\\
Thus, we first define the nonzero entries in row $v=ln+1$ per local triangle and then add them up. 
\begin{itemize}
\item Triangle $T_{2(l-1)(n-1)+1}$ of type a, has node $v$ as the third local vertex. 
 Thus row $v$ of the matrix $R$ has the entry $\dfrac{h}{6}\;\hat{k}$ in columns  $nl-n+1, v = nl+1$, and $nl+2$. 
 \item Triangle $T_{2l(n-1)+1}$ of type a, has node $v$ as the first local vertex. 
 Thus row $v$ of the matrix $R$ has the entry $-\dfrac{h}{6}\;\hat{k}$ in columns  $v = nl+1, nl+n+1, $, and $nl+n+2$. 
 \item Triangle $T_{2l(n-1)+2}$ has node $v$ is the first local vertex. 
 Thus row $v$ of the matrix $R$ has the entry $0$ in columns  $v = nl+1, nl+2 $, and $nl+n+2 $. 
\end{itemize}
Thus the nonzero entries in row $v = nl+1$ of $R$ for $l = 1,2,...,n-2$ are in columns  \\
$\begin{array}{lll}
 nl-n+1 : & \dfrac{h}{6}\;\hat{k}&= \dfrac{h}{6}\;\hat{k} \nonumber \vspace{3mm}\\
 v = nl+1 : & \dfrac{h}{6}\;\hat{k}-\dfrac{h}{6}\;\hat{k}&= 0\nonumber \vspace{3mm}\\
  nl+2 : &\dfrac{h}{6}\;\hat{k}& =\dfrac{h}{6}\;\hat{k} \nonumber \vspace{3mm}\\
  nl+n+1 : & -\dfrac{h}{6}\;\hat{k}& = -\dfrac{h}{6}\;\hat{k}\nonumber \vspace{3mm}\\
   nl+n+2 : &-\dfrac{h}{6}\;\hat{k}& = -\dfrac{h}{6}\;\hat{k}\nonumber \vspace{-1mm}
\end{array}$
\item[\bf b)] the right boundary with index $v = ln$ and $l = 2,3,...,n-1$, that belong to triangles\\ $T_{2(l-1)(n-1)}, T_{2(l-1)(n-1)-1}$, and $T_{2l(n-1)}$.\\
Thus, we first define the nonzero entries in row $v=ln$ per local triangle and then add them up. 
\begin{itemize}
\item Triangle $T_{2(l-1)(n-1)-1}$ of type a, has node $v$ is the second local vertex. 
 Thus row $v$ of the matrix $R$ has the entry $ 0$ in columns  $nl-n-1, nl-1$, and $ v = nl$. 
  \item Triangle $T_{2(l-1)(n-1)}$ of type b, has node $v$ as the third local vertex. 
 Thus row $v$ of the matrix $R$ has the entry $\dfrac{h}{6}\;\hat{k}$ in columns  $nl-n-1, nl-n $, and $v = nl $. 
  \item Triangle $T_{2l(n-1)} $ of type b, has node $v$ is the second local vertex . 
 Thus row $v$ of the matrix $R$ has the entry $-\dfrac{h}{6}\;\hat{k}$ in columns  $ nl-1, v = nl$, and $nl+n $. 
\end{itemize}
Thus the nonzero entries in row $v = nl$ of $R$ for $l = 2,3,...,n-1$ are in columns  \\
$\begin{array}{lll}
 nl-n-1 : &\dfrac{h}{6}\;\hat{k} &=  \dfrac{h}{6}\;\hat{k}\nonumber \vspace{3mm}\\
 nl-n : &\dfrac{h}{6}\;\hat{k} &= \dfrac{h}{6}\;\hat{k} \nonumber \vspace{3mm}\\
 nl-1 : &-\dfrac{h}{6}\;\hat{k} &= -\dfrac{h}{6}\;\hat{k}\nonumber \vspace{3mm}\\
 v = nl : &\dfrac{h}{6}\;\hat{k}-\dfrac{h}{6}\;\hat{k}&=0  \nonumber \vspace{3mm}\\
  nl+n : &-\dfrac{h}{6}\;\hat{k}& = -\dfrac{h}{6}\;\hat{k}\nonumber
\end{array}$
\item[\bf c)] the lower boundary with index $v = i$ and $i = 2,3,...,n-1$, that belong to triangles\\ $T_{2(i-1)}, T_{2(i-1)+1}$, and $T_{2(i-1)+2}$.\\
Thus, we first define the nonzero entries in row $v=i$ per local triangle and then add them up. 
\begin{itemize}
 \item Triangle $T_{2(i-1)}$ of type b, has node $v = i$ as the second local vertex. 
 Thus row $v = i$ of the matrix $R$ has the entry $-\dfrac{h}{6}\;\hat{k}$ in columns  $i-1, v=i $, and $i+ n $. 
\item Triangle $T_{2(i-1)+1}$ of type a, has node $v = i$ as the first local vertex.
 Thus row $v = i$ of the matrix $R$ has the entry $-\dfrac{h}{6}\;\hat{k} $ in columns  $ v = i, i+n$, and $ i+n+1$. 
 \item Triangle $T_{2(i-1) +2}$ of type b, has node $v = i$ as the first local vertex. 
 Thus row $v = i$ of the matrix $R$ has the entry $ 0$ in columns  $v=i, i+1 $, and $i+ n+1$. 
\end{itemize}
Thus the nonzero entries in row $v = i$ of $R$ for $i = 2,3,...,n-1$ are in columns \\ 
$\begin{array}{lll}
 i-1 : &-\dfrac{h}{6}\;\hat{k} &= -\dfrac{h}{6}\;\hat{k}   \nonumber \vspace{3mm}\\
 v = i : &- \dfrac{h}{6}\;\hat{k}-\dfrac{h}{6}\;\hat{k}&= -\dfrac{h}{3}\;\hat{k}\nonumber \vspace{3mm}\\
     i+n : & -\dfrac{h}{6}\;\hat{k} -\dfrac{h}{6}\;\hat{k}&=-\dfrac{h}{3}\;\hat{k} \nonumber \vspace{3mm}\\
   i+n+1 : &-\dfrac{h}{6}\;\hat{k} &= -\dfrac{h}{6}\;\hat{k} \nonumber \vspace{3mm}\\
\end{array}$
\item[\bf d)] the upper boundary with index $v = n(n-1)+i$ and $i = 2,3,...,n-1$, that belong to triangles\\ $T_{2(n-1)(n-2)+2(i-1)-1}, T_{2(n-1)(n-2)+2(i-1)}$, and $T_{2(n-1)(n-2)+2(i-1)+1}$.
\\
Thus, we first define the nonzero entries in row $v=n(n-1)+i$ per local triangle and then add them up. 
\begin{itemize}
\item Triangle $T_{2(n-1)(n-2)+2(i-1)-1}$ of type a, has node $v =  n(n-1)+i$ as the second local vertex. 
 Thus row $v = n(n-1)+i$ of the matrix $R$ has the entry $0$ in columns  $ n(n-2)+i-1,  n(n-1)+i-1$, and $v = n(n-1)+i$. 
 \item Triangle $T_{2(n-1)(n-2)+2(i-1)} $of type b, has node $v = n^2-n+i$ as the third local vertex.
 Thus row $v = i$ of the matrix $R$ has the entry $\dfrac{h}{6}\;\hat{k}$ in columns  $n(n-2)+i-1 ,n(n-2)+i $, and $v = n(n-1)+i$. 
\item Triangle $T_{2(n-1)(n-2)+2(i-1)+1}$ of type a, has node $v =  n(n-1)+i$ as the third local vertex. 
 Thus row $v = n(n-1)+i$ of the matrix $R$ has the entry $\dfrac{h}{6}\;\hat{k}$ in columns  $ n(n-2)+i, v = n(n-1)+i$, and $n(n-1)+i+1$. 
\end{itemize}
Thus the nonzero entries in row $v = n(n-1)+i$ of $R$ for $i = 2,3,...,n-1$, are in columns \\ 
$\begin{array}{lll}
 n(n-2)+i-1 : &\dfrac{h}{6}\;\hat{k} & = \dfrac{h}{6}\;\hat{k}   \nonumber \vspace{3mm}\\
 n(n-2)+i : &\dfrac{h}{6}\;\hat{k}+\dfrac{h}{6}\;\hat{k} & = \dfrac{h}{3}\;\hat{k}   \nonumber \vspace{3mm}\\
 v =  n(n-1)+i : & \dfrac{h}{6}\;\hat{k}+\dfrac{h}{6}\;\hat{k} &=\dfrac{h}{3}\;\hat{k} \nonumber \vspace{3mm}\\
 n(n-1)+ i+1 : & \dfrac{h}{6}\;\hat{k}& = \dfrac{h}{6}\;\hat{k}\nonumber \vspace{3mm}\\
\end{array}$
\end{enumerate} 

\noindent \textbf{3- Green Corner Nodes:}

There are 2 green corners with index :
\begin{enumerate}
\item[\bf a)] $v = 1$ that belongs to triangles $T_{1}$, and $T_{2}$. First, we define the nonzero entries in row $v=1$ per local triangle.
\begin{itemize}\item Triangle $T_{1}$ of type a, has node $v =  1$ is the first local vertex. 
 Thus row $v = 1$ of the matrix $R$ has the entry $-\dfrac{h}{6}\;\hat{k}$ in columns  $v = 1, n+1 $, and $n+2$. 
 \item Triangle $T_{2}$ of type b, has node $v =  1$ as the first local vertex of triangle $T_{1}$. 
 Thus row $v = 1$ of the matrix $R$ has the entry $0$ in columns  $v = 1,  2$, and $n+2$. 
\end{itemize}
Thus the nonzero entries in row $v = 1$ of $R$  are in columns  \\
$\begin{array}{lll}
 v=1 : &-\dfrac{h}{6}\;\hat{k}& \nonumber \vspace{3mm}\\
 n+1 : &-\dfrac{h}{6}\;\hat{k}& \nonumber \vspace{3mm}\\
  n+2 : & -\dfrac{h}{6}\;\hat{k}& \nonumber \vspace{3mm}\\
\end{array}$
\item[\bf b)] $v = n^2$ that belongs to triangles $T_{2(n-1)(n-1)- 1}$, and $T_{2(n-1)(n-1)}$.\\ First, we define the nonzero entries in row $v=n^2$ per local triangle and then add them up. 
\begin{itemize}\item Triangle $T_{2(n-1)(n-1)-1}$ of type a, has node $v =  n^2$ as the second local vertex of triangle $T_{2(n-1)(n-1)-1}$. 
 Thus row $v = n^2$ of the matrix $R$ has the entry $ 0$ in columns  $n(n-1)-1, n^2-1$, and $v = n^2$. 
 \item Triangle $T_{2(n-1)(n-1)}$ of type n, has node $v =  n^2$ as the third local vertex. 
 Thus row $v = n^2$ of the matrix $R$ has the entry $\dfrac{h}{6}\;\hat{k}$ in columns  $n(n-1)-1, n(n-1)$, and $v=n^2$. 
\end{itemize}
Thus the nonzero entries in row $v = n^2$ of $R$  are in columns  \\
$\begin{array}{lll}
 n(n-1)-1 : &\dfrac{h}{6}\;\hat{k}&\nonumber \vspace{3mm}\\
 n(n-1) : &\dfrac{h}{6}\;\hat{k}&  \nonumber \vspace{3mm}\\
  v=n^2 : &\dfrac{h}{6}\;\hat{k}& \nonumber
\end{array}$
\end{enumerate}

\noindent \textbf{4- Blue Corner Nodes:}

There are 2 blue corners with index :
\begin{enumerate}
\item[\bf a)] $v = n$ that belongs to triangle $T_{2(n-1)}$ of type b where node $v =  n$ is the second local vertex. 
 Thus row $v = n$ of the matrix $R$ has the entry $-\dfrac{h}{6}\;\hat{k}$ in columns  $n-1,v = n$, and $2n$. 
\item[\bf b)] $v = n(n-1)+1$ that belongs to triangle $T_{2(n-2)(n-1)+ 1}$ of type a,
 where node $v =  n(n-1) + 1$ is the third local vertex. 
 Thus row $v = n(n-1) + 1$ of the matrix $R$ has the entry $\dfrac{h}{6}\;\hat{k}$ in columns  $n(n-2)+1,\\ v = n(n-1)+1$, and $n(n-1)+2$. 
\end{enumerate}
\noindent For $n=5$, the matrix $R$ corresponding to the mesh in Figure \ref{fig:mesh2}, is given by $  \dfrac{h}{6} \hat{k} A$ where  $A = $ 
 
$\hspace{-25mm}\left[
\begin{array}{ccccc|ccccc|ccccc|ccccc|ccccc}
-1 & 0 & 0&0&0& -1 &-1&0&0&0&0&0&0&0&0&0&0&0&0&0&0&0&0&0&0\\
-1&-2&0&0&0&0&-2&-1&0&0&0&0&0&0&0&0&0&0&0&0&0&0&0&0&0\\
0&-1&-2&0&0&0&0&-2&-1&0&0&0&0&0&0&0&0&0&0&0&0&0&0&0&0\\
0&0&-1&-2&0&0&0&0&-2&-1&0&0&0&0&0&0&0&0&0&0&0&0&0&0&0\\
0&0&0&-1&-1&0&0&0&0&-1&0&0&0&0&0&0&0&0&0&0&0&0&0&0&0 \\ \hline
1&0&0&0&0&0&1&0&0&0&-1&-1&0&0&0&0&0&0&0&0&0&0&0&0&0\\
1&2&0&0&0&-1&0&1&0&0&0&-2&-1&0&0&0&0&0&0&0&0&0&0&0&0\\
0&1&2&0&0&0&-1&0&1&0&0&0&-2&-1&0&0&0&0&0&0&0&0&0&0&0\\
0&0&1&2&0&0&0&-1&0&1&0&0&0&-2&-1&0&0&0&0&0&0&0&0&0&0\\
0&0&0&1&1&0&0&0&-1&0&0&0&0&0&-1&0&0&0&0&0&0&0&0&0&0\\  \hline
0&0&0&0&0&1&0&0&0&0&0&1&0&0&0&-1&-1&0&0&0&0&0&0&0&0\\
0&0&0&0&0&1&2&0&0&0&-1&0&1&0&0&0&-2&-1&0&0&0&0&0&0&0\\
0&0&0&0&0&0&1&2&0&0&0&-1&0&1&0&0&0&-2&-1&0&0&0&0&0&0\\
0&0&0&0&0&0&0&1&2&0&0&0&-1&0&1&0&0&0&-2&-1&0&0&0&0&0\\
0&0&0&0&0&0&0&0&1&1&0&0&0&-1&0&0&0&0&0&-1&0&0&0&0&0\\  \hline
0&0&0&0&0&0&0&0&0&0&1&0&0&0&0&0&1&0&0&0&-1&-1&0&0&0\\
0&0&0&0&0&0&0&0&0&0&1&2&0&0&0&-1&0&1&0&0&0&-2&-1&0&0\\
0&0&0&0&0&0&0&0&0&0&0&1&2&0&0&0&-1&0&1&0&0&0&-2&-1&0\\
0&0&0&0&0&0&0&0&0&0&0&0&1&2&0&0&0&-1&0&1&0&0&&-2&-1\\
0&0&0&0&0&0&0&0&0&0&0&0&0&1&1&0&0&0&-1&0&0&0&0&0&-1\\ \hline 
0&0&0&0&0&0&0&0&0&0&0&0&0&0&0&1&0&0&0&0&1&1&0&0&0\\
0&0&0&0&0&0&0&0&0&0&0&0&0&0&0&1&2&0&0&0&0&2&1&0&0\\
0&0&0&0&0&0&0&0&0&0&0&0&0&0&0&0&1&2&0&0&0&0&2&1&0\\
0&0&0&0&0&0&0&0&0&0&0&0&0&0&0&0&0&1&2&0&0&0&0&2&1\\
0&0&0&0&0&0&0&0&0&0&0&0&0&0&0&0&0&0&1&1&0&0&0&0&1\\
\end{array} \right]
$

\noindent In general, the $R$ matrix is an $n^2 \times n^2$ block tridiagonal matrix with at most 6 nonzero entries per row and 6 per column.
$$R =\dfrac{h}{6} \hat{k}\left[ \begin{array}{cccccc} A_{1,1} & A_{1,2}& 0&0&\cdots&0\\
A_{2,1} & A_{2,2}& A_{2,3}&0&\cdots&0\\
0&\ddots&\ddots&\ddots&\ddots&\vdots\\
\vdots&\ddots&A_{j,k}&A_{j,j}&A_{j,i}&0\\
0&\cdots&0&A_{i,j}&A_{i,i}&A_{i,n}\\
0&\cdots&0&0&A_{n,i}&A_{n,n}\\
\end{array}\right]
$$
where $i = n-1, j = n-2, k= n-3$, and the block matrices are of size $n\times n$ with the following sparsity patterns:
\begin{itemize}
\item $A_{1,1}$ is a lower bidiagonal matrix with $2n-1$ nonzero entries.
 \item $A_{n,n} = - A_{1,1}^T$.
\item $A_{i,i}$ for $i = 2,3,.., n-1$ are tridiagonal matrices with zero diagonal entries with $A(i,i)^T = -A(i,i)$.
\item $A_{i+1,i} = -A_{1,1}$ for $i = 1,2,3,.., n-1$
\item $A_{i,i+1}=A_{1,1}^T$ for $i = 1,2,3,.., n-1$ 
\end{itemize}
Thus, $R$ has a total of $(n-2)(2n-2)+(2n)(2n-1) = 2n^2-6n+4+4n^2-2n = 6n^2 - 8n + 4$ nonzero entries.
 
 \noindent As for the nonzero entries in $R$, their absolute values are less than or equal to $\dfrac{h}{3}\;|\hat{k}|$.
%
%
\subsubsection{The Matrix $R$ assuming Periodic Boundary Conditions}\label{sec:PR}
Given the meshing of Figure \ref{fig:mesh2} and assuming periodic boundary conditions, where the values of U are equal at the upper and lower red vertices, the left and right red boundary vertices, and the corner vertices i.e. 
\begin{eqnarray}
U_{kn+1}&=&U_{kn+n}, \qquad \qquad \qquad \qquad \quad \;\; for \;\; k = 1,2,..,n-2\\
U_{i}&=&U_{i+n(n-1)}, \qquad \qquad \qquad \qquad for \;\; i = 2, 3,..,n-1\\ 
U_1 &=& U_n \; = \; U_{1+n(n-1)} \; = \; U_{n^2}
\end{eqnarray} and assuming that our domain is torus shaped and that the vector U is of size $(n-1)^2$, $$U = [U_1,.., U_{n-1},U_{n+1},...,U_{2n-1}, U_{2n+1},...,U_{3n-1}, ......, U_{n(n-1)-1}]$$ 
then $R$ is an $(n-1)^2\times (n-1)^2$ sparse matrix with at most 6 nonzero entries per row, as discussed below. This ``periodic" $R$ matrix, can be obtained from the general one described in the previous section by merging/adding the rows corresponding to equal $U$ entries and also the columns.

 Note that since rows/columns $n, 2n, 3n , ..., (n-1)n$ of the general $R$ matrix are merged with other rows/columns, then the indices have to be reindexed to get the corresponding rows/columns of the ``periodic" matrix $R$ as such:
\begin{eqnarray}
\left[ 1, \, \cdots , \, {n-1}, \, {n+1}, \, \cdots, \, {2n-1}, \right.&{2n+1},&\left. \cdots, \, {3n-1}, \, \cdots \cdots,\, {(n-1)n-1} \right] \nonumber\\
& \downarrow &  \label{reindex2}\\
\left[ 1, \, \cdots , \, {n-1}, \, {n}, \, \cdots, \, {2n-2}, \right.&{2n-1},&\left. \cdots, \, {3n-3}, \, \cdots \cdots,\, {(n-1)^2} \right] \nonumber
\end{eqnarray}
The entries in the matrix $R$ that will be modified are the ones corresponding to the green left corner, lower and left red boundary nodes, and the upper and left boundary black nodes in Figure \ref{fig:mesh2}.
\vspace{2mm}

\noindent \textbf{1- Green Left Node:}
Assuming that the vertices $1$, $n$, $n(n-1) + 1$ and $n^2$ coincide, then the first row of the ``periodic" $R$ will be the sum of the entries in rows  $1$, $n$, $n(n-1) + 1$ and $n^2$  of the general $R$:
\begin{enumerate}
\item[\bf a)] Row $v = 1$ of the general matrix $R$ has the entry $ -\dfrac{h}{6}\;\hat{k}$  in columns $1, n+1, n+2$. 
\item[\bf b)] Row $v = n$ of the general matrix $R$ has the entry $-\dfrac{h}{6}\;\hat{k}$ in columns  $n-1,v = n = 1$, and $2n = n+1$.    
\item[\bf c)] Row $v = n(n-1) + 1$ of the general matrix $R$ has the entry has the entry $\dfrac{h}{6}\;\hat{k}$ in columns $n(n-2)+1, v = n(n-1)+1 = 1$, and $n(n-1)+2 = 2$.  
\item[\bf d)] Row $v = n^2$ of the general matrix $R$ has the entry $ \dfrac{h}{6}\;\hat{k}$  in columns $n(n-1)-1,
 n(n-1) = n(n-2)+1,  v=n^2=1$.
\end{enumerate}
\noindent Adding up these rows, we get row $1$ of the ``periodic" $R$ with the following nonzero entries:\vspace{1mm}

$\begin{array}{lll}
    1: &  -\dfrac{h}{6}\;\hat{k} -\dfrac{h}{6}\;\hat{k}+\dfrac{h}{6}\;\hat{k}+\dfrac{h}{6}\;\hat{k}&=0\nonumber \vspace{3mm}    \end{array}$ 
    
\begin{minipage}{0.2\textwidth}
 $\begin{array}{ll}
    2: &\dfrac{h}{6}\;\hat{k} \nonumber \vspace{3mm}\\
         n-1: &-\dfrac{h}{6}\;\hat{k} \nonumber\vspace{3mm}
    \end{array}$ 
  \end{minipage}
  \begin{minipage}{0.35\textwidth}
 $\begin{array}{lll}
    n+1: &-\dfrac{h}{6}\;\hat{k}-\dfrac{h}{6}\;\hat{k}&=-\dfrac{h}{3}\;\hat{k}\nonumber \vspace{3mm}\\
        n+2: &-\dfrac{h}{6}\;\hat{k} &\nonumber \vspace{3mm}
  \end{array}$ 
  \end{minipage}
   \begin{minipage}{0.3\textwidth}
 $\begin{array}{lll}
    n(n-2)+1: &\dfrac{h}{6}\;\hat{k}+\dfrac{h}{6}\;\hat{k}&=\dfrac{h}{3}\;\hat{k}\nonumber \vspace{3mm}\\
    n(n-1)-1: & \dfrac{h}{6}\;\hat{k}&\nonumber\vspace{3mm}
  \end{array}$ 
  \end{minipage}
  
  Recall that these column indices have to be reindex by \eqref{reindex2}.\vspace{2mm}
  
  \noindent \textbf{2- Lower Red Boundary Nodes:}
  Assuming that the vertices $i$ and ${i+n(n-1)}$   coincide for  $i = 2, 3,..,n-1$, then the rows $i$ will be the sum of the entries in rows $i$ and ${i+n(n-1)}$ of the general $S(U)$:
  \begin{enumerate}
\item[\bf a)] Row $i$ of the general matrix $R$ for $i = 2,3,...,n-1$ has nonzero entries in columns  \vspace{2mm}\\
\begin{minipage}{0.3\textwidth}
$\begin{array}{ll}
 i-1 : &-\dfrac{h}{6}\;\hat{k}    \nonumber \vspace{3mm}\\
 i+n : &-\dfrac{h}{3}\;\hat{k} \nonumber \vspace{3mm}\\
\end{array}$
\end{minipage}
\begin{minipage}{0.3\textwidth}
$\begin{array}{ll}
v = i : & -\dfrac{h}{3}\;\hat{k}\nonumber \vspace{3mm}\\
   i+n+1 : &-\dfrac{h}{6}\;\hat{k} \nonumber \vspace{3mm}\\
\end{array}$
 \end{minipage}

\item[\bf b)] Row $v = n(n-1)+i$ of the general matrix $R$ for $i = 2,3,...,n-1$, has nonzero entries in columns \\
\begin{minipage}{0.4\textwidth}
$\begin{array}{ll}
 n(n-2)+i-1 : &\dfrac{h}{6}\;\hat{k}    \nonumber \vspace{3mm}\\
 n(n-2)+i : & \dfrac{h}{3}\;\hat{k}   \nonumber \vspace{3mm}\\
\end{array}$
 \end{minipage}
\begin{minipage}{0.5\textwidth}
$\begin{array}{ll}
 v =  n(n-1)+i=i : & \dfrac{h}{3}\;\hat{k} \nonumber \vspace{3mm}\\
 n(n-1)+ i+1=i+1 : &  \dfrac{h}{6}\;\hat{k}\nonumber \vspace{3mm}\\
\end{array}$
 \end{minipage}
\end{enumerate}
Adding up these 2 rows, we get row $i$ of the ``periodic" $R$ with the following nonzero entries for $i = 2,3,...,n-1$:\\
\begin{minipage}{0.25\textwidth}
 $\begin{array}{ll}
 i-1:& -\dfrac{h}{6}\;\hat{k}   \nonumber \vspace{3mm}\\
 i+1:& \dfrac{h}{6}\;\hat{k} \nonumber \vspace{3mm}\\
 \end{array}$
 \end{minipage}
\begin{minipage}{0.3\textwidth}
 $\begin{array}{ll}
 i+n:& -\dfrac{h}{3}\;\hat{k}  \nonumber \vspace{3mm}\\
 i+n+1:&-\dfrac{h}{6}\;\hat{k}  \nonumber \vspace{3mm}\\
 \end{array}$
 \end{minipage}
 \begin{minipage}{0.3\textwidth}
 $\begin{array}{ll}
 n(n-2)+i-1 :&\dfrac{h}{6}\;\hat{k}   \nonumber \vspace{3mm} \\
 n(n-2)+i:&\dfrac{h}{3}\;\hat{k}  \nonumber \vspace{3mm}\\
 \end{array}$
 \end{minipage}
  $\begin{array}{lll}i:&\dfrac{h}{3}\;\hat{k} -\dfrac{h}{3}\;\hat{k} & =0 \nonumber \vspace{3mm}\end{array}$

 \noindent Recall that these column indices have to be reindex by \eqref{reindex2}.
 Note that for $i=n-1$ we get the following nonzero entries:\\
 
 \begin{minipage}{0.25\textwidth}
 $\begin{array}{ll}
 n-2:& -\dfrac{h}{6}\;\hat{k}   \nonumber \vspace{3mm}\\
 n=1:& \dfrac{h}{6}\;\hat{k} \nonumber \vspace{3mm}\\
 \end{array}$
 \end{minipage}
\begin{minipage}{0.3\textwidth}
 $\begin{array}{ll}
 2n-1:& -\dfrac{h}{3}\;\hat{k}  \nonumber \vspace{3mm}\\
 2n=n+1:&-\dfrac{h}{6}\;\hat{k}  \nonumber \vspace{3mm}\\
 \end{array}$
 \end{minipage}
 \begin{minipage}{0.3\textwidth}
 $\begin{array}{ll}
 n(n-2)+n-2 :&\dfrac{h}{6}\;\hat{k}   \nonumber \vspace{3mm} \\
 n(n-2)+n-1:&\dfrac{h}{3}\;\hat{k}  \nonumber \vspace{3mm}\\
 \end{array}$
 \end{minipage}

  \noindent \textbf{3- Left Red Boundary Nodes:}
  Assuming that the vertices $ln+1$ and $ln+n$ coincide for  $l = 1, 2,..,n-2$, then the corresponding rows $ln+1-l$ of the ``periodic" $R$  will be the sum of the entries in rows $ln+1$ and $ln+n$ of the general $R$:
  \begin{enumerate}
\item[\bf a)] Row $ln+1$ of the general matrix $R$ for $l = 1, 2,..,n-2$,  has nonzero entries in columns  \\
\begin{minipage}{0.45\textwidth}
$\begin{array}{ll}
 nl-n+1 : & \dfrac{h}{6}\;\hat{k} \nonumber \vspace{3mm}\\
  nl+2 : &\dfrac{h}{6}\;\hat{k} \nonumber \vspace{3mm}\\
\end{array}$
\end{minipage}
\begin{minipage}{0.45\textwidth}
$\begin{array}{ll}
  nl+n+1 : &  -\dfrac{h}{6}\;\hat{k}\nonumber \vspace{3mm}\\
   nl+n+2 : & -\dfrac{h}{6}\;\hat{k}\nonumber \vspace{-1mm}
\end{array}$
\end{minipage}

\item[\bf b)] Row $(l+1)n$ of the general matrix $R$ for $l = 1,2,...,n-2$, has nonzero entries in columns \\
\begin{minipage}{0.45\textwidth}
$\begin{array}{ll}
 nl-1 : &\dfrac{h}{6}\;\hat{k} \nonumber \vspace{3mm}\\
 nl = nl-n+1 : &\dfrac{h}{6}\;\hat{k} \nonumber \vspace{3mm}\\
\end{array}$
\end{minipage}
\begin{minipage}{0.45\textwidth}
$\begin{array}{ll}
 nl+n-1 : &-\dfrac{h}{6}\;\hat{k} \nonumber \vspace{3mm}\\
  nl+n+n = nl+n+1 : &-\dfrac{h}{6}\;\hat{k}\nonumber
\end{array}$
\end{minipage}

\end{enumerate}
Adding up these 2 rows, we get row $nl+1-l$ of the ``periodic" $R$ with the following nonzero entries for $l = 1,2,...,n-2$:\\
\begin{minipage}{0.45\textwidth}
 $\begin{array}{lll}
 nl-n+1 :&\dfrac{h}{6}\;\hat{k} +\dfrac{h}{6}\;\hat{k}  & = \dfrac{h}{3}\;\hat{k} \vspace{3mm}\\
 nl+n+1:& -\dfrac{h}{6}\;\hat{k}-\dfrac{h}{6}\;\hat{k}  & = -\dfrac{h}{3}\;\hat{k} \vspace{3mm}\\
 \end{array}$
\end{minipage}
\begin{minipage}{0.3\textwidth}
 $\begin{array}{ll}
 nl-1:&\dfrac{h}{6}\;\hat{k}  \vspace{3mm}\\
 nl+2:& \dfrac{h}{6}\;\hat{k}    \vspace{3mm}\\
 \end{array}$
 \end{minipage}
 \begin{minipage}{0.3\textwidth}
 $\begin{array}{ll}
 nl+n-1:&-\dfrac{h}{6}\;\hat{k}   \vspace{3mm}\\
 nl+n+2:&-\dfrac{h}{6}\;\hat{k}    \vspace{3mm}\\
 \end{array}$
 \end{minipage}

 \noindent Recall that these column indices have to be reindex by \eqref{reindex2}.
 Note that for $l=n-2$ we get the following nonzero entries:\\
 \begin{minipage}{0.45\textwidth}
 $\begin{array}{lll}
 n(n-2)-n+1 :& \dfrac{h}{3}\;\hat{k} \vspace{3mm}\\
 n(n-2)+n+1 = 1:&-\dfrac{h}{3}\;\hat{k} \vspace{3mm}
 \end{array}$
\end{minipage}
\begin{minipage}{0.3\textwidth}
 $\begin{array}{ll}
 n(n-2)-1:&\dfrac{h}{6}\;\hat{k}  \vspace{3mm}\\
 n(n-2)+2:& \dfrac{h}{6}\;\hat{k}    \vspace{3mm}
 \end{array}$
 \end{minipage}
 \begin{minipage}{0.3\textwidth}
 $\begin{array}{ll}
 n(n-2)+n-1:&-\dfrac{h}{6}\;\hat{k}   \vspace{3mm}\\
 n(n-2)+n+2=2:&-\dfrac{h}{6}\;\hat{k}    \vspace{3mm}
 \end{array}$
 \end{minipage}
%
 
  \noindent \textbf{4- Right Black Boundary Nodes:}
  Rows $(l+1)n-l-1$  of the ``periodic" $R$ matrix correspond to rows $(l+1)n-1$ of the general $R$ for  $l = 1,..,n-3$ with nonzero entries  in columns:\\  
   \begin{minipage}{0.25\textwidth}
$\begin{array}{ll}
 nl-2 : &  \dfrac{h}{6}\hat{k}\nonumber \vspace{3mm}\\
nl-1 : &\dfrac{h}{3}\hat{k}  \nonumber \vspace{3mm}
\end{array}$
  \end{minipage}
 \begin{minipage}{0.3\textwidth}
$\begin{array}{ll}
 nl+n-2 : & -\dfrac{h}{6}\hat{k} \nonumber \vspace{3mm}\\
  nl+2n-1 : &  -\dfrac{h}{3}\hat{k} \nonumber \vspace{3mm}
\end{array}$
  \end{minipage}
  \begin{minipage}{0.3\textwidth}
$\begin{array}{ll}
  nl+n = nl+1 : & \dfrac{h}{6}\hat{k}\nonumber \vspace{3mm}\\
    nl+2n = nl+n+1: & -\dfrac{h}{6}\hat{k}\nonumber
\end{array}$
  \end{minipage}\\
Recall that these column indices have to be reindex by \eqref{reindex2}.\vspace{2mm}

   \noindent \textbf{5- Upper Black Boundary Nodes:}
  Rows $(n-1)(n-2)+i$  of the ``periodic" $R$ matrix correspond to rows $n(n-2)+i$ of the general $R$ for  $i = 2,3,..,n-1$ with nonzero entries in columns:\\ 
   \begin{minipage}{0.3\textwidth}
  $\begin{array}{ll}
n(n-3)+i-1   : & \dfrac{h}{6}\hat{k}\nonumber \vspace{3mm}\\
 n(n-3)+i : & \dfrac{h}{3}\hat{k} \nonumber \vspace{3mm}
\end{array}$
   \end{minipage}
  \begin{minipage}{0.3\textwidth}
  $\begin{array}{ll}
 n(n-2)+i-1 : &-\dfrac{h}{6}\hat{k}\nonumber \vspace{3mm}\\ 
  n(n-2)+i+1 : & \dfrac{h}{6}\hat{k} \nonumber \vspace{3mm}
\end{array}$
   \end{minipage}
  \begin{minipage}{0.3\textwidth}
  $\begin{array}{ll}
  n(n-1)+i = i : &-\dfrac{h}{3}\hat{k}     \nonumber \vspace{3mm}\\
    n(n-1)+i+1 = i+1 : & -\dfrac{h}{6}\hat{k}   \nonumber
\end{array}$
   \end{minipage}

  Recall that these column indices have to be reindex by \eqref{reindex2}.
 Note that for $i=n-1$ we get the following nonzero entries:\\
  \begin{minipage}{0.3\textwidth}
 $\begin{array}{ll}
  n(n-2)-2: &\dfrac{h}{6}\hat{k} \nonumber \vspace{3mm}\\
 n(n-2)-1: & \dfrac{h}{3}\hat{k}  \nonumber \vspace{3mm}
\end{array}$
  \end{minipage}
  \begin{minipage}{0.3\textwidth}
   $\begin{array}{ll}
 n(n-1)-2 : & -\dfrac{h}{6}\hat{k} \nonumber \vspace{3mm}\\
    n-1: &  -\dfrac{h}{3}\hat{k}   \nonumber \vspace{3mm}
\end{array}$
  \end{minipage}
  \begin{minipage}{0.3\textwidth}
   $\begin{array}{ll}
   n(n-1) = n(n-2)+1  : &  \dfrac{h}{6}\hat{k} \nonumber \vspace{3mm}\\
  n-1+1 = n = 1 : & -\dfrac{h}{6}\hat{k}  \nonumber
\end{array}$
  \end{minipage}
  
 \noindent For $n=5$, the matrix $R$ corresponding to the mesh in Figure \ref{fig:mesh2} has the following sparsity pattern 
  $$R = \dfrac{h}{6}\;\hat{k}\;\left[
\begin{array}{cccc|cccc|cccc|cccc}
0&1&0&-1&-2&-1&0&0&0&0&0&0&2&0&0&1\\
-1&0&1&0&0&-2&-1&0&0&0&0&0&1&2&0&0\\
0&-1&0&1&0&0&-2&-1&0&0&0&0&0&1&2&0\\
1&0&-1&0&-1&0&0&-1&0&0&0&0&0&0&1&2\\
\hline
2&0&0&1&0&1&0&-1&-2&-1&0&0&0&0&0&0\\
1&2&0&0&-1&0&1&0&0&-2&-1&0&0&0&0&0\\
0&1&2&0&0&-1&0&1&0&0&-2&-1&0&0&0&0\\
0&0&1&2&1&0&-1&0&-1&0&0&-2&0&0&0&0\\
 \hline
0&0&0&0&2&0&0&1&0&1&0&-1&-2&-1&0&0\\
0&0&0&0&1&2&0&0&-1&0&1&0&0&-2&-1&0\\
0&0&0&0&0&1&2&0&0&-1&0&1&0&0&-1&-1\\
0&0&0&0&0&0&1&2&1&0&-1&0&-1&0&0&-2\\
\hline 
-2&-1&0&0&0&0&0&0&2&0&0&1&0&1&0&-1\\
0&-2&-1&0&0&0&0&0&1&2&0&0&-1&0&1&0\\
0&0&-2&-1&0&0&0&0&0&1&2&0&0&-1&0&1\\
-1&0&0&-2&0&0&0&0&0&0&1& 2&1&0&-1&0\\
\end{array} \right]
$$ 
 \noindent In general, the $R$ matrix is an $(n-1)^2 \times (n-1)^2$ block tridiagonal matrix with 2 additional blocks in the upper right  and lower left corner. Moreover, it is a skew-symmetric matrix ($R^T = -R$) with zeros on the diagonal, and 6 nonzero entries per row of the form $\dfrac{h}{6}\hat{k} \alpha$, 6 nonzero entries per column, of the form $\dfrac{h}{6}\hat{k} \alpha$ where $\alpha = -2,-1,-1,1,1,2$.
 
$$R =\dfrac{h}{6}\hat{k} \left[ \begin{array}{cccccc} A_{1,1} & A_{1,2}& 0&\cdots&0&A_{1,n-1}\\
A_{2,1} & A_{2,2}& A_{2,3}&0&\cdots&0\\
0&\ddots&\ddots&\ddots&\ddots&\vdots\\
\vdots&\ddots&A_{j,m}&A_{j,j}&A_{j,i}&0\\
0&\cdots&0&A_{i,j}&A_{i,i}&A_{i,l}\\
A_{l,1}&0&\cdots&0&A_{l,i}&A_{l,l}\\
\end{array}\right]
$$
where $i = n-2, j = n-3, l= n-1, m=n-4$, and the $3(n-1)$ nonzero block matrices are of size $(n-1)\times (n-1)$ with $2(n-1)$ nonzero entries each, and the following sparsity patterns:
\begin{itemize}
\item $A_{i,i}$ for $i = 1,2,.., n-1$ are such that, $A_{i,i}(j,j+1)=1$, $A_{i,i}(j+1,j)=-1$, for $j=1,2,..., n-2$ $A_{i,i}(1,n-1) = -1$, and $A_{i,i}(n-1,1) = 1$.
\item $A_{1,n-1} = A_{i+1,i}$ for $i = 1,2,.., n-2$ are lower bidiagonal matrices ($A_{i+1,i}(j,j) =2, A_{i+1,i}(j+1,j) =1$), with $A_{i+1,i}(1,n-1) = 1$,
\item  $A_{n-1,1} = A_{i,i+1}$ for $i = 1,.., n-2$ are upper bidiagonal matrices ($A_{i,i+1}(j,j) =-2, A_{i,i+1}(j,j+1) =-1$), with $A_{i,i+1}(n-1,1) = -1$,
\end{itemize}
Thus, $R$ has a total of $3(n-1)2(n-1) = 6(n-1)^2$ nonzero entries.
 As for the nonzero entries in $R$, their absolute values are less than or equal to $\dfrac{h}{3}\;|\hat{k}|$, similary to the general matrix $R$. 

\end{document}